\documentclass[smallextended,natbib,runningheads,nospthms]{svjour3}
\journalname{Annals of the Institute of Statistical Mathematics}
\smartqed  
\usepackage[dvips]{graphicx}
\usepackage{amsthm}
\usepackage[body={6.3in,9in}]{geometry}
\usepackage{amsmath}
\usepackage{amsfonts}
\usepackage{amssymb}
\usepackage{amscd}
\usepackage{mathrsfs}
\usepackage[all]{xypic}
\usepackage{url}
\usepackage{epsfig}
\usepackage{longtable}

\newtheorem{thm}{Theorem}[section]
\newtheorem{cor}[thm]{Corollary}
\newtheorem{lem}[thm]{Lemma}
\newtheorem{prop}[thm]{Proposition}
\newtheorem{conj}[thm]{Conjecture}

\theoremstyle{definition}

\theoremstyle{definition}
\newtheorem{example}[thm]{Example}

\theoremstyle{remark}

\newcommand{\C}{\mathbb{C}}
\newcommand{\N}{\mathbb{N}}

\newcommand{\R}{\mathbb{R}}

\newcommand{\s}{\mathbb{S}}
\newcommand{\PP}{\mathbb{P}}

\begin{document}

\title{Multivariate Gaussians, Semidefinite Matrix Completion, \\
 and Convex Algebraic Geometry}



\author{Bernd Sturmfels         \and
        Caroline Uhler 
\thanks{BS is supported in part by NSF grants DMS-0456960 and DMS-0757236.
CU is supported by an International Fulbright Science
and Technology Fellowship.}}


\institute{Bernd Sturmfels \at
              Department of Mathematics \#3840 \\
              University of California at Berkeley \\
              970 Evans Hall \\
              Berkeley, CA 94720-3840, U.S.A. \\
              \email{bernd@math.Berkeley.edu}         
           \and
           Caroline Uhler \at
              Department of Statistics \#3860 \\
              University of California at Berkeley \\
              345 Evans Hall \\
              Berkeley, CA 94720-3860, U.S.A. \\
              \email{cuhler@stat.Berkeley.edu}
}

\date{ }

\maketitle

\begin{abstract}

We study multivariate normal models that are described by linear constraints
on the inverse of the covariance matrix. Maximum likelihood estimation
for such models leads to the problem of maximizing the determinant
function over a spectrahedron, and to the problem of
characterizing the image of the positive definite cone under
an arbitrary linear projection.
These problems at the interface of statistics and optimization
are here examined
from the perspective of convex algebraic geometry.

\end{abstract}

\section{Introduction}
\label{intro}
Every positive definite $m \times m$-matrix $\Sigma $ is the {\em covariance matrix} of a multivariate
normal distribution on $\R^m$. Its inverse matrix $K = \Sigma^{-1}$ is also positive definite and
known as the {\em concentration matrix} of the distribution.
We study statistical models for multivariate normal distributions on $\R^m$,
 where the concentration matrix can be written as a linear combination
\begin{equation}
\label{con_model}
K \,\, = \,\, \lambda_1 K_1 + \lambda_2 K_2 + \cdots +\lambda_d K_d
\end{equation}
of some fixed linearly independent symmetric matrices $K_1,\ldots ,K_d$.
Here, $\lambda_1,\lambda_2, \ldots ,\lambda_d$ are unknown real coefficients.
It is assumed that $K$ is positive definite for some choice
of $\lambda_1,\lambda_2, \ldots ,\lambda_d$.
Such statistical models, which we call {\em linear concentration models},
 were introduced by \cite{A}.

Let $\s^m $ denote the vector space of real symmetric $m \times m$-matrices.
We identify $\s^m$ with its dual space via the
 inner product  $\,\langle A, B \rangle  := {\rm trace}(A \cdot B)$.
The cone $\s^m_{\succeq 0}$ of positive semidefinite matrices is
a full-dimensional self-dual cone in $\s^m$. Its interior is the open cone
$\s^m_{\succ 0}$ of positive definite matrices.
We define a linear concentration model to be any non-empty set of the form
$$\,\mathcal{L}^{-1}_{\succ 0} \,\,:= \,
\,\bigl\{\,\Sigma \in \s^m_{\succ 0} \,:\,\Sigma^{-1} \in \mathcal{L} \,\bigr\},$$
where $\mathcal{L}$ is a linear subspace of $\s^m$.
Given a basis $K_1, \dots ,K_d$ of the subspace $\mathcal{L}$
as in (\ref{con_model}),
the basic statistical problem is to estimate the parameters
$\lambda_1,\ldots ,\lambda_d$ when $n$ observations $X_1,\dots ,X_n$ are
drawn from a multivariate normal distribution $\mathcal{N}(\mu,\Sigma)$,
whose covariance matrix $\Sigma = K^{-1}$ is
in the model $ \mathcal{L}^{-1}_{\succ 0}$.
The $n$ observations $X_i$
and their mean $\bar{X}$ are summarized in the sample covariance matrix
\[ S \,\,= \,\, \frac{1}{n}\sum_{i=1}^{n}(X_i-\bar{X}) (X_i-\bar{X})^T
 \,\,\, \in \,\, \s^m_{\succeq 0}.
\]
In our model, we make no assumptions on the mean vector $\mu$
and always use the sample mean $\bar{X}$ as estimate for $\mu$.
Thus, we are precisely in the situation of
\cite[Prop.~2.1.12]{DSS}, with $\Theta_2 = \mathcal{L}^{-1}$.
The log-likelihood function for the linear concentration model (\ref{con_model}) equals
\begin{equation} \label{loglike}
 {\rm log} \, {\rm det} (K) - \langle S, K\rangle \,\,\,=\,\,\,
\log\det\biggl(\sum_{j=1}^{d}\lambda_j K_j \biggr)-\sum_{j=1}^{d}\lambda_j \langle S, K_j \rangle
\end{equation}
times the constant $n/2$. This is a strictly concave
function on the relatively open cone
$\,\s^m_{\succ 0} \cap \mathcal{L}$. If a maximum
(i.e. the maximum likelihood estimate or MLE) exists,
then it is attained by a unique matrix $\hat K$ in
$\,\s^m_{\succ 0} \cap \mathcal{L}$.  Its inverse $\hat{\Sigma} = \hat{K}^{-1}$
is uniquely determined by the linear equations
\begin{equation}
\label{Dequations}
\langle \hat{\Sigma}, K_j \rangle \,= \, \langle S, K_j \rangle \qquad
{\rm for} \,\,\,\, j=1,2,\ldots ,d.
\end{equation}
This characterization follows from the statistical theory of {\em exponential families} \cite[\S 5]{B}.
In that theory, the scalars $\lambda_1,\ldots,\lambda_d$ are the
{\em canonical parameters} and
$\langle S, K_1 \rangle ,\ldots, \langle S, K_d \rangle $ are the
{\em  sufficient statistics} of the exponential family (\ref{con_model}).
For a special case see \cite[Theorem~2.1.14]{DSS}.

We consider the set of all covariance matrices whose sufficient statistics are given by the matrix $S$.
This set is a {\em spectrahedron}.
It depends only on $S$ and $\mathcal{L}$, and it is denoted
$$ {\rm fiber}_\mathcal{L}(S) \quad = \quad
\bigl\{ \Sigma \in \s^m_{\succ 0} \,: \,
\langle \Sigma, K \rangle \,  = \,
\langle S, K \rangle \,\,
\hbox{for all} \,\, K \in \mathcal{L} \bigr\}. $$
The MLE exists for a sample covariance matrix $S$ if and only if
${\rm fiber}_\mathcal{L}(S)$ is non-empty.
If ${\rm rank}(S) < m$ then it can happen that the fiber is empty,
in which case the MLE does not exist for $(\mathcal{L},S)$.
Work of \cite{Buh} and \cite{BJT} addresses this issue for
graphical models; see Section \ref{sec:graphical} below.

Our motivating statistical problem is to
identify conditions on the pair $(\mathcal{L},S)$ that ensure
the existence of the MLE.
 This will involve studying the geometry of
the semi-algebraic set $\mathcal{L}^{-1}_{\succ 0}$ and of
the algebraic function
$\, S \mapsto \hat{\Sigma}\,$
which takes a sample covariance matrix to its MLE in
$\mathcal{L}^{-1}_{\succ 0}$.

\begin{example} \label{OldLabsPictures}
We illustrate the concepts introduced so far by way of a small
explicit example whose geometry is visualized in Fig.~1.
Let $m=d=3$ and let $\mathcal{L}$ be the
real vector space spanned by
$$
K_1 \, = \, \begin{pmatrix}
1 & 0 & 0 \\
0 & 1 & 1 \\
0 & 1 & 1
\end{pmatrix} , \quad
K_2 \, = \, \begin{pmatrix}
1 & 0 & 1 \\
0 & 1 & 0 \\
1 & 0 & 1
\end{pmatrix} \quad \hbox{and} \quad
K_3 \, = \, \begin{pmatrix}
1 & 1 & 0 \\
1 & 1 & 0 \\
0 & 0 & 1
\end{pmatrix}.
$$
The linear concentration model (\ref{con_model}) consists of
all positive definite matrices of the form
\begin{equation}
\label{exmodel}
K \,\, = \,\,
\begin{pmatrix}
\lambda_1 + \lambda_2 + \lambda_3 & \lambda_3 & \lambda_2 \\
\lambda_3 & \lambda_1 + \lambda_2 + \lambda_3 & \lambda_1 \\
\lambda_2 & \lambda_1 & \lambda_1 + \lambda_2 + \lambda_3
\end{pmatrix}.
\end{equation}
Given a sample covariance matrix $S = (s_{ij})$, the sufficient statistics are
$$
t_1 = {\rm trace}(S) + 2 s_{23} \,, \,\,\,
t_2 = {\rm trace}(S) + 2 s_{13} \,, \,\,\,
t_3 = {\rm trace}(S) + 2 s_{12} .
$$
If $S \in  \s^3_{\succ 0} $ then ${\rm fiber}_\mathcal{L}(S)$
is an open $3$-dimensional convex body whose boundary
is a cubic surface. This is the {\em spectrahedron}
shown on the left in Fig.~1. The MLE $\hat{\Sigma}$ is the unique
matrix of maximum determinant in the spectrahedron
${\rm fiber}_\mathcal{L}(S)$. Here is an explicit
algebraic formula for
the MLE $\hat{\Sigma}  = (\hat{s}_{ij})$:
First, the matrix entry $\hat{s}_{33}$ is determined (e.g.~using {\em Cardano's formula}\footnote{
{\tt www.literka.addr.com/mathcountry/algebra/quartic.htm}}) from the equation
\begin{eqnarray*}
0 & = & 240\,\hat{s}_{33}^4\, + \, (-32 t_1-32 t_2-192 t_3) \hat{s}_{33}^3
+(-8 t_1^2+16 t_1 t_2+16 t_1 t_3-8 t_2^2+16 t_2 t_3+32 t_3^2) \hat{s}_{33}^2 \\
&&
+(8 t_1^3-8 t_1^2 t_2-8 t_1 t_2^2+8 t_2^3) \hat{s}_{33} \,
- \, 4 t_1^3 t_3-6 t_1^2 t_2^2+4 t_1^2 t_3^2 + 4 t_1 t_2^3
+4 t_1 t_2^2 t_3+4 t_1^2 t_2 t_3-t_2^4
\\ && -4 t_2^3 t_3+4 t_2^2 t_3^2-8
t_1 t_2 t_3^2-t_1^4+4 t_1^3 t_2.
\end{eqnarray*}
Next, we read off $\hat{s}_{23}$ from
\begin{eqnarray*}
-24 \,(t_1^2-2 t_1 t_2+t_2^2-t_3^2) \,\hat{s}_{23} & = &
 120 \hat{s}_{33}^3
- (16 t_1 + 16 t_2+36 t_3)  \hat{s}_{33}^2
         +(2  t_1^2-4 t_1 t_2+2  t_2^2-8  t_3^2) \hat{s}_{33} -6 t_1^3 \\&&
         +18 t_1^2 t_2+t_1^2 t_3-18 t_1 t_2^2-2 t_1 t_2 t_3+10 t_1 t_3^2+
6 t_2^3+t_2^2 t_3-2 t_2 t_3^2-4 t_3^3.
\end{eqnarray*}
Then we read off $\hat{s}_{22}$ from
\begin{eqnarray*}
- 24\,(t_1-t_2) \, \hat{s}_{22} & = &
60 \hat{s}_{33}^2
\,+\,(4 t_1 - 20 t_2 - 24 t_3) \hat{s}_{33}
+ 24 (t_1 - t_2 - t_3) \hat{s}_{23} \\&&
-\,11 t_1^2 + 10 t_1 t_2 + 10 t_1 t_3 + t_2^2- 2 t_2 t_3 - 4 t_3^2.
\end{eqnarray*}
Finally,  we obtain the first row of $\hat{\Sigma}$ as follows:
$$ \hat{s}_{13} \,\,=\,\, \hat{s}_{23}- t_1/2+ t_2/2, \quad\quad
 \hat{s}_{12} \,\,= \,\, \hat{s}_{23}- t_1/2 + t_3/2, \quad\quad
 \hat{s}_{11} \,\,=\,\, t_1 - \hat{s}_{33} - 2 \hat{s}_{23} - \hat{s}_{22} . $$
\begin{center}
\begin{figure*}[bt]
\begin{tabular}{c c c}
\parbox{5.7cm}{\includegraphics[width=5.3cm]{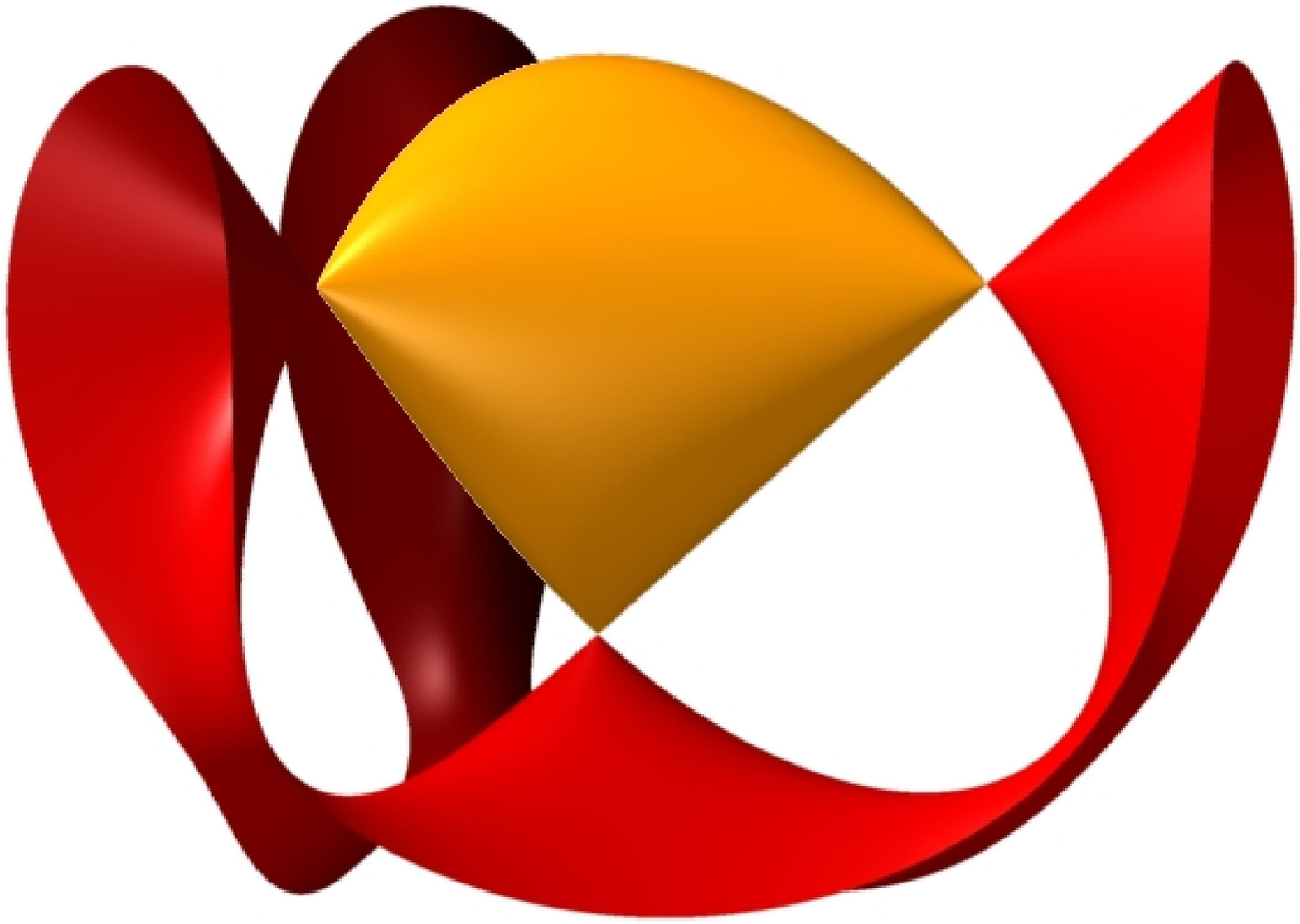}} &
\parbox{4.7cm}{\includegraphics[width=4.2cm]{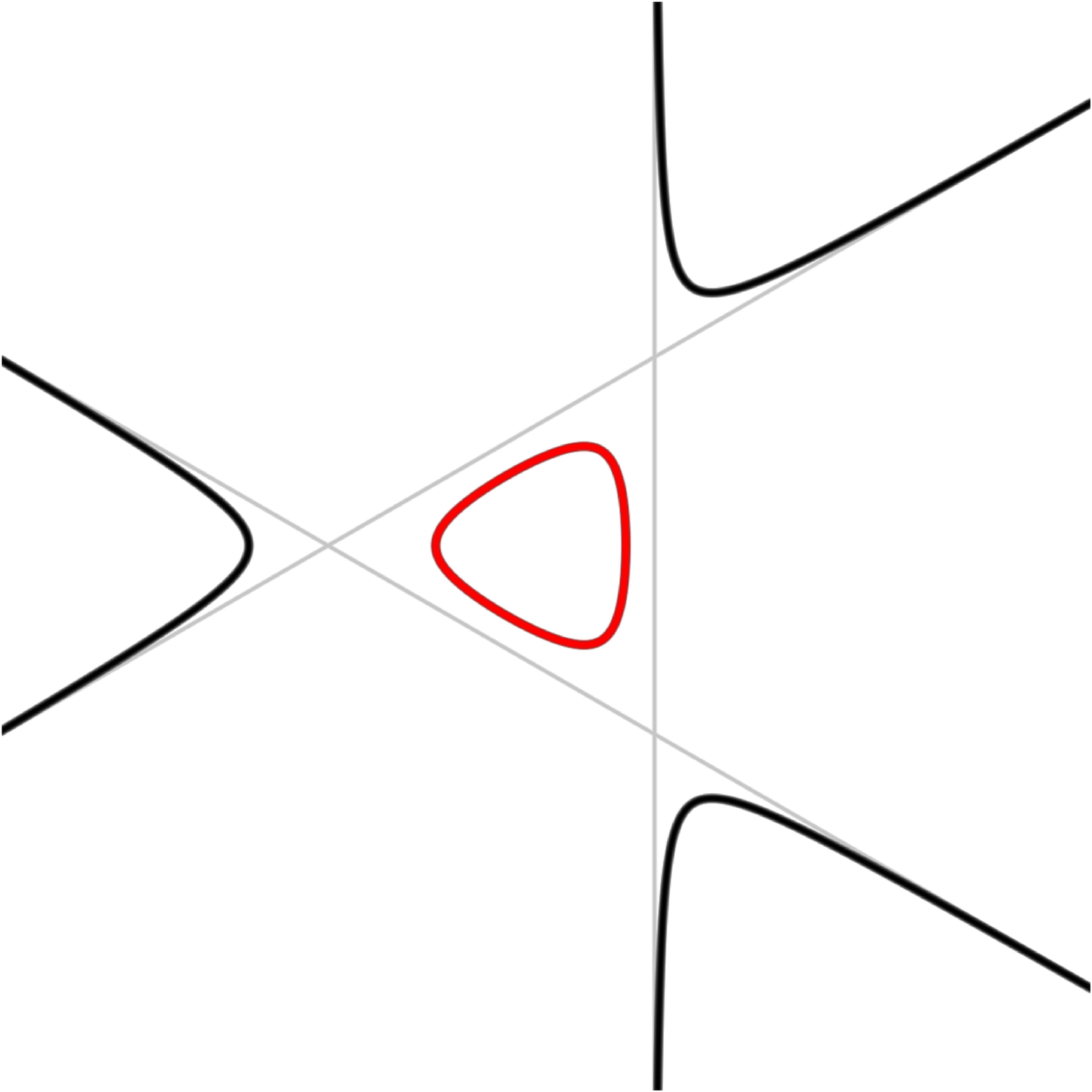}} &
\parbox{4.9cm}{\includegraphics[width=4.8cm]{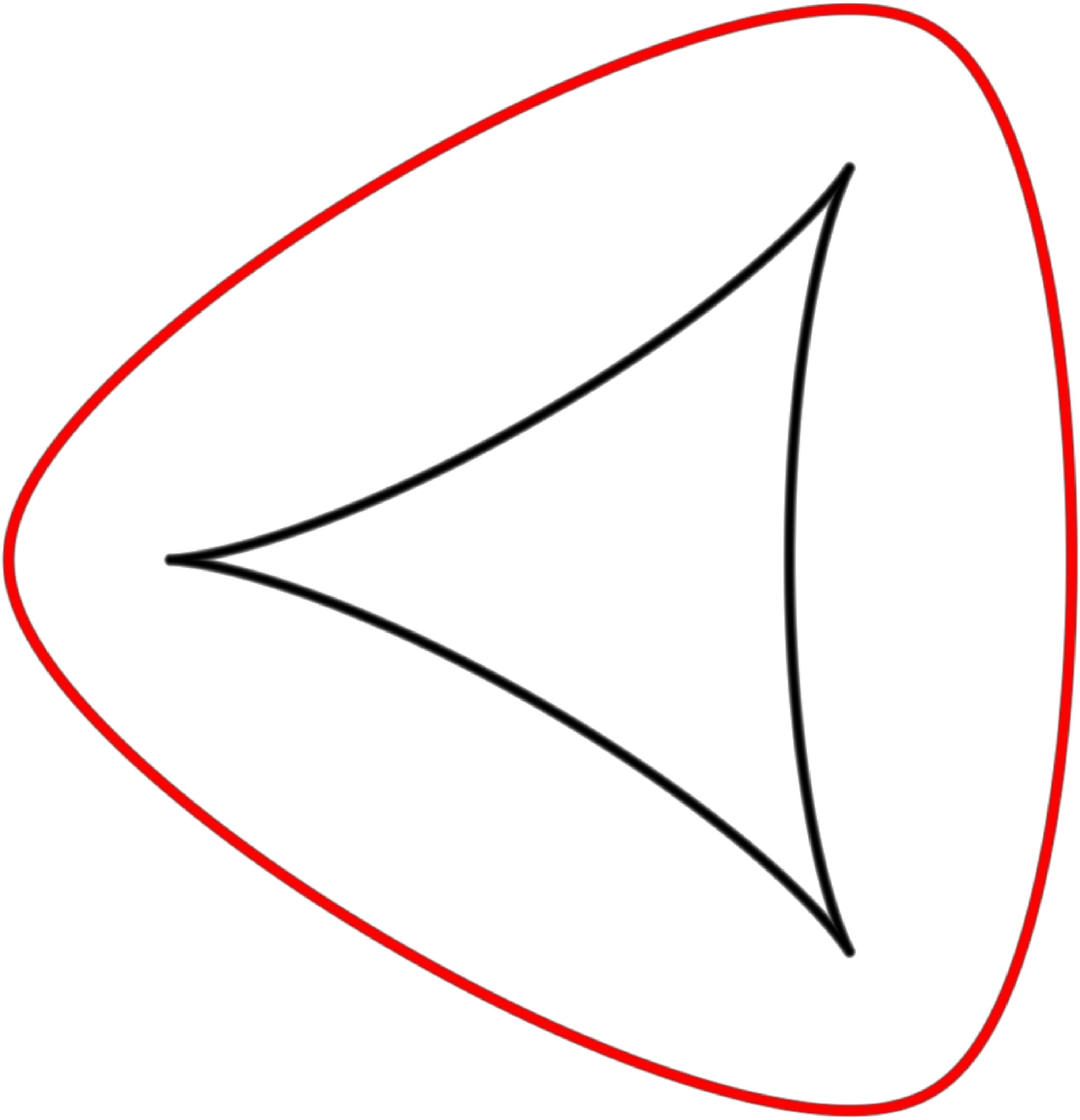}}
\end{tabular}
\caption{Three figures, taken from \cite{NRS}, illustrate
Example \ref{OldLabsPictures}.
These figures show the spectrahedron ${\rm fiber}_\mathcal{L}(S)$ (left), a
cross section of the spectrahedral cone $\mathcal{K}_{\mathcal{L}}$ (middle), and a
cross section of its dual cone $\mathcal{C}_{\mathcal{L}}$ (right).}
\label{fig_labs}
\end{figure*}
\end{center}

\vskip -.3cm

The MLE $\hat{\Sigma} = (\hat{s}_{ij})$ is an algebraic function of
degree $4$ in the sufficient statistics $(t_1,t_2,t_3)$.
In short, the model (\ref{exmodel}) {\em has ML degree~$4$}.
We identify our model with the  subvariety $\mathcal{L}^{-1}$
of projective space $\PP^5$ that is
parametrized by this algebraic function.
 The ideal of polynomials vanishing on  $\mathcal{L}^{-1}$ equals
\begin{eqnarray*}
P_\mathcal{L} &\,\,= \,\, & \langle \,
 s_{13}^2- s_{23}^2- s_{11} s_{33} + s_{22} s_{33}\,,\,\,\,
 s_{12}^2- s_{11} s_{22}- s_{23}^2+ s_{22} s_{33} \,,\,\, \\ && \,\,\,
  s_{12} s_{13} - s_{13} s_{22}- s_{11} s_{23}+ s_{12} s_{23}+ s_{13} s_{23}+ s_{23}^2
  - s_{12} s_{33}- s_{22} s_{33} \, , \\ && \,\,\,
 s_{11} s_{13} {-} s_{13} s_{22} {-} s_{11} s_{23} {+} s_{22} s_{23} {-} s_{11} s_{33}
 {-} 2 s_{12} s_{33} {-} s_{13} s_{33} {-} s_{22} s_{33} {-} s_{23} s_{33} {+} s_{33} ^2, \\ &&\,\,\,
      s_{11} s_{12}{-} s_{11} s_{22} {-} s_{12} s_{22} {-} 2 s_{13} s_{22} {+} s_{22}^2
      {-} s_{11} s_{23} {-} s_{22} s_{23} {-} s_{12}  s_{33} {-} s_{22}  s_{33} {+} s_{23} s_{33}, \\ && \,\,\,
s_{11}^2 - 2 s_{11} s_{22} - 4 s_{13} s_{22} + s_{22}^2 - 4 s_{11} s_{23}
      -2 s_{11} s_{33} - 4 s_{12} s_{33} {-} 2 s_{22} s_{33} {+} s_{33}^2  \,\,\,\rangle.
\end{eqnarray*}
The domain of the maximum likelihood map $\,(t_1,t_2,t_3) \mapsto \hat{\Sigma}\,$ is
the cone of sufficient statistics
$\mathcal{C}_\mathcal{L}$ in $\R^3$. The polynomial
$H_\mathcal{L}$ which vanishes on the boundary of this convex cone
 has degree six. It equals
\begin{eqnarray*}H_\mathcal{L} & \,\,= \,\,&
t_1^6-6 t_1^5 t_2+19 t_1^4 t_2^2-28 t_1^3 t_2^3
+19 t_1^2 t_2^4-6 t_1 t_2^5+t_2^6-6 t_1^5 t_3
  +   14 t_1^4 t_2 t_3  - 24 t_1^3 t_2^2 t_3
-24 t_1^2 t_2^3 t_3 \\&&
     +14 t_1 t_2^4 t_3-6 t_2^5 t_3+19 t_1^4 t_3^2-24 t_1^3 t_2  t_3^2
      +106 t_1^2 t_2^2 t_3^2-24 t_1 t_2^3 t_3^2+19 t_2^4 t_3^2
     -28 t_1^3 t_3^3
-24 t_1^2 t_2 t_3^3\\&&  -24 t_1 t_2^2 t_3^3
     -28 t_2^3 t_3^3+19 t_1^2 t_3^4+14 t_1 t_2 t_3^4
     +19 t_2^2 t_3^4-6 t_1 t_3^5-6 t_2 t_3^5+t_3^6.
\end{eqnarray*}
The sextic curve $\{H_\mathcal{L} = 0\}$  in $\PP^2$ is
shown on the right in Fig.~1.
 It is dual to the cubic curve $\{{\rm det}(K) = 0\}$,
shown in the middle of Fig.~1.
The cone over the convex region enclosed by the
red part of that cubic curve
is the set $\,\mathcal{K}_\mathcal{L} = \s^3_{\succ 0} \cap \mathcal{L}\,$
of concentration matrices in our model (\ref{con_model}).
\qed
\end{example}

This paper is organized as follows. In Section \ref{sec_duality} we
formally define the objects $\mathcal{K}_\mathcal{L}$,
$\mathcal{C}_\mathcal{L}$, $P_\mathcal{L}$ and $H_\mathcal{L}$,
which already appeared in Example \ref{OldLabsPictures}, and
we derive three guiding questions that constitute the
main thread of this paper. These questions are answered
for generic linear spaces $\mathcal{L}$ in Subsection \ref{subsec_generic}.
That subsection is written for algebraists, and readers from statistics
or optimization can skip it at their first reading.
In Section \ref{sec_diagonal} we answer our three questions for
diagonal concentration models, using results from geometric combinatorics.
Section \ref{sec:graphical} deals with Gaussian graphical models,
which are the most prominent linear concentration models.
We resolve our three questions for chordal graphs, then for chordless cycles, and finally for wheels and all graphs with five or less vertices. We conclude this paper with a study of colored Gaussian graphical models in Section \ref{colored_graphs}. These are special Gaussian graphical models with additional linear restrictions on the concentration matrix given by the graph coloring.
\vskip 0.5cm

\section{Linear Sections, Projections and Duality}
\label{sec_duality}

Convex algebraic geometry is concerned with the geometry of
real algebraic varieties
and semi-algebraic sets that arise in convex optimization, especially
in semidefinite programming. A fundamental problem is to
study convex sets that arise as
 linear sections and projections of the
cone of positive definite matrices $\s^m_{\succ 0}$.
 As we saw in the Introduction, this problem arises naturally
 when studying maximum likelihood estimation in
linear concentration models for Gaussian random variables.
In particular, the issue of estimating a covariance matrix
from the sufficient statistics can be seen as an extension of the
familiar semidefinite matrix completion problem
\citep{BJT, GJSW}. In what follows, we develop
an algebraic and geometric framework
for systematically addressing such problems.

\subsection{Derivation of three guiding questions}
\label{subsec_derivation}

As before, we fix a linear subspace
$\mathcal{L}$ in the real vector space $\s^m$ of symmetric $m\times m$-matrices,
and we fix a basis  $\{K_1,\ldots,K_d\}$ of  $\mathcal{L}$.
The {\em cone of concentration matrices} is the relatively open cone
$$ \mathcal{K}_\mathcal{L} \,\,=\,\,
\mathcal{L} \,\cap \,\s^m_{\succ 0}. $$
We assume throughout that $\mathcal{K}_\mathcal{L}$ is non-empty.
Using the basis $K_1,\ldots , K_d$ of $\mathcal{L}$, we can identify $\mathcal{K}_\mathcal{L}$ with
\begin{equation}
\label{K_with_basis}
 \mathcal{K}_\mathcal{L} \,\,=\,\,
\bigl\{ (\lambda_1,\ldots,\lambda_d) \in \R^d \,:\,
\sum_{i=1}^d \lambda_i K_i \,\,\hbox{is positive definite} \,\bigr\}.
\end{equation}
This is a non-empty open convex cone in $\R^d$.
The orthogonal complement $\mathcal{L}^\perp$ of $\mathcal{L}$
is a subspace of dimension $\,\binom{m+1}{2} - d$ in $\s^m$,
so that $\s^m/\mathcal{L}^\perp \,\simeq \, \R^d $,
and we can consider the canonical map
$$ \pi_\mathcal{L} \,: \, \s^m \rightarrow \s^m/{\mathcal L}^\perp . $$
This is precisely the linear map which takes a sample covariance matrix $S$ to its canonical sufficient statistics.
The chosen basis of $\mathcal{L}$ allows us to identify this map with
\begin{equation}
\label{pi_with_basis}
\pi_\mathcal{L} \,: \,\s^m \rightarrow \R^d, \,\, S \mapsto \bigl( \langle S,K_1 \rangle,
\ldots,\langle S, K_d \rangle \bigr).
\end{equation}
We write $\mathcal{C}_\mathcal{L}$ for the image of
the positive-definite cone $\s^m_{\succ 0}$ under the  map $\pi_\mathcal{L}$.
We call $\mathcal{C}_\mathcal{L}$ the {\em cone of sufficient statistics}.
The following result explains the duality
between the two red curves in Fig.~1.

\begin{prop} \label{prop:duality}
The cone of sufficient statistics is the convex dual to the cone of concentration matrices.
The basis-free version of this duality states
\begin{equation}
\label{duality_without_basis}
 \mathcal{C}_\mathcal{L} \,=\,\bigl\{ \,S \in \s^m/\mathcal{L}^\perp \,: \langle S, K \rangle > 0
\,\,\,\hbox{for all} \,\,\,K \in \mathcal{K}_\mathcal{L} \bigr\}.
\end{equation}
The basis-dependent version of this duality,
in terms of (\ref{K_with_basis})
and (\ref{pi_with_basis}), states
\begin{equation}
\label{duality_with_basis} \mathcal{C}_\mathcal{L} \,=\,
 \bigl\{\, (t_1,\ldots,t_d) \in \R^d \,: \,
\sum_{i=1}^d t_i \lambda_i > 0 \,
\,\,\hbox{for all} \,\,(\lambda_1,\ldots,\lambda_d) \in \mathcal{K}_\mathcal{L} \bigr\}.
\end{equation}
\end{prop}

\begin{proof}
Let $\mathcal{K}_\mathcal{L}^\vee$ denote the
right-hand side of (\ref{duality_without_basis}) and let $M = \binom{m+1}{2}$.
Using the fact that the $M$-dimensional convex
cone $\s^m_{\succ 0}$ is self-dual, general duality theory
for convex cones implies
$$
\mathcal{K}_\mathcal{L}^\vee \,\,\, = \,\,\,
( \s^m_{\succ 0} \,\cap \, \mathcal{L} )^\vee \,\,\, = \,\,\,
(\s^m_{\succ 0}\,+\, \mathcal{L}^\perp)/\mathcal{L}^\perp \,\,\,\, = \,\,\,\,\mathcal{C}_\mathcal{L}. $$
To derive (\ref{duality_with_basis}) from (\ref{duality_without_basis}), we pick any
basis $U_1,\ldots,U_M$ of $\s^m$ whose first $d$ elements
 serve as the dual basis to $K_1,\ldots,K_d$, and whose last $\binom{m+1}{2}-d$
 elements span $\mathcal{L}^\perp$. Hence
 $\,\langle U_i, K_j \rangle = \delta_{ij}\,$ for all $i,j$.
Every matrix $U$ in $\s^m$ has a unique representation
$\, U \, = \, \sum_{i=1}^M t_i U_i $, and its image under the map
(\ref{pi_with_basis}) equals $\,\pi_\mathcal{L} (U) \, = \, (t_1,\ldots,t_d)$.
For any matrix $\,K = \sum_{j=1}^d \lambda_j K_j\,$ in $\mathcal{L}$
we have $\langle U,K\rangle = \sum_{i=1}^d t_i \lambda_i$, and
this expression is positive for all $K \in \mathcal{K}_{\mathcal{L}}$
if and only if $(t_1,\ldots,t_d)$ lies in $\mathcal{C}_{\mathcal{L}}$.
\end{proof}

It can be shown that both the cone  $\mathcal{K}_\mathcal{L}$ of concentration matrices
and its dual, the cone $\mathcal{C}_\mathcal{L}$  of sufficient statistics,
are well-behaved in the following sense. Their topological closures
are closed convex cones that  are dual to each other, and they are
obtained by respectively intersecting and projecting the closed cone $\s^m_{\succeq 0}$
of positive semidefinite matrices. In symbols,
these closed semi-algebraic cones satisfy
\begin{equation}
\label{closedcone}
 \overline{\mathcal{K}_\mathcal{L}} \,\, = \,\, \mathcal{L} \,\cap \, \s^m_{\succeq 0}
\qquad \hbox{and} \qquad
\overline{\mathcal{C}_\mathcal{L}} \,\, = \,\, \pi_{\mathcal{L}} (\s^m_{\succeq 0}).
\end{equation}
One of our objectives will be to explore the geometry of their boundaries
$$
\partial \mathcal{K}_\mathcal{L} \,\, := \,\,
\overline{\mathcal{K}_\mathcal{L}}   \backslash \mathcal{K}_\mathcal{L}
\qquad \hbox{and} \qquad
\partial \mathcal{C}_\mathcal{L} \,\, := \,\,
\overline{\mathcal{C}_\mathcal{L}}   \backslash \mathcal{C}_\mathcal{L}.
$$
These are convex algebraic hypersurfaces in $\R^d$,
as seen in Example \ref{OldLabsPictures}.
The statistical theory of exponential families implies the
following corollary concerning the geometry of their interiors:

\begin{cor} \label{cor1}
The map $ \, K \mapsto T = \pi_{\mathcal{L}}(K^{-1}) \, $ is a homeomorphism
between the dual pair of open cones $\, \mathcal{K}_\mathcal{L} \,$
and $\,\mathcal{C}_\mathcal{L}$. The inverse map
$T \mapsto K$ takes the sufficient statistics to the
MLE of the concentration matrix. Here,
 $K^{-1}$ is the unique maximizer of the determinant over the
spectrahedron $\pi_\mathcal{L}^{-1}(T) \cap \s^m_{\succ 0} $.
 \end{cor}

One natural first step in studying this picture
is to simplify it by passing to the complex numbers $\C$. This allows
us to relax various inequalities over the real numbers $\R$ and to work
with  varieties over the algebraic closed field $\C$.
We thus identify our model $\mathcal{L}^{-1}$ with its Zariski closure
in the $(\binom{m+1}{2}-1)$-dimensional complex projective space $\PP(\s^m)$.
Let $P_\mathcal{L}$ denote the homogeneous prime ideal
 of all polynomials in $\R[\Sigma] = \R[s_{11}, s_{12}, \ldots, s_{mm}]$
that vanish on $\mathcal{L}^{-1}$. One way to compute $P_\mathcal{L}$ is to
eliminate the entries of an indeterminate symmetric $m {\times}m$-matrix $K$ from the following system of equations:
\begin{equation}
\label{howtogetP}
 \Sigma \cdot K \,=\,{\rm Id}_m \,\, ,\, \, \,\,\,
K \in \mathcal{L}.
\end{equation}
Given a sample covariance matrix $S$, its
maximum likelihood estimate $\hat \Sigma$ can be computed
algebraically, as in Example \ref{OldLabsPictures}.
We do this by solving the following zero-dimensional system
of polynomial equations:
\begin{equation}
\label{symmetricMLE}
\Sigma \cdot K \,=\,{\rm Id}_m \,,\,\,\,
K \in \mathcal{L}\,,\,\,\,
\Sigma - S \in \mathcal{L}^\perp.
\end{equation}
In the present paper we focus on
the systems (\ref{howtogetP}) and (\ref{symmetricMLE}).
Specifically, for various classes
of linear concentration  models $ \mathcal{L}$,
we seek to answer the following three guiding questions.
Example \ref{OldLabsPictures} served to introduce
these three questions. Many more examples will be
featured throughout our discussion.

\smallskip

\noindent {\bf Question 1}. What can be said about the geometry of the
$(d{-}1)$-dimensional projective variety $\mathcal{L}^{-1}$?
What is the {\em degree} of this variety, and what are the minimal generators
of its prime ideal $P_\mathcal{L}$?

\smallskip

\noindent {\bf Question 2}. The map taking a sample covariance matrix $S$ to
its maximum likelihood estimate $\hat \Sigma$ is an algebraic function.
Its degree is the ML degree of the model $\mathcal{L}$.
See \cite[Def.~2.1.4]{DSS}.
Can we find a formula for this ML degree? Which
models $\mathcal{L}$ have their ML degree  equal to  $1$?

\smallskip

\noindent {\bf Question 3}.
The Zariski closure of the boundary
$\,\partial \mathcal{C}_{\mathcal{L}}\,$ of the cone of sufficient
statistics $\, \mathcal{C}_{\mathcal{L}} \,$
is a hypersurface in the complex projective space $\PP^{d-1}$.
What is the defining polynomial $H_{\mathcal{L}}$ of this hypersurface?

\subsection{Generic linear concentration models}
\label{subsec_generic}

In this subsection we examine the case when $\mathcal{L}$
is a generic subspace of dimension $d$ in $\s^m$. Here
 ``generic'' is understood in the sense of algebraic geometry.
In terms of the model representation (\ref{con_model}),
this means that the matrices $K_1, \ldots, K_d$ were chosen at random.
This is precisely the hypothesis made by \cite{NRS}, and one of our
goals is to explain the connection of Questions 1-3
to that paper.

To begin with, we establish the result
that the two notions of degree coincide in the generic case.

\begin{thm} \label{deg=MLdeg}
The ML degree of the model (\ref{con_model}) defined by
a generic linear subspace $\mathcal{L}$ of dimension $d$ in $\s^m$
equals the degree  of the
projective variety $\mathcal{L}^{-1}$.
That degree is denoted
$\phi(m,d)$ and it
satisfies
\begin{equation}
\label{eq:duality}
 \phi(m,d) \,\, = \,\,\,\phi\left(m \,, \binom{m+1}{2}+1-d\,\right).
\end{equation}
\end{thm}

We calculated the ML degree $\phi(m,d)$ of the generic model $\mathcal{L}$
for all matrix sizes up to $m = 6$:
$$
 \begin{array}{c | c c c c c c c c c c c c c c }
       d  & \quad 1 \quad & \quad 2 \quad & \quad 3 \quad & \quad 4 \quad & \quad 5 \quad & \quad 6 \quad & \quad 7 \quad & \quad 8 \quad & \quad 9 \quad & \quad 10 \quad & \quad 11 \quad & \quad 12 \quad & \quad 13 \quad & \\\hline
\phi(3,d) & 1 & 2 & 4 & 4 & 2 & 1 &   &   &   &    &    &    & & \\
\phi(4,d) & 1 & 3 & 9 & 17 & 21 & 21 & 17 & 9 & 3 & 1  &    &  & &  \\
\phi(5,d) & 1 & 4 & 16 & 44 & 86 & 137 & 188 & 212 & 188 & 137 & 86 & 44 & 16
& \cdots \\
\phi(6,d) & 1& 5 & 25 & 90& 240 & 528 & 1016 & 1696 & 2396  & 2886 & 3054 & 2886 & 2396 & \cdots
\end{array}
$$
This table was computed with the software
{\tt Macaulay2}\footnote{www.math.uiuc.edu/Macaulay2/}, using
the commutative algebra techniques discussed in the proof of
Theorem \ref{deg=MLdeg}. At this point, readers from statistics
are advised to skip the algebraic technicalities
in the rest of this section and to go straight to
Section 4 on graphical models.

\smallskip

The last three entries in each row follow from B\'ezout's Theorem
because $P_\mathcal{L}$ is a complete intersection
when the codimension of $\mathcal{L}^{-1}$ in $\PP(\s^m)$
is at most two.
Using the duality relation (\ref{eq:duality}), we conclude
$$ \phi(m,d) \,\,=\,\, (m-1)^{d-1} \quad \,\, \hbox{for} \,\,\, d = 1,2,3. $$
When $\mathcal{L}^{-1}$ has codimension $3$, it is
the complete intersection defined
by three generic linear combinations of the comaximal minors. From this
complete intersection we must remove the variety
of $m \times m$-symmetric matrices of rank $\leq m-2$, which has
also codimension $3$ and has degree $\binom{m+1}{3}$. Hence:
$$ \phi(m,4) \,\,\,=\,\,\,  (m-1)^3-\binom{m+1}{3} \,\,\,= \,\,\,
\frac{1}{6}(5 m - 3) (m - 1) (m - 2) . $$
When $d$ is larger than $4$, this approach leads to
a problem in {\em residual intersection theory}.
A formula due to \cite{Stk}, rederived in recent work by
\cite{CEU} on this subject, implies that
$$ \phi(m,5) \,\,\,=\,\,\, \frac{1}{12} (m-1)(m-2)(7m^2-19 m+6). $$
For any fixed dimension $d$, our ML degree $\phi(m,d)$ seems to be a
polynomial function
of degree $d-1$ in $m$, but it gets progressively more challenging to
compute explicit formulas for these polynomials.

\begin{proof}[Proof of Theorem \ref{deg=MLdeg}]
Let $I$ be the ideal in the polynomial ring
$\,\R[\Sigma] = \R[s_{11}, s_{12}, \ldots, s_{mm}]\,$
that is generated by the $(m{-}1) \times (m{-}1)$-minors of
the symmetric $m {\times} m$-matrix $\Sigma = (s_{ij})$.
\cite{Kot} proved that the Rees algebra
$ \mathcal{R}(I)$ of the ideal $I$ is
equal to the symmetric algebra of $I$. Identifying the
generators of $I$ with the entries of another
symmetric matrix of unknowns $K = (k_{ij})$, we represent
this Rees algebra as
$ \mathcal{R}(I)\, = \, \R[\Sigma, K]/J $
where the ideal $J$ is obtained by eliminating the unknown
$t$ from the matrix equation $\Sigma \cdot K = t \cdot {\rm Id}_m$.
The presentation ideal
$\,J = \langle \Sigma \cdot K \,-\,t\,{\rm Id}_m \rangle \,\cap \,
\R[\Sigma,K]\,$
is  prime and it is homogeneous with respect to the
natural $\N^2$-grading on the polynomial ring $\R[\Sigma,K]$.
Its variety $V(J)$ in $\,\PP^{M-1} \times \PP^{M-1}\,$
is the closure of the set of pairs of symmetric matrices
that are inverse to each other.
Here $ M = \binom{m+1}{2}$.
Both the dimension and the codimension of $V(J)$ is equal to $M-1$.

We now make use of the notion of multidegree
introduced in the text book of \cite{CCA}.
Namely, we consider the {\em bidegree} of the Rees algebra
$ \mathcal{R}(I) =  \R[\Sigma, K]/J $ with respect to its
$\N^2$-grading. This bidegree
 is a homogeneous polynomial
in two variables $x$ and $y$ of degree $M-1$.
Using notation as in \cite[Def.~8.45]{CCA} and
\cite[Thm.~10]{NRS},
we claim that
\begin{equation}
\label{bidegree}
 \mathcal{C}\bigl(\mathcal{R}(I);x,y \bigr) \quad = \quad
\sum_{d=1}^M \phi(m,d) \,x^{M-d} \,y^{d-1} .
\end{equation}
Indeed, the coefficient of $\,x^{M-d} \,y^{d-1} \,$
in the expansion of $\mathcal{C}\bigl(\mathcal{R}(I);x,y \bigr)$
equals the cardinality of
the finite variety $\,V(J) \cap
(\mathcal{M} \times \mathcal{L}) $, where
$\mathcal{L}$ is a generic plane
of dimension $d-1$ in the second factor
$\PP^{M-1}$ and $\mathcal{M}$ is a generic plane
of dimension $M-d$ in the first factor $\PP^{M-1}$.
We now take $\mathcal{M}$ to be the specific plane which is
spanned by the image of $\mathcal{L}^\perp$ and one
extra generic point $S$, representing  a random
sample covariance matrix.
Thus our finite variety is precisely the same as
the one described by the affine equations in  (\ref{symmetricMLE}),
and we conclude that its cardinality equals the
ML degree $\,\phi(m,d)$.

Note that $\,V(J) \cap (\PP^{M-1} \times \mathcal{L}) \,$
can be identified with the variety $V(P_\mathcal{L})$ in
$\PP^{M-1}$. The argument in the previous paragraph relied
on the fact that $P_\mathcal{L}$ is Cohen-Macaulay,
which allowed us to chose any subspace $\mathcal{M}$
for our intersection count provided it is disjoint
from $V(P_\mathcal{L})$ in $\PP^{M-1}$. This proves
that $\phi(m,d)$ coincides with the degree of $V(P_\mathcal{L})$.
The Cohen-Macaulay property of $P_\mathcal{L}$ follows
from a result of \cite{HVV} together with the aforementioned work of \cite{Kot} which shows that the ideal
$I$ has sliding depth.
Finally, the duality (\ref{eq:duality}) is obvious
for the coefficients of the bidegree (\ref{bidegree})
of the Rees algebra $ \mathcal{R}(I)$ since its presentation ideal $J$
is symmetric under swapping $K$ and $\Sigma$.
\end{proof}

We now come to our third question, which is  to determine the
Zariski closure $V(H_\mathcal{L})$ of the boundary of the cone
$\,\mathcal{C}_{\mathcal{L}} \,=\,\mathcal{K}_{\mathcal{L}}^\vee$.
Let us assume  now that $\mathcal{L}$ is any $d$-dimensional
linear subspace of $\s^m$, not necessarily generic.
The Zariski closure of $\partial \mathcal{K}_{\mathcal{L}}$
is the hypersurface $\,\{{\rm det}(K) = 0\} \,$ given by the
vanishing of the determinant of $ \, K = \sum_{i=1}^d \lambda_i K_i$.
This determinant is a polynomial of degree $d$ in
$\lambda_1,\ldots,\lambda_d$. Our task is to compute the
{\em dual variety} in the sense of projective algebraic geometry of each
irreducible component of this hypersurface.
See~\cite[\S 5]{NRS} for basics on projective duality.
We also need to compute the
dual variety for its singular locus, and for the
singular locus of the singular locus, etc.

Each singularity stratum encountered along the way needs
to be decomposed into irreducible components, whose
duals need to be examined.
If such a component has a real point that lies in
$\partial \mathcal{K}_{\mathcal{L}}$ and if its dual
variety is a hypersurface
then that hypersurface appears in $H_\mathcal{L}$.
How to run this procedure in practice is shown
in Example \ref{ex_5graph}.
For now, we summarize the construction informally as follows.

\begin{prop} \label{fuzzy}
Each irreducible hypersurface in the Zariski closure of
$\,\partial \mathcal{C}_{\mathcal{L}}\,$ is  the
projectively dual variety to some irreducible component
of the hypersurface $\{{\rm det}(K) = 0\}$,
or it is dual to some irreducible variety
further down in the   singularity stratification
of the hypersurface $\{{\rm det}(K) = 0\} \subset \PP^{d-1}$.
\end{prop}

The singular stratification of $\{{\rm det}(K) = 0\}$
can be computed by applying primary decomposition to the
ideal of $p {\times} p$-minors of $K$ for $1 \leq p \leq m$.
If $I$ is any minimal prime of such a determinantal ideal
then its dual variety is computed as follows.
Let $c = {\rm codim}(I)$ and consider the
Jacobian matrix of $I$.
The rows of the Jacobian matrix are the derivatives of the generators of $I$
with respect to the unknowns $\lambda_1,\ldots,\lambda_d$.
Let $J$ be the ideal generated by $I$ and the $c {\times} c$-minors
of the matrix formed by augmenting the Jacobian matrix by the extra row
$(t_1,t_2,\ldots,t_d)$. We saturate $J$ by the $c {\times} c$-minors
of the Jacobian, and thereafter we compute the elimination ideal
$\,J \cap \R[t_1,t_2,\ldots,t_d]$. If this elimination ideal is principal,
we retain its generator. The desired polynomial $H_\mathcal{L}$
is the product of these principal generators, as $I$ runs over
all such minimal primes whose variety has a real point on
 the convex hypersurface $\partial \mathcal{K}_{\mathcal{L}}$.

Proposition \ref{fuzzy} is visualized also in Fig.~\ref{dualpillow} below. Let us now apply this result in the case
when the subspace $\mathcal{L}$ is generic of dimension $d$.
The ideal of $p {\times} p$-minors of $K$ defines a subvariety of
$\PP^{d-1}$, which is irreducible whenever it is positive-dimensional
(by Bertini's Theorem).  It is known from \cite[Prop.~5]{NRS} that
the dual variety to that determinantal variety is a hypersurface if and only if
\begin{equation}
\label{pataki}
  \binom{m-p+2}{2} \,\leq \, d-1
\qquad \hbox{and} \qquad
\binom{p}{2} \,\leq\, \binom{m+1}{2} -d+1.
\end{equation}
Assuming that these inequalities hold, the
dual hypersurface is defined by an irreducible
homogeneous polynomial whose degree we denote by
$\delta(d-1,m,p-1)$.  This notation is consistent with
\cite{NRS} where this number is called the
{\em algebraic degree of semidefinite programming (SDP)}.

\begin{cor}
For a generic $d$-dimensional subspace $\mathcal{L}$ of
$\s^m$, the polynomial $H_{\mathcal{L}}$ is the product
of irreducible polynomials of degree $\delta(d-1,m,p-1)$.
That number is  the algebraic degree of semidefinite
programming. Here $p$ runs over
integers that satisfy (\ref{pataki}) and
 $\partial \mathcal{K}_{\mathcal{L}}$ contains a
matrix of rank $p-1$.
\end{cor}

\section{Diagonal Matrices, Matroids and Polytopes}
\label{sec_diagonal}

This section concerns the case when $\mathcal{L}$ is
a $d$-dimensional space consisting only of
diagonal matrices in $\s^m$. Here, the set
$\,\mathcal{L}^{-1}_{\succ 0}\,$ of covariance matrices in
the model also consists of diagonal
matrices only, and we may restrict our considerations to
the space $\R^m$ of diagonal matrices in $\s^m$.
Thus, throughout this section, our ambient space is
$\R^m$, and we identify $\R^m$ with its dual vector space via the
standard inner product $\langle u,v \rangle = \sum_{i=1}^m u_i v_i$.
We fix any $d \times m$-matrix $A$ whose rows space equals
$\mathcal{L}$, and we assume
that $\mathcal{L} \cap \R^m_{>0} \not= \emptyset$.
We consider the
induced projection of the open positive orthant
\begin{equation}
\label{A_map}
 \,\pi : \R^m_{> 0} \rightarrow \R^d \,,\,\, x \mapsto A x .
\end{equation}
Since $\mathcal{L} = {\rm rowspace}(A)$ contains a strictly positive vector,
the image of $\pi$ is a pointed polyhedral cone,
namely $\mathcal{C}_{\mathcal{L}} = {\rm pos}(A)$
is the cone spanned by the columns of $A$.
Each fiber of $\pi$ is a bounded convex polytope,
and maximum likelihood estimation amounts to finding
a distinguished point $\hat x$ in that fiber.

The problem of characterizing the existence of the MLE in this situation
amounts to a standard problem of
{\em geometric combinatorics} (see e.g.~\cite{Zie}), namely,
to computing the
facet description of the convex polyhedral cone spanned by
the columns of $A$.
For a given vector $t \in \R^d$ of sufficient statistics, the
maximum likelihood estimate exists in this diagonal concentration model
if and only if $t$ lies in the interior of the cone ${\rm pos}(A)$.
This happens if and only if all facet inequalities are strict for~$t$.

This situation is reminiscent of {\em Birch's Theorem} for toric models
in algebraic statistics \cite[Theorem 1.10]{ASCB}, and,
indeed, the combinatorial set-up for deciding the existence of the MLE
is identical to that for toric models. For a statistical perspective see
\cite{EFRS}. However, the algebraic structure here is not that of
 toric models, described in \cite[\S 1.2.2]{ASCB}, but that of the
linear models in \cite[\S 1.2.1]{ASCB}.

Our model here is not toric but it
is the coordinatewise reciprocal of an open polyhedral cone:
$$ \mathcal{L}^{-1}_{> 0}  \,\, = \,\, \bigl\{ \, u \in \R^m_{> 0} : \,
u^{-1}  = ( u_1^{-1}, u_2^{-1},\ldots, u_m^{-1}) \in \mathcal{L} \,\bigr\}. $$
As in Section 2, we view its Zariski closure $\mathcal{L}^{-1}$
as a subvariety in complex projective space:
$$ \mathcal{L}^{-1} \,\, = \,\, \bigl\{ \, u \in \PP^{m-1} \, : \,
u^{-1}  = ( u_1^{-1}, u_2^{-1},\ldots, u_m^{-1}) \in \mathcal{L} \,\bigr\}. $$
Maximum likelihood estimation means intersecting
the variety $\mathcal{L}^{-1}$ with the fibers of $\pi$.

\begin{example}
Let $m=4$, $d=2$ and take $\mathcal{L}$ to be the row space of the matrix
$$\,A \,\,= \,\,
\begin{pmatrix} 3 && 2 && 1 && 0 \\ 0 && 1 && 2 && 3 \end{pmatrix}.$$
The corresponding statistical model consists of all
multivariate normal distributions on $\R^4$  whose
concentration matrix has the diagonal form
$$ K \quad = \quad \begin{pmatrix}3 \lambda_1 & 0 &  0 & 0 \\
0 &  2 \lambda_1 + \lambda_2 & 0 & 0 \\
0 & 0 & \lambda_1 +2 \lambda_2 & 0 \\
0 & 0 & 0 & 3 \lambda_2
\end{pmatrix}.
 $$
Our variety $\mathcal{L}^{-1}$ is the curve in $\PP^3$
parametrized by the inverse diagonal matrices
which we write as
$\,K^{-1} = {\rm diag}(x_1,x_2,x_3,x_4)$.
The prime ideal $P_\mathcal{L}$ of this curve
 is generated by three quadratic equations:
$$
 x_2 x_3 - 2 x_2 x_4+x_3 x_4 \,=\, 2 x_1 x_3 - 3 x_1 x_4 + x_3 x_4 \,=\,
     x_1 x_2 - 3 x_1 x_4 + 2 x_2 x_4 \,\,\, = \,\,\, 0 .
$$
Consider any sample covariance matrix $S = (s_{ij})$, with sufficient statistics
$$ t_1 \,=\, 3 s_{11} + 2 s_{22} + s_{33}  \,> \,0\quad \hbox{and} \quad
     t_2 \, = \, s_{22} + 2 s_{33} + 3 s_{44} \, > \, 0. $$
The MLE for these sufficient statistics is the unique positive solution
$\hat x$ of the
three quadratic equations above, together with the two linear equations
$$
 3 x_1 + 2 x_2 + x_3 \,= \,t_1 \quad \hbox{and} \quad
     x_2 + 2 x_3 + 3 x_4 \, = \, t_2 .
$$
We find that $\hat x$ is
 an algebraic function
of degree $3$ in the sufficient statistics $(t_1,t_2)$, so
the ML degree of the model $K$ equals $3$.
This is consistent  with formula (\ref{binobeta}) below, since
$\binom{4-1}{2-1} = 3$. \qed
\end{example}

We now present the solutions to our three guiding problems for
arbitrary $d$-dimensional subspaces $\mathcal{L}$ of the space $\R^m$ of
$m {\times} m$-diagonal matrices.
The degree of the projective variety $\mathcal{L}^{-1}$ and its prime ideal
$P_{\mathcal{L}}$ are known from work of \cite{Ter}
and its refinements due to \cite{PS}.
 Namely, the degree of
$\, \mathcal{L}^{-1}\,$ equals the {\em beta-invariant} of the
rank $m-d$ matroid on $[m] = \{1,2,\ldots,m\}$
associated with $\mathcal{L}$.
We denote this beta-invariant by $\beta(\mathcal{L})$.
For matroid basics see \cite{Wh}.

The beta-invariant $\beta(\mathcal{L})$ is known to equal
the number of bounded regions in the $(m-d)$-dimensional
hyperplane arrangement (cf.~\cite{Zas}) obtained by intersecting the
affine space $u + \mathcal{L}^\perp$ with the $m$
coordinate hyperplanes $\{x_i = 0\}$.
Here $u$ can be any generic vector in $\R^m_{> 0}$. One of these regions,
 namely the one containing $u$, is   precisely the fiber of $\pi$.
If $\mathcal{L}$ is a generic $d$-dimensional linear subspace of
$\R^m$, meaning that the above matroid is the uniform matroid, then
the beta-invariant equals
\begin{equation}
\label{binobeta}
 \beta(\mathcal{L}) \,=\,  \binom{m-1}{d-1}.
 \end{equation}
For non-generic subspaces $\mathcal{L}$,
this binomial coefficient is always an upper bound for $\beta(\mathcal{L})$.

\begin{thm}
\label{thm:beta} {\rm (\cite{PS, Ter})}
The degree of the projective variety $\mathcal{L}^{-1}$ equals
the beta-invariant $\beta(\mathcal{L})$.
Its prime ideal $P_\mathcal{L}$ is generated
by the homogeneous polynomials
\begin{equation}
\label{circuits} \sum_{i \in {\rm supp}(v)} \, v_i \cdot \prod_{j\not=i} x_j
\end{equation}
where $v$ runs over all non-zero vectors
of minimal support in $\mathcal{L}^\perp$.
\end{thm}

For experts in combinatorial commutative algebra, we note that
\cite{PS} actually prove the following stronger results.
The homogeneous polynomials (\ref{circuits}) form a
{\em universal Gr\"obner basis} of $P_{\mathcal{L}}$.
The initial monomial ideal of $P_{\mathcal{L}}$
with respect to any term order is the Stanley-Reisner ideal
of the corresponding
{\em broken circuit complex} of the matroid of $\mathcal{L}$. Hence the Hilbert series of
$P_{\mathcal{L}}$ is the rational function obtained by
dividing the
 h-polynomial $h(t)$ of the broken circuit complex
 by $(1-t)^d$. In particular,
the degree of $P_\mathcal{L}$ is
the number $h(1) = \beta(\mathcal{L})$ of broken circuit bases \cite[\S 7]{Wh}.

We next consider Question 2 in the diagonal case.
The maximum likelihood map takes each vector $t$ in the cone of sufficient statistics
$\,\mathcal{C}_{\mathcal{L}} =
{\rm pos}(A)\,$ to a point of its fiber, namely:
\begin{equation}
\label{argmax}
 \hat x  \,\, = \,\, {\rm argmax} \bigl\{
\sum_{i=1}^m {\rm log}(x_i) \,\,:\,\,
 x \in \R^m_{>0} \,\,\hbox{and} \,\, A x = t \bigr\} .
\end{equation}
 This is the unique point in the polytope
 $\,\pi^{-1}(t) = \{ x \in \R^m_{>0} \,\,\hbox{and} \,\, A x = t \}\, $ which maximizes the product
 $x_1 x_2 \cdots x_n$ of the coordinates. It is also the unique point in
 $\pi^{-1}(t)$ that lies in the reciprocal linear variety $\mathcal{L}^{-1}$.
 In the linear programming literature, the point $\hat x$ is known as the
{\em analytic center} of the polytope $\pi^{-1}(t)$.
In Section 2 we discussed the extension
of this concept
from linear programming to semidefinite programming:
the {\em analytic center of a spectrahedron} is the unique point $\hat \Sigma$ at which the determinant function attains its maximum.
For an applied perspective see \cite{VBW}.

For any linear subspace $\mathcal{L}$ in $\s^m$, the algebraic degree of the
maximum likelihood map $t \mapsto \hat \Sigma$ is always less than or equal
to the degree of the projective variety $\mathcal{L}^{-1}$.
We saw in Theorem \ref{deg=MLdeg} that these degrees
are equal for generic $\mathcal{L}$.
We next show that the same conclusion holds for
diagonal subspaces $\mathcal{L}$.

\begin{cor}
\label{MLbeta}
The ML degree of any diagonal linear concentration model
$\mathcal{L} \subset \R^m  \subset \s^m$
is equal to the beta-invariant $\beta(\mathcal{L})$ of the
corresponding matroid of rank $m-d$ on
$\{1,2,\ldots,m\}$.
\end{cor}

\begin{proof}
The beta-invariant $\beta(\mathcal{L})$ counts
the bounded regions in the arrangement of hyperplanes
arising from the given facet description of the polytope $\pi^{-1}(t)$.
{\em Varchenko's Formula} for linear models, derived in
\cite[Theorem 1.5]{ASCB}, states that the optimization problem
(\ref{argmax}) has precisely one real critical point in each bounded region,
and that there are no other complex critical points.
\end{proof}

A fundamental question regarding the ML degree of any
class of algebraic statistical models is to characterize those models
which have ML degree one. These are the models
whose maximum likelihood estimator is a rational function in
the sufficient statistics \cite[\S 2.1]{DSS}.
In the context here, we have the following characterization
of matroids whose beta-invariant $\beta(\mathcal{L})$ equals one.

\begin{cor} The ML degree $\beta(\mathcal{L})$
of a diagonal linear concentration model $\mathcal{L}$  is equal to one
if and only if the matroid of $\mathcal{L}$ is the
graphic matroid of a series-parallel graph.
\end{cor}

\begin{proof}
The equivalence of series-parallel and $\beta = 1$
first appeared in \cite[Theorem 7.6]{Bry}.
\end{proof}

We now come to Question 3 which concerns the
duality of convex cones in Proposition \ref{prop:duality}.
In the diagonal case, the geometric view on this
duality is as follows.
The cone of sufficient statistics equals
$\,\mathcal{C}_{\mathcal{L}} = {\rm pos}(A)\,$ and its convex dual
is the cone
$\,\mathcal{K}_{\mathcal{L}} =  {\rm rowspace}(A) \cap \mathcal{L}$.
Both cones are convex, polyhedral, pointed, and have dimension $d$.
By passing to  their cross sections with suitable affine
hyperplanes, we can regard the
two cones $\,\mathcal{C}_\mathcal{L}\,$ and
$\,\mathcal{K}_\mathcal{L}\,$ as a dual pair of $(d-1)$-dimensional
convex polytopes.

The hypersurface $\{{\rm det}(K) = 0\}$ is a union of $m$
hyperplanes. The strata in its singularity stratification,
discussed towards the end of Section 2,
correspond to the various faces $F$ of the polytope $\mathcal{K}_\mathcal{L}$. The dual
variety to a face $F$ is the complementary
face of the dual polytope $\mathcal{C}_\mathcal{L}$, and hence the
codimension of that dual variety equals one if and only if $F$
is a vertex (= $0$-dimensional face) of $\mathcal{K}_\mathcal{L}$.
This confirms that the polynomial $H_\mathcal{L}$ sought in Question 3
is the product of all facet-definining linear forms of $\mathcal{C}_\mathcal{L}$.

Corollary \ref{cor1} furnishes a homeomorphism
$\, u \mapsto A u^{-1}\,$ from the interior of the polytope $\mathcal{K}_{\mathcal{L}}$
onto the interior of its dual polytope $\mathcal{C}_{\mathcal{L}}$.
The inverse to the rational function  $\,u \mapsto A u^{-1} \,$ is
an algebraic function whose degree is the beta-invariant $\beta(\mathcal{L})$.
This homeomorphism is the natural generalization, from simplices
to arbitrary polytopes, of the classical
{\em Cremona transformation} of projective geometry.
We close this section with a nice $3$-dimensional example
which illustrates this homeomorphism.

\begin{example}[How to morph a cube into an octahedron] \
Fix  $m=8$, $d=4$, and  $\mathcal{L}$ the row space of
$$ A \quad = \quad
\begin{pmatrix}
& 1 &          - 1 & \,\,\,0\, & \phantom{-}0 & \,\,\,0\, & \phantom{-}0 & \\
& 0 & \phantom{-}0 & \,\,\,1\, &          - 1 & \,\,\,0\, & \phantom{-}0 & \\
& 0 & \phantom{-}0 & \,\,\,0\, & \phantom{-}0 & \,\,\,1\, &         -  1 & \\
& 1 & \phantom{-}1 & \,\,\,1\, & \phantom{-}1 & \,\,\,1\, & \phantom{-}1 & \\
\end{pmatrix} . $$
We identify the cone
$\,\mathcal{K}_{\mathcal{L}} \,=\, {\rm rowspace}(A) \cap \R_{>0}^6\,$
with $\,\{\lambda \in \R^4 \,:\, \lambda \cdot A > 0\}$.
This is the cone over the $3$-cube,
which is obtained by setting $\lambda_4 = 1$.
The dual cone $\,\mathcal{C}_{\mathcal{L}} = {\rm pos}(A)\,$ is
spanned by the six columns of the matrix $A$. It is the cone
over the octahedron, which is obtained by setting $t_4= 1$.

We write the homeomorphism $u \mapsto  A u^{-1}$ between
these two four-dimensional
cones in terms of the coordinates of $\lambda $ and $t$.
Explicitly, the equation $\,t = A \cdot (\lambda A)^{-1} \,$
translates into the scalar equations:
\begin{eqnarray*}
& t_1 \,\,
= & \frac{1}{\lambda_4 + \lambda_1} - \frac{1}{\lambda_4 - \lambda_1} ,\\
& t_2 \,\,= & \frac{1}{\lambda_4 + \lambda_2} - \frac{1}{\lambda_4 - \lambda_2} ,\\
& t_3 \,\,= & \frac{1}{\lambda_4 + \lambda_3} - \frac{1}{\lambda_4 - \lambda_3} ,\\
& t_4 \,\, =  & \frac{1}{\lambda_4 + \lambda_1} + \frac{1}{\lambda_4 - \lambda_1} +
  \frac{1}{\lambda_4 + \lambda_2} + \frac{1}{\lambda_4 - \lambda_2} +
 \frac{1}{\lambda_4 + \lambda_3} + \frac{1}{\lambda_4 - \lambda_3} .
\end{eqnarray*}
Substituting $\lambda_4= 1$, we get the bijection
$\,(\lambda_1,\lambda_2,\lambda_3) \mapsto ( t_1/t_4,t_2/t_4,t_3/t_4)\,$
between  the open cube $(-1,+1)^3$ and the
open octahedron
$\{\,t \in \R^3 \,:\, |t_1|+|t_2|+|t_3| < 1\}$.
The inverse map $\,t \mapsto \lambda\,$ is an
algebraic function of degree $\beta(\mathcal{L}) = 7 $.
That the ML degree of this model is $7$ can be seen as follows. The fibers
$\pi^{-1}(t)$ are the convex polygons
which can be obtained
from a regular hexagon by parallel displacement of its six
edges. The corresponding arrangement of six lines
has $7$ bounded regions.
\qed
\end{example}

\section{Gaussian Graphical Models}
\label{sec:graphical}

An {\em undirected Gaussian graphical model} arises when
the subspace $\mathcal{L}$ of $\s^m$ is defined by the vanishing of
some off-diagonal entries of the concentration matrix $K$.
We fix a graph $G = ([m],E)$ with vertex set $[m] = \{1,2,\ldots,m\}$
and whose edge set $E$ is assumed to contain all self-loops.
A basis for $\mathcal{L}$ is the set
$\{K_{ij}\mid (i,j)\in E\}$ of matrices $K_{ij}$
with a single 1-entry in position $(i,j)$
and 0-entries in all other positions. We shall use
the notation $\,\mathcal{K}_G, \mathcal{C}_G, P_G \,$
for the objects
$\,\mathcal{K}_\mathcal{L}, \mathcal{C}_\mathcal{L}, P_\mathcal{L} $,
respectively. Given a sample covariance matrix $S$,
the set $\,{\rm fiber}_G(S)\,$
consists of all positive definite matrices $\,\Sigma\in\s^m_{\succ 0}\,$ with
$$\Sigma_{ij}=S_{ij} \qquad \textrm{for all }\, (i,j)\in E.$$

The cone of concentration matrices $\mathcal{K}_G$
is important for semidefinite matrix completion problems.
Its closure was denoted $\mathcal{P}_G$ by \cite{L1, L2}.
The dual cone $\,\mathcal{C}_G \,$
consists of all partial matrices $T \in \R^E$ with entries in
positions $(i,j)\in E$, which can be extended to a full positive definite
matrix. So, maximum likelihood estimation in Gaussian graphical models
corresponds to the classical positive definite matrix completion problem
\citep{BJT, GJSW, BJL, L2}.
In this section we investigate our three guiding questions,
first for chordal graphs, next for the chordless $m$-cycle $C_m$,
then for all graphs with five or less vertices, and finally
for the $m$-wheel $W_m$.

\subsection{Chordal graphs}
\label{chordal_graphs}

A graph $G$ is {\em chordal} (or {\em decomposable}) if
every induced $m$-cycle in $G$ for $m \geq 4$ has a chord. A
theorem due to \cite{GJSW} fully resolves
Question 3 when $G$ is chordal. Namely, a  partial matrix $T\in \R^E$
lies in the cone $\mathcal{C}_{G}$ if and only if all principal minors
$T_{CC}$ indexed by cliques $C$ in $G$ are positive definite.
The ``only if'' direction in this statement is
true for all graphs $G$, but the ``if'' direction holds
only when $G$ is chordal. This result
is equivalent to the characterization of chordal graphs as those
that have {\em sparsity order} equal to one, i.e., all
extreme rays of $\mathcal{K}_G$ are matrices of rank one.
We refer to \cite{AHMR} and \cite{L1} for details.
From this characterization of chordal graphs
in terms of sparsity order,
we infer the following description
of the Zariski closure of the boundary
of $\mathcal{C}_G$.

\begin{prop}
For a chordal graph $G$, the defining polynomial $\,H_{G}\,$
of $\,\partial \mathcal{C}_G\,$ is equal to
$$ H_G \,\,\,\,=\prod_{C \,\textrm{maximal} \atop
\textrm{clique of }G} \det(T_{CC}).$$
\end{prop}

We now turn to Question 2 regarding the ML degree of
a Gaussian graphical model $G$.
This number is here simply denoted by $\textrm{ML-degree}(G)$.
Every chordal graph is a clique sum of complete graphs.
We shall prove that the ML degree is multiplicative
with respective to taking clique sums.

\begin{lem} \label{MLdegree_cliquesum}
Let $G$ be a clique sum of $n$ graphs $G_1, \ldots , G_n$.
Then the following equality holds:
$$\textrm{ML-degree}\,(G)\,\,\,= \,\,\,
\prod_{i=1}^{n} \textrm{ML-degree}\,(G_i).$$
\end{lem}

\begin{proof}
We first prove this statement for $n=2$.  Let $G$ be a graph which
can be decomposed in disjoint subsets $(A, B, C)$ of the vertex
set $V$, such that $C$ is a clique and separates $A$ from $B$.
Let $G_{[W]}$ denote the induced subgraph on a vertex subset $W\subset V$.
So, we wish to prove:
\begin{equation}
\label{mldegfactor1}
\textrm{ML-degree}\,(G) \,\,\, = \,\,\,
\textrm{ML-degree}\,(G_{[A\cup C]})\cdot \textrm{ML-degree}\,(G_{[B\cup C]}).
\end{equation}
Given a generic matrix $\,S\in\mathbb{S}^m$, we fix
$\,\Sigma\in\mathbb{S}^m\,$ with entries $\,\Sigma_{ij}=S_{ij}\,$
for $\,(i,j)\in E\,$ and unknowns
$\,\Sigma_{ij}=z_{ij}\,$ for $\,(i,j)\notin E$.
The ML degree of $G$ is the number of complex solutions to
the equations
\begin{equation}
\label{mldegfactor2}
(\Sigma^{-1})_{ij}=0 \qquad \textrm{for all } (i,j)\notin E.
\end{equation}
Let $K = \Sigma^{-1}$ and denote by $K^1=(\Sigma_{[A\cup C]})^{-1}$
(respectively, $K^2=(\Sigma_{[B\cup C]})^{-1}$) the inverse of the submatrix
of $\Sigma$ corresponding to the induced subgraph on $A\cup C$
(respectively, $B\cup C$).
Using Schur complements, we can see that these matrices
are related by the following block structure:
$$ K=\begin{pmatrix} K^1_{AA} & K^1_{AC} & 0 \\ K^1_{CA} & K_{CC} & K^2_{CB} \\ 0 & K^2_{BC} & K^2_{BB} \end{pmatrix},
\quad
K^1=\begin{pmatrix} K^1_{AA} & K^1_{AC} \\ K^1_{CA} & K^1_{CC} \end{pmatrix}, \quad
K^2=\begin{pmatrix} K^2_{CC} & K^2_{CB} \\ K^2_{BC} & K^2_{BB} \end{pmatrix}.
$$
This block structure reveals that,
when solving the system (\ref{mldegfactor2}),
one can solve for the variables $z_{ij}$
corresponding to missing edges in the subgraph
$A\cup C$ independently from the variables over $B\cup C$ and $A\cup B$.
This implies the equation  (\ref{mldegfactor1}).
Induction yields the theorem for $n \geq 3$.
\end{proof}

The following theorem characterizes chordal graphs
in terms of their ML degree. It extends the equivalence
of parts (iii) and (iv) in \cite[Thm.~3.3.5]{DSS}
from discrete to Gaussian models.

\begin{thm}
A graph $G$ is chordal if and only if
$\textrm{ML-degree}(G) = 1$.
\end{thm}

\begin{proof}
The if-direction follows from Lemma
\ref{MLdegree_cliquesum} since every
chordal graph is a clique sum of complete
graphs, and a complete graph trivially has
ML degree one. For the only-if direction suppose
that $G$ is a graph that is not chordal.
Then $G$ contains the chordless cycle
$C_m$ as an induced subgraph for some $m \geq 4$.
It is easy to see that the ML degree of any
graph is bounded below by that of any induced
subgraph. Hence what we must prove is that the
chordless cycle $C_m$ has strictly positive ML degree.
This is precisely the content of Lemma \ref{MLispositive} below.
\end{proof}

We now come to Question 1 which concerns the
homogeneous prime ideal $P_G$ that defines the
Gaussian graphical model as a subvariety of $\PP(\s^m)$.
Fix a symmetric $m {\times} m $-matrix of unknowns $\Sigma = (s_{ij})$
and let $\Sigma_{ij}$ denote the
comaximal minor obtained by deleting the $i$th row
and the $j$th column from $\Sigma$.
We shall define several ideals
in $\R[\Sigma]$ that approximate $P_G$.
The first is the saturation
\begin{equation}
\label{idealbysaturation}
 P'_G\,\,=\,\, \bigl( \,
\langle \,{\rm det}(\Sigma_{ij}) \mid (i,j)\in E \,\rangle : \langle
{\rm det}(\Sigma)\rangle^\infty \,\bigr).
\end{equation}
This ideal is contained in the desired prime ideal, i.e.~$P'_G \subseteq P_G$.
The two ideals have the same radical,
but it might happen that they are not equal. One disadvantage
of the ideal $P'_G$ is that the saturation step (\ref{idealbysaturation})
is computationally expensive and terminates only for very small graphs.

A natural question is whether the prime ideal $P_G$ can be constructed
easily from the prime ideals $P_{G_1}$ and $P_{G_2}$ when
$G$ is a clique sum of two smaller graphs $G_1 $ and $G_2$.
As in the proof of Lemma \ref{MLdegree_cliquesum}, we
partition $[m] = A \cup B \cup C$,
where $G_1$ is the induced subgraph on $A \cup C$,
and $G_2$ is the induced subgraph on $B \cup C$.
If $|C| = c$ then we say that $G$ is a {\em $c$-clique sum} of $G_1$ and $G_2$.

The following ideal is contained in $P_G$ and defines
the same algebraic variety in the open cone $\s^m_{\succ 0}$:
\begin{equation}
\label{cliquesumideal}
 P_{G_1} \,+\,P_{G_2} \,+ \,\bigl\langle  \,
\hbox{\rm $(c{+}1) {\times} (c{+}1)$-minors of
$\Sigma_{A \cup C,B\cup C} $} \,\bigr\rangle.
\end{equation}
One might guess that  (\ref{cliquesumideal}) is equal to $P_G$, at least
up to radical, but  this fails for $c \geq 2$.
Indeed we shall see in Example \ref{ex:seth} that
the variety of (\ref{cliquesumideal})
can have extraneous components on the boundary
$\,\s^m_{\succeq 0} \backslash \s^m_{\succ 0} \,$
of the semidefinite cone. We do conjecture, however,
that this equality holds for $c \leq 1$. This is easy to prove
for $c = 0$ when $G$ is disconnected and is the disjoint union of
$G_1$ and $G_2$.
The case $c=1$ is considerably more delicate.
At present, we do not have a proof that (\ref{cliquesumideal})
is prime for $c=1$,
but we believe that even a lexicographic Gr\"obner basis for $P_G$
can be built by taking the union of such
Gr\"obner basis for $P_{G_1}$ and $P_{G_2}$ with the
$2 \times 2$-minors of $\Sigma_{A \cup C,B\cup C} $.
This conjecture would imply the following.

\begin{conj} \label{quadraticgraphs}
The prime ideal $P_G$ of an undirected Gaussian graphical model
is generated in degree $\leq 2$ if and only if each
connected component of the graph $G$ is a $1$-clique sum of complete graphs.
In this case, $P_G$ has a Gr\"obner basis consisting of
entries of $\Sigma$ and $2\times 2$-minors of $\Sigma$.
\end{conj}

This conjecture is an extension of the results
and conjectures for (directed) trees in \cite[\S 5]{Sul}.
Formulas for the degree of $P_G$ when $G$ is a tree are
found in \cite[Corollaries 5.5 and 5.6]{Sul}.
The ``only if'' direction in the first sentence of
Conjecture \ref{quadraticgraphs} can be shown as follows.
If $G$ is not chordal then it contains an $m$-cycle
($m \geq 4$) as an induced subgraph, and, this gives rise
to cubic generators  for $P_G$, as seen in Subsection 4.2 below.
If $G$ is chordal but is not a $1$-clique sum of complete graphs,
then its decomposition involves a $c$-clique sum for some $c \geq 2$,
and the right hand side of (\ref{cliquesumideal})
contributes a minor of size $c+1 \geq 3$ to
the minimal generators of $P_G$. The algebraic structure of
chordal graphical models
is more delicate in the Gaussian case then in the discrete
case, and there is no Gaussian analogue to the
characterizations of chordality in (i) and (ii) of \cite[Theorem 3.3.5]{DSS}.
This is highlighted by the following example
which was suggested to us by Seth Sullivant.

\begin{example} \label{ex:seth}
Let $G$ be the graph on $m=7$ vertices
consisting of the triangles
$\{i,6,7\}$ for $i = 1,2,3,4,5$.
Then $G$ is chordal because it is the
$2$-clique sum of these five triangles.
The ideal $P_G$ is minimally generated by
$105$ cubics and one quintic. The cubics
are spanned by the $3 \times 3$-minors
of the matrices $\Sigma_{A \cup C, B \cup C}$
where $C = \{6,7\}$ and $\{A,B\}$ runs over all
unordered partitions of $\{1,2,3,4,5\}$.
These minors do not suffice to define
the variety $V(P_G)$ set-theoretically.
For instance, they vanish whenever the last two rows and columns
of $\Sigma$ are zero. The additional quintic generator of $P_G$ equals
\begin{eqnarray*}
&  \, s_{12} s_{13} s_{24} s_{35} s_{45} - s_{12} s_{13} s_{25} s_{34} s_{45} - s_{12} s_{14} s_{23} s_{35} s_{45}
 + s_{12} s_{14} s_{25} s_{34} s_{35} \\ &
 + s_{12} s_{15} s_{23} s_{34} s_{45} - s_{12} s_{15} s_{24} s_{34} s_{35}
 + s_{13} s_{14} s_{23} s_{25} s_{45} - s_{13} s_{14} s_{24} s_{25} s_{35}
\\ & - s_{13} s_{15} s_{23} s_{24} s_{45}
 + s_{13} s_{15} s_{24} s_{25} s_{34} + s_{14} s_{15} s_{23} s_{24} s_{35} - s_{14} s_{15} s_{23} s_{25} s_{34}.
\end{eqnarray*}
This polynomial is the {\em pentad} which is relevant for
factor analysis \cite[Example 4.2.8]{DSS}.
\qed
\end{example}

Given an undirected graph $G$ on $[m]$, we define its
{\em Sullivant-Talaska ideal} ${\rm ST}_G$ to be
the  ideal in $\R[\Sigma]$ that is generated by the
following collection of minors of $\Sigma $. For any
submatrix $\Sigma_{A,B}$ we include in ${\rm ST}_G$ all
$c \times c$-minors of $\Sigma_{A,B}$ provided $c$ is the smallest
cardinality of a set $C$ of vertices that separates
$A$ from $B$ in $G$. Here, $A, B $ and $C$ need not be disjoint,
and separation means that any
path from a node in $A$ to a node in $B$ must pass through
a node in $C$. \cite{ST} showed
that the generators of ${\rm ST}_G$ are precisely those
subdeterminants of $\Sigma$ that lie in $P_G$,
and both ideals cut out the same variety in the positive
definite cone $\s^m_{\succ 0}$. However, generally their varieties differ
on the boundary of that cone, even for chordal graphs $G$,
as seen in Example  \ref{ex:seth}. In our experiments,
we found that ${\rm ST}_G$ can often be computed quite fast,
and it frequently coincides with the
desired prime ideal $P_G$.

\subsection{The chordless $m$-cycle}
\label{subsec_cycle}

We next discuss  Questions 1, 2, and 3 for the simplest non-chordal graph,
namely, the $m$-cycle $C_m$. Its Sullivant-Talaska ideal
${\rm ST}_{C_m}$ is generated by the $3 \times 3$-minors of the
submatrices $\Sigma_{A,B}$ where $A = \{i,i{+}1,\ldots,j{-}1,j\} $,
$B = \{j,j{+}1,\ldots,i{-}1,i\}$, and $ |i-j| \geq 2$.
Here $\{A,B\}$ runs over all diagonals in the $m$-gon, and
indices are understood modulo $m$. We conjecture that
\begin{equation}
\label{PCm}
 P_{C_m} \,\,\, = \,\,\, {\rm ST}_{C_m}  .
\end{equation}
We computed the ideal $P_{C_m}$ in
{\tt Singular}\footnote{{\tt www.singular.uni-kl.de/}}
for small $m$. The following table lists the results:
$$
\begin{array}{l | c c c c c c}
m & \quad 3\quad  & \,\quad 4\,\quad  & \,\quad 5 \,\quad & \,\quad 6 \,\quad & \,\quad 7 \,\quad & \,\quad 8 \,\quad \\\hline
\textrm{dimension } $d$ & 6 & 8 & 10 & 12 & 14 & 16 \\
\textrm{degree} & 1 & 9 & 57 & 312 & 1578 &  7599 \\
\textrm{ML-degree} & 1 & 5 & 17  & 49  & 129 & 321 \\
\textrm{minimal generators (degree:number)} & 0 & 3\!:\!2& 3\!:\!15 &3\!:\!63 &
3\!:\!196 & 3\!:\!504
 \end{array}
 $$
In all cases in this table, the minimal generators consist of
cubics only, which is consistent with the conjecture (\ref{PCm}).
For the degree of the Gaussian $m$-cycle we  conjecture the following formula.
\begin{conj}
The degree of the projective variety $V(P_{C_m}) $
associated with the $m$-cycle equals
$$\frac{m+2}{4}\binom{2m}{m}-3\cdot 2^{2m-3}.$$
\end{conj}

Regarding Question 2, the following formula
was conjectured in \cite[\S 7.4]{DSS}:
\begin{equation}
\label{MLconj}
 \hbox{\rm ML-degree}(C_m) \,\,\, = \,\,\,
(m-3) \cdot 2^{m-2}+1, \qquad {\rm for} \,\,\,m\geq 3.
\end{equation}
This quantity is an algebraic complexity measure for the
following matrix completion problem. Given real numbers
$x_i $ between $-1$ and $+1$, fill up the partially specified
symmetric $m {\times} m$-matrix
\begin{equation}
\label{partialmatrix} \begin{pmatrix}
    1 & x_1 &   ?   & ?     &  \cdots & ? & x_m \\
x_1 & 1     & x_2 & ?     &  \cdots & ? &  ?  \\
 ?    & x_2 & 1     & x_3 &  ? & \ddots  & ? \\
 ?   &   ?    & x_3 &   1   & x_4 & \ddots & \vdots \\
 \vdots & \vdots & & \ddots & \ddots & \ddots & ?  \\
   ?  & ?       & ? &     ?   & x_{m-2} & 1 & x_{m-1} \\
   x_m & ?   & ?  &    ?   &   ?     & x_{m-1} & 1
  \end{pmatrix}
 \end{equation}
to make it positive definite.
We seek the unique fill-up that maximizes the determinant.
The solution to this convex optimization problem
is an algebraic function of $x_1,x_2,\ldots,x_m$
whose degree equals $\hbox{\rm ML-degree}(C_m)$.
We do not know how to prove
(\ref{MLconj}) for $m \geq 9$. Even
the following lemma is not easy.

\begin{lem}
\label{MLispositive}
The ML-degree of the cycle $C_m$ is strictly larger than $1$ for $m \geq 4$.
\end{lem}

\begin{proof}[Sketch of Proof]
We consider the special case of (\ref{partialmatrix})
when all of the parameters are equal:
\begin{equation}
\label{allequal}
 x \,:= \, x_1 = x_2 = \cdots = x_m .
\end{equation}
Since the logarithm of the determinant is a concave function,
the solution to our optimization problem is fixed under
the symmetric group of the $m$-gon, i.e., it is a
symmetric {\em circulant matrix} $\Sigma_m$.
Hence there are only  $\lfloor \frac{m-2}{2} \rfloor $
distinct values for the question marks in
(\ref{partialmatrix}), one for each of the symmetry class of
long diagonals in the $m$-gon.  We denote these unknowns by
$\,s_1,s_2,\ldots,s_{\lfloor \frac{m-2}{2}\rfloor}\,$
where $s_i$ is the unknown on the $i$-th circular off-diagonal.
For instance, for $m=7$, the circulant matrix we seek has two unknown
entries $s_1$ and $s_2$, and it looks like this:
$$  \Sigma_7 \quad = \quad
\begin{pmatrix}
 1 & x & s_1 & s_2 & s_2 & s_1 & x \\
 x & 1 & x & s_1 & s_2 & s_2 & s_1 \\
 s_1 & x & 1 & x & s_1 & s_2 & s_2  \\
 s_2 & s_1 & x & 1 & x & s_1 & s_2  \\
 s_2 & s_2 & s_1 & x & 1 & x & s_1  \\
 s_1 & s_2 & s_2 & s_1 & x & 1 & x  \\
 x & s_1 & s_2 & s_2 & s_1 & x & 1
\end{pmatrix}
$$
The key observation is that the
determinant of the circular symmetric
matrix $\Sigma_m$ factors into a product
of $m$ linear factors with real coefficients,
one for each $m$th root of unity. For example,
$$ {\rm det}(\Sigma_7) \,\,\,= \,\,
  \prod_{w:w^7=1} \! \bigl(\,
1 + (w+w^6)\cdot x + (w^2+w^5)\cdot s_1 + (w^3+w^4) \cdot s_2 \,\bigr).  $$
Thus, for fixed $x$, our problem is to maximize a product of linear forms.
By analyzing the critical equations, obtained by taking
logarithmic derivatives of ${\rm det}(\Sigma_m)$, we can show that
the optimal solution $\,(\hat s_1, \hat s_2, \ldots,
\hat s_{\lfloor \frac{m-2}{2}\rfloor})\,$
is not a rational function in $x$.
For example, when $m=7$, the solution $(\hat s_1, \hat s_2)$
is an algebraic function of degree $3$ in $x$. Its explicit
representation is
$$
   \hat s_1  = \frac{x^2+\hat s_2 x-\hat s_2^2 - \hat s_2}{1-x}
\qquad \hbox{and} \qquad
         \hat s_2^3 + (1-2 x) \hat s_2^2 + (-x^2+x-1) \hat s_2 + x^3 = 0
$$
A detailed proof, for arbitrary $m$, will
appear in the PhD dissertation of the second author.
\end{proof}

We now come to our third problem, namely to giving
an algebraic description of the cone of
sufficient statistics, denoted
$\,{\mathcal C}_m  := \mathcal{C}_{C_m} $.
This is a full-dimensional open convex cone  in $\R^{2m}$.
The coordinates on $\R^{2m}$ are $s_{11}, s_{22}, \ldots, s_{mm}$
and $x_1=s_{12}$, $x_2 = s_{23}, \ldots, x_m = s_{m1}$.
We consider
$$\,\mathcal{C}'_m \,\, := \,\, \mathcal{C}_m \,\cap \,
\bigl\{ s_{11} = s_{22} = \cdots = s_{mm} = 1 \bigr\} . $$
This is a full-dimensional
open bounded spectrahedron in $\R^m$. It consists of all
 $(x_1,\ldots,x_m)$ such that
(\ref{partialmatrix}) can be filled up to a
positive definite matrix. The $2 \times 2$-minors
of (\ref{partialmatrix}) imply that $\mathcal{C}'_m$ lies
in the cube $(-1,1)^m = \{|x_i| < 1\}$. The issue
is to identify further constraints.
We note that any description of the
$m$-dimensional spectrahedron $\mathcal{C}'_m$ leads to a description of the
$2m$-dimensional cone $ \mathcal{C}_m$ because
a vector $s \in \R^{2m}$ lies in $ \mathcal{C}_m$
if and only if the vector $x \in \R^m$ with
the following coordinates lies in $\mathcal{C}'_m$:
\begin{equation}
\label{eq:x_i}
 x_i \,\,\,= \,\, \frac{s_{ij}}{\sqrt{s_{ii} s_{jj}}}  \qquad
\hbox{for}\,\, i=1,2,\ldots,m
\end{equation}
\cite{BJT} gave a
beautiful polyhedral description of the spectrahedron $\mathcal{C}'_m$.
The idea is to replace each $x_i$ by its arc-cosine,
that is, to substitute $x_i = {\rm cos}(\phi_i)$ into (\ref{partialmatrix}).
Remarkably, the image of the spectrahedron $\mathcal{C}'_m$ under this
transformation is a convex polytope. Explicit linear inequalities
in the angle coordinates $\phi_i$ describing the facets
of this polytope are given in~\cite{BJT}.

To answer Question 3, we  take the cosine-image
of any of these facets and compute its Zariski closure.
This leads to the following trigonometry problem.
Determine the unique (up to scaling) irreducible
polynomial  $\Gamma'_m$ which is
obtained by rationalizing the equation
\begin{equation}
\label{rationalizing}
 x_1 \,\, = \,\, {\rm cos} \biggl(  \sum_{i=2}^m {\rm arccos}(x_i) \biggr).
\end{equation}
We call $\Gamma'_m$ the {\em $m$-th cycle polynomial}. Interestingly,
$\Gamma'_m$ is invariant under all permutations of the $m$ variables
$x_1,x_2,\ldots,x_m$.
We also define the {\em homogeneous $m$-th cycle polynomial}
$\Gamma_m$
to be the numerator of the image of $\Gamma'_m$ under
the substitution (\ref{eq:x_i}). The first
cycle polynomials arise for $m=3$:
$$
\Gamma'_3 \,\, = \,\,
{\rm det} \begin{pmatrix}
1 & x_1 & x_3 \\
x_1 & 1 & x_2 \\
x_3 & x_2 & 1 \end{pmatrix}
\qquad \hbox{and} \qquad
\Gamma_3 \,\, = \,\,
{\rm det} \begin{pmatrix}
s_{11} & s_{12} & s_{13} \\
s_{12} & s_{22} & s_{23} \\
s_{13} & s_{23} & s_{33}
\end{pmatrix}.
$$
The polyhedral characterization of $\mathcal{C}_m$ given in \cite{BJT}
translates into the following theorem.

\begin{thm}
The Zariski closure of  the boundary of the cone $\mathcal{C}_m$,
$m \geq 4$,
is defined by the polynomial
$$ H_{C_m}(s_{ij}) \quad = \quad
\Gamma_m(s_{ij}) \cdot
(s_{11} s_{22} - s_{12}^2) \cdot
(s_{22} s_{33} - s_{23}^2) \,\cdots \,
(s_{mm} s_{11} - s_{1m}^2).
$$
\end{thm}

To compute the cycle polynomial $\Gamma'_m$, we
iteratively apply the  sum formula for the cosine,
$$ {\rm cos}(a+b) \,\,=\,\, {\rm cos}(a) \cdot {\rm cos}(b) - {\rm sin}(a) \cdot {\rm sin}(b) ,$$
and we then use the following relation to write (\ref{rationalizing})
as an algebraic expression in $x_1,\ldots,x_n$:
$$ {\rm sin}\bigl({\rm arccos}(x_i) \bigr) \,\, = \,\, \sqrt{1-x_i^2 \,\,}. $$
Finally, we eliminate the square roots (e.g. by using resultants)
to get the polynomial $\Gamma'_m$.

For example, the cycle polynomial for the square
$(m=4)$ has degree $6$ and has $19$ terms:
$$
\Gamma'_4 \,\, = \,\,
4 \! \sum_{i < j < k} \!\! x_i^2 x_j^2 x_k^2 \,-\, 4 x_1x_2 x_3 x_4 \! \sum_i x_i^2 \,+\, \sum_i x_i^4
\,-\, 2 \sum_{i <j} x_i^2 x_j^2 \, +\, 8 x_1 x_2 x_3 x_4.
$$
By substituting (\ref{eq:x_i}) into this expression and taking the
numerator, we obtain the homogeneous cycle polynomial $\Gamma_4$
which has degree $8$.
Here is a table summarizing what we know about
the expansions of these cycle polynomials. Note that
$\Gamma'_m$ and $\Gamma_m$ have different degrees but
the same number of terms.
$$
\begin{array}{l | c c c c c c c c c}
m &                                                     \quad 3\quad  & \quad 4\quad    & \quad 5\quad   & \quad 6\quad  & \quad 7\quad  & \quad 8\quad  & \quad 9\quad  & \quad 10\quad  & \quad 11\quad  \\ \hline
{\rm degree}(\Gamma'_m) &     3 &  6 & 15 & 30 & 70 & 140 & 315 & 630 & 1260   \\
{\rm degree}(\Gamma_m) &      3 &  8 & 20 & 48 & 112 & 256 & 576 & 1280 & 2816  \\
\# \hbox{of terms} & 5 & 19 & 339  & 19449 & ?& ? & ? & ? & ? \\
\end{array}
$$
The degree of the $m$-th cycle polynomial $\Gamma'_m$ grows roughly like
$2^m$, but we do not know an exact formula. However, for the homogeneous
cycle polynomial $\Gamma_m$ we predict the following behavior.

\begin{conj}
The degree of the homogeneous $m$-th cycle polynomial $\Gamma_m$ equals
$\,m\cdot 2^{m-3}$.
\end{conj}

There is another way of defining and computing the
cycle polynomial $\Gamma_m$, without any reference to
trigonometry or semidefinite programming.
Consider the prime ideal generated by the
$3 \times 3$-minors of
the generic symmetric $m {\times} m$-matrix $\Sigma = (s_{ij})$.
Then $\langle \Gamma_m \rangle$ is the principal ideal obtained
by eliminating all unknowns $s_{ij}$ with $|i-j| \geq 2$.
Thus, geometrically, vanishing of  the
homogeneous polynomial $\Gamma_m$ characterizes
partial matrices on the $m$-cycle $C_m$ that can be completed to
a matrix of rank $\leq 2$. Similarly, vanishing
of $\Gamma'_m$ characterizes partial matrices (\ref{partialmatrix})
that can be completed to rank $\leq 2$.

Independently of the work of \cite{BJT},
a solution to the problem of characterizing the cone $\mathcal{C}_m$
appeared in the same year in the statistics literature, namely by \cite{Buh}.
For statisticians, the cone $\mathcal{C}_m$ is the set of partial sample
covariance matrices on the $m$-cycle for which the MLE exists.

\subsection{Small graphs, suspensions and wheels}
\label{subsec_wheels}

We next examine  Questions 1, 2 and 3 for all graphs with at most
five vertices. In this analysis
we can restrict ourselves to connected graphs only.
Indeed, if $G$ is the disjoint union of two graphs $G_1$ and $G_2$
then the prime ideal $P_G$ is obtained from $P_{G_1}$ and $P_{G_2}$ as in
(\ref{cliquesumideal}) with $c=0$,  the ML-degrees multiply by
Lemma \ref{MLdegree_cliquesum},
and the two dual cones both decompose as direct products:
$$   \mathcal{C}_G \,\, = \,\, \mathcal{C}_{G_1} \times \mathcal{C}_{G_2}
\quad \hbox{and} \quad
\mathcal{K}_G \,\, = \,\, \mathcal{K}_{G_1} \times \mathcal{K}_{G_2}. $$
Chordal graphs were dealt with in Section \ref{chordal_graphs}.
We now consider
connected non-chordal graphs with $m\leq 5$ vertices.
There are seven such graphs, and in Table \ref{table_smallgraphs} we
summarize our findings for these
 seven graphs. In the first two rows of Table \ref{table_smallgraphs}
we find the
$4$-cycle and the $5$-cycle which were discussed
in Subsection \ref{subsec_cycle}. As an illustration we examine
in detail the graph in the second-to-last row of Table
\ref{table_smallgraphs}.

\begin{table}
\caption{Our three guiding questions
for all non-chordal graphs with $m\leq 5$ vertices. Column 4 reports the degrees of the minimal generators together with the number of occurrence (degree:number). The last column lists the degrees of the irreducible factors of the polynomial $H_G$ that defines
the Zariski closure of the boundary of $\mathcal{C}_G$. For each
factor we report in lowercase the rank of the concentration matrices
defining its dual irreducible component in
the boundary of $\mathcal{K}_G$.}
\label{table_smallgraphs}
\smallskip
\begin{center}
\begin{tabular}{l | c | c | c | c | c }
Graph $G$ & dim $d$ & deg $P_G$ & mingens $P_G$ & ML-deg & deg $H_G$\\ \hline
\parbox{1.7cm}{\rule[-1mm]{0pt}{11mm}\includegraphics[width=0.98cm]{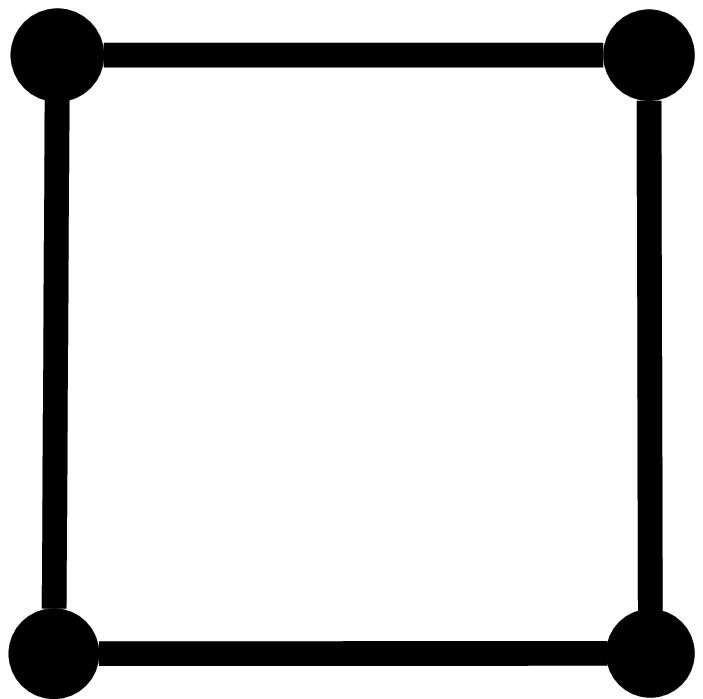}} & 8 & 9 & 3:2 &  5 & $4\cdot 2_1 + 8_2$\\
\parbox{1.7cm}{\rule[-1mm]{0pt}{11mm}\includegraphics[width=1.1cm]{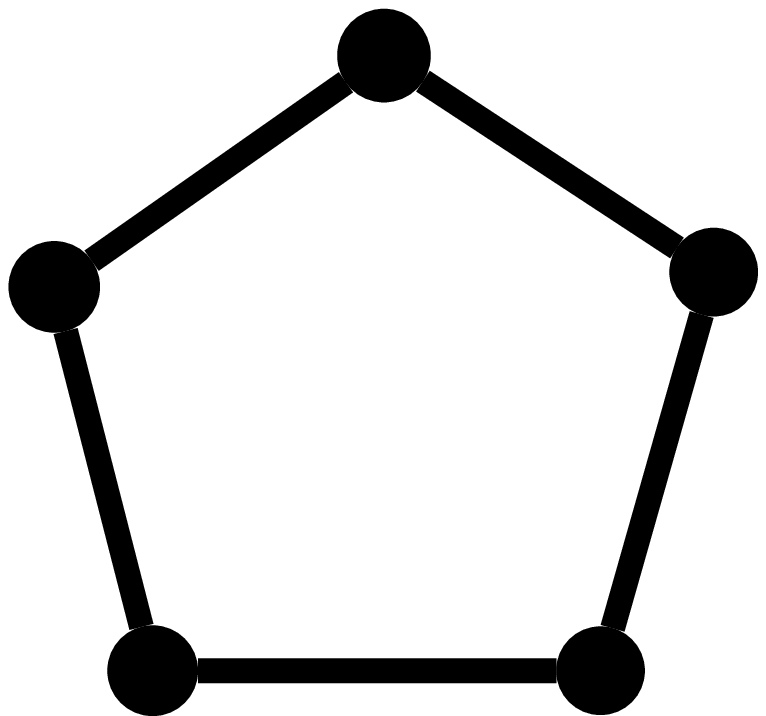}} & 10 & 57 & 3:15 & 17 & $5\cdot 2_1 + 20_3$\\
\parbox{1.7cm}{\rule[-1mm]{0pt}{11mm}\includegraphics[width=1.8cm]{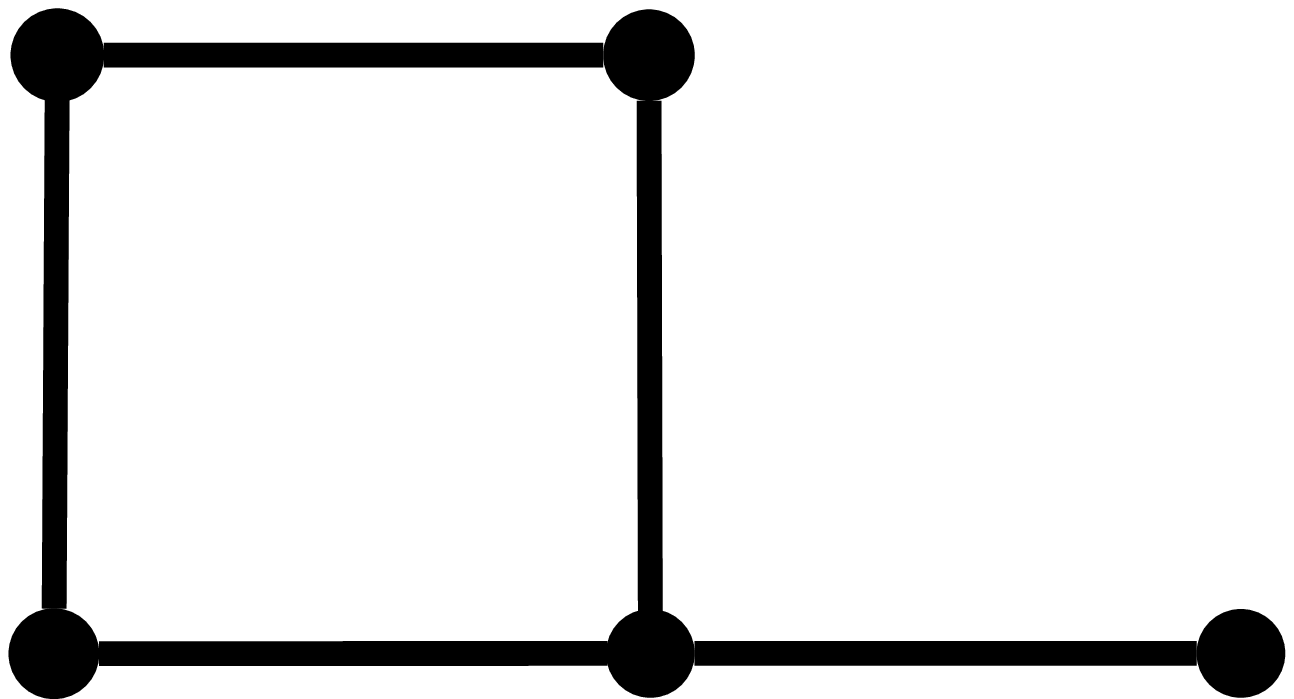}} & 10 & 30 & 2:6, 3:4 & 5 & $5\cdot 2_1 + 8_2$\\
\parbox{1.7cm}{\rule[-1mm]{0pt}{11mm}\includegraphics[width=1.65cm]{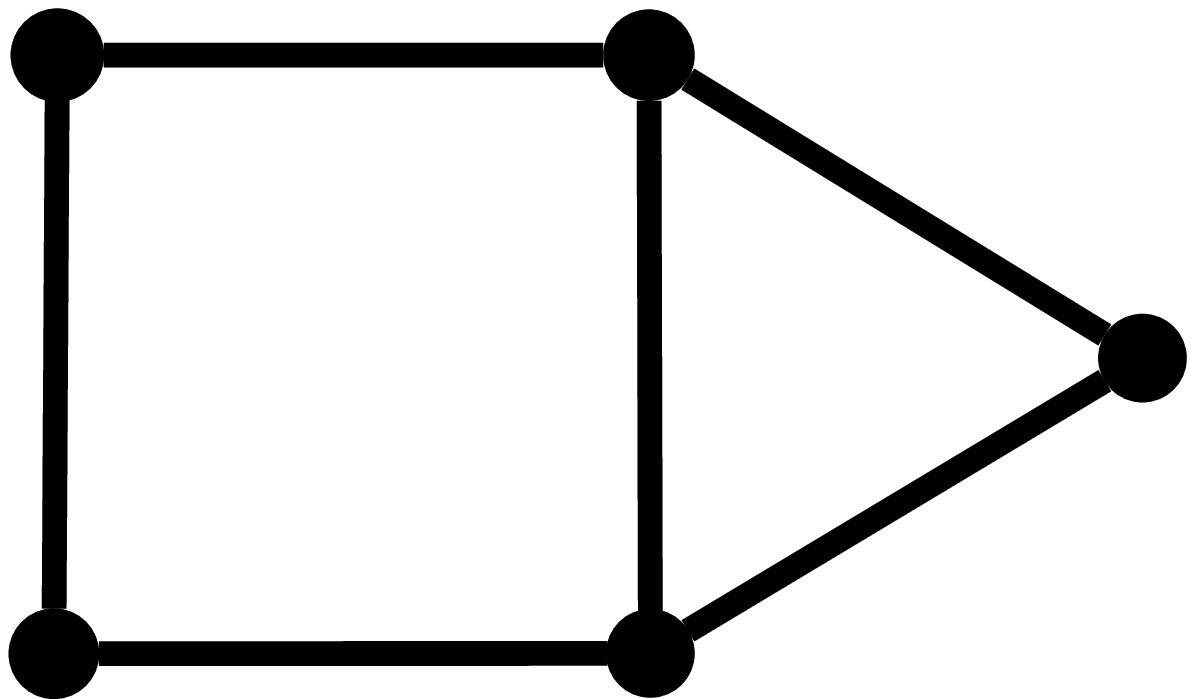}} & 11 & 31 & 3:10 & 5 &  $3\cdot 2_1+3_1+8_2$\\
\parbox{1.7cm}{\rule[-1mm]{0pt}{11mm}\includegraphics[width=1.65cm]{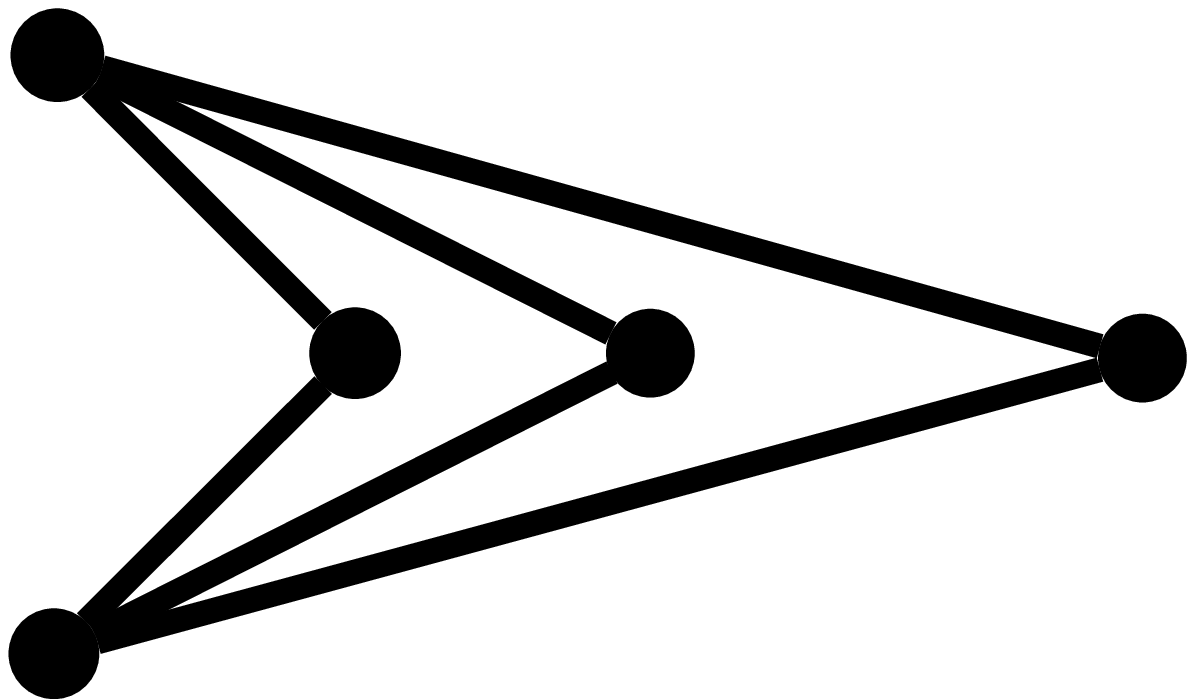}} & 11 & 56 & 3:7, 4:1 & 7 &  $6\cdot 2_1+3\cdot 8_2$\\
\parbox{1.7cm}{\rule[-1mm]{0pt}{11mm}\includegraphics[width=1.65cm]{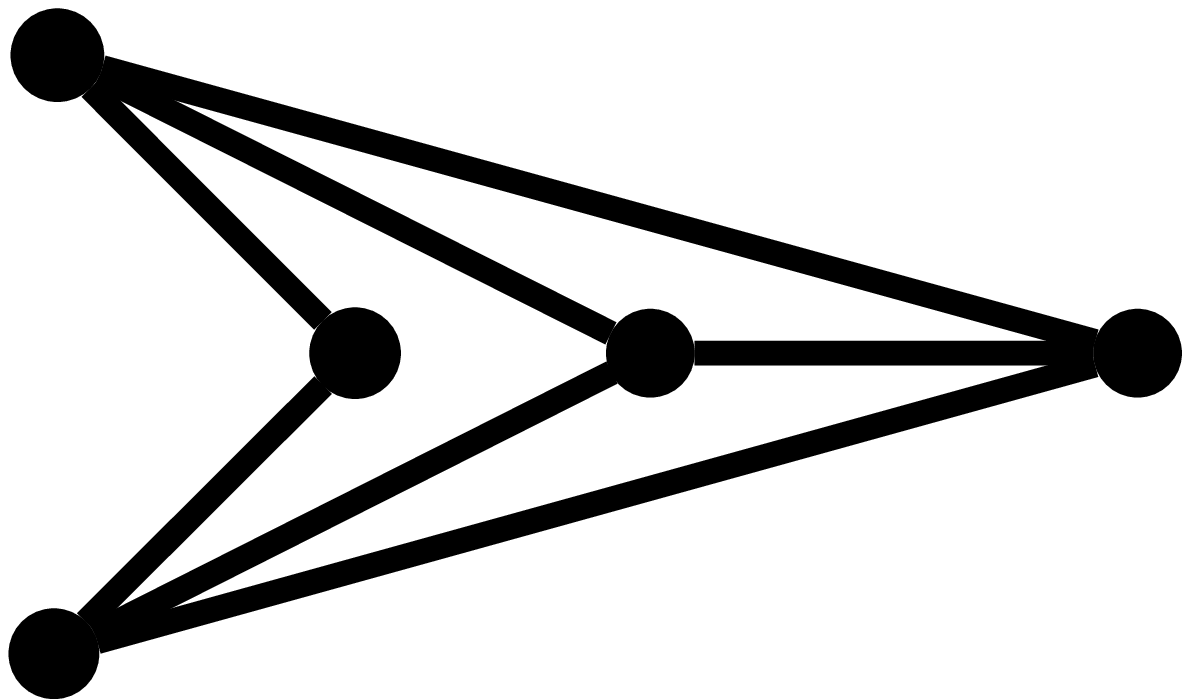}} & 12 & 24 & 3:4, 4:1 & 5 &  $2\cdot 2_1 + 2\cdot 3_1 + 10_2$\\
\parbox{1.7cm}{\rule[-1mm]{0pt}{11mm}\includegraphics[width=1.65cm]{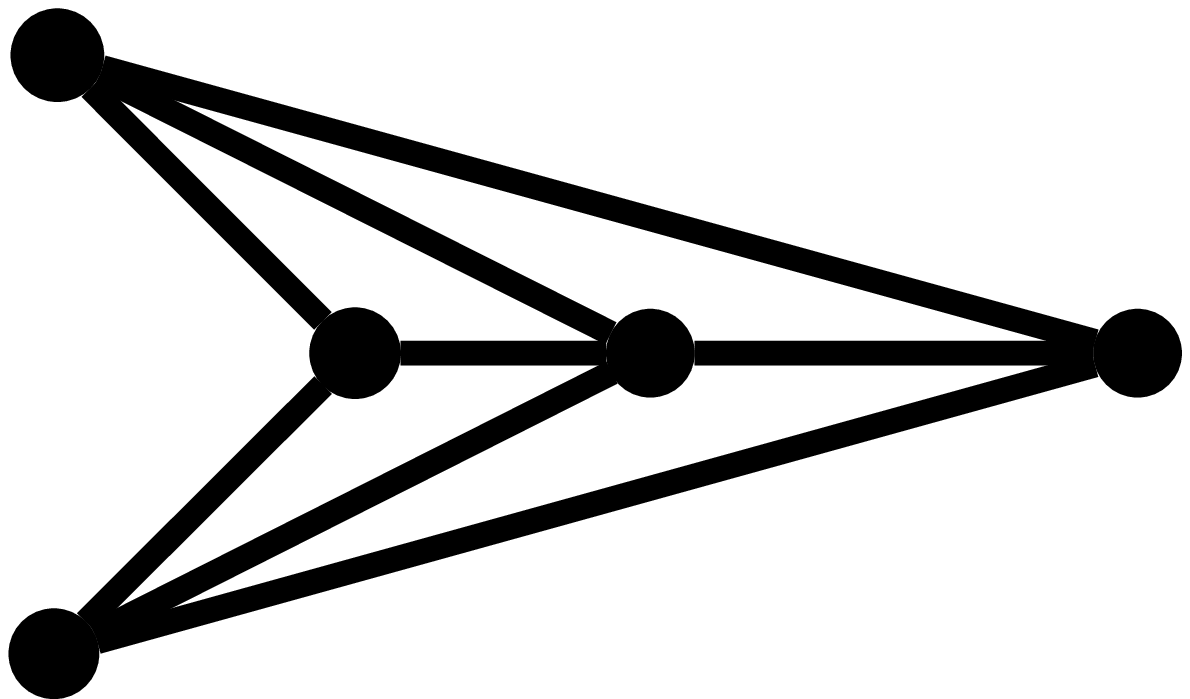}} & 13 & 16 & 4:2 & 5 & $4\cdot 3_1 + 12_2$
\end{tabular}
\end{center}
\end{table}

\begin{figure*}[!h]
\centering
\includegraphics[width=3cm]{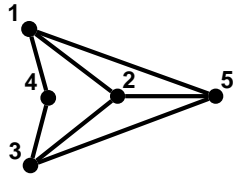}
\caption{A Gaussian graphical model on five vertices
and seven edges having
dimension $d=12$.}
\label{ex_graph5}
\end{figure*}

\begin{example}
\label{ex_5graph}
The graph in Fig.~\ref{ex_graph5}
defines the Gaussian graphical model with concentration matrix
$$ K \quad = \quad
\begin{pmatrix}
\lambda_1 & \lambda_6 & 0 & \lambda_9 & \lambda_{10} \\
\lambda_6 & \lambda_2 & \lambda_7 & 0 & \lambda_{11} \\
0 & \lambda_7 & \lambda_3 & \lambda_8 & \lambda_{12} \\
\lambda_9 & 0 & \lambda_8 & \lambda_4 & 0 \\
\lambda_{10} & \lambda_{11} & \lambda_{12} & 0 & \lambda_5
\end{pmatrix}.
$$
We wish to describe the boundary of the cone $\mathcal{C}_G$
by identifying the irreducible factors in its defining polynomial
$H_G$. We first use the {\tt Matlab} software
{\tt CVX}\footnote{{\tt www.stanford.edu/$\sim$boyd/cvx/}}, which
is specialized in convex optimization, to find the ranks of
all concentration matrices $K$ that are extreme rays
in the boundary of $\mathcal{K}_G$. Using {\tt CVX}, we
maximize random linear functions over the compact spectrahedron
$\, \overline{\mathcal{K}_G} \cap \{ {\rm trace}(K) = 1\}$,
and we record the ranks of the optimal matrices. We found the
possible matrix ranks to be $1$ and $2$, which agrees
with the constraints $\,2 \leq p \leq 3\,$
seen in (\ref{pataki}) for generic subspaces
$\mathcal{L}$ with $m=5$ and $d=12$.

We next ran the software {\tt Singular} to
compute the minimal primes of the ideals of $p\times p$-minors of $K$
for $p=2$ and $p=3$, and thereafter we computed their dual ideals
in $\mathbb{R}[t_1,t_2,\ldots ,t_{12}]$ using
{\tt Macaulay2}. The latter step was done
using the procedure with Jacobian matrices described in
Subsection \ref{subsec_generic}. We only retained dual
ideals that are principal. Their generators
are the candidates for factors of $H_G$.

The variety of rank one matrices $K$ has four irreducible components.
Two of those components correspond to the edges
$(3,4)$ and $(1,4)$ in Fig.~\ref{ex_graph5}. Their dual ideals are
generated by the quadrics
$$p_1=4t_3t_4-t_8^2 \qquad \hbox{and} \qquad p_2=4t_1t_4-t_9^2.$$
The other two irreducible components of the variety of rank one concentration matrices correspond to the 3-cycles $(1,2,5)$ and $(2,3,5)$ in the graph. Their dual ideals are generated by the cubics
$$ p_3= 4t_1t_2t_5-t_5t_6^2-t_2t_{10}^2+t_6t_{10}t_{11}-t_1t_{11}^2 \quad
\, \hbox{and} \quad \,
p_4=4t_2t_3t_5-t_5t_7^2-t_3t_{11}^2+t_7t_{11}t_{12}-t_2t_{12}^2.$$
The variety of rank two  matrices $K$ has two irreducible components.
One corresponds to the chordless 4-cycle $(1,2,3,4)$ in the graph and
its dual ideal is generated by $p_5=\Gamma_4$, which is of degree $8$.
The other component consists of rank two matrices $K$
for which rows $2$ and $5$ are linearly dependent.
The polynomial $p_6$ that defines the dual ideal consists of 175 terms and has degree 10.

The polynomial $H_G$ is the product of those principal
generators $p_i$ whose hypersurface meets
$\,\partial \mathcal{C}_G$.  We again used {\tt CVX}
to check which of the six components
actually contribute extreme rays in
$\,\partial \mathcal{K}_G$.
We found that only one of the six components to be missing, namely that
corresponding to the chordless 4-cycle $(1,2,3,4)$.
This means that $p_5$ is not a factor of $H_G$, and we conclude
\begin{equation}
H_G \,\,=\,\,p_1 p_2 p_3 p_4 p_6 \quad \hbox{and}
\quad deg(H_G)\,\,= \,\,2\cdot 2_1 \,+\, 2\cdot 3_1 \,+\, 10_2.
\end{equation}
Concerning Question 1 we note that the ideal
$P_G$ is minimally generated by the four $3 {\times} 3$-minors
of $\Sigma_{1235,134}$ the determinant of $\Sigma_{1245,2345}$,
and for Question 2 we note that the ML degree is five
because the MLE can be derived from the MLE of the $4$-cycle
obtained by contracting the edge $(2,5)$.
\qed
\end{example}

The graph in the last row of  Table \ref{table_smallgraphs} is the
wheel $W_4$. It is obtained from the cycle
$C_4$ in the first row by
connecting all four vertices to a new fifth vertex.
 We see in   Table \ref{table_smallgraphs} that the ML degree $5$
is the same for both graphs, the two cubic
generators of $P_{C_4}$ correspond to the
two quartic generators of $P_{W_4}$, and there
is a similar correspondence between the irreducible
factors of the dual polynomials $H_{C_4}$ and $H_{W_4}$.
In the remainder of this section we shall offer an explanation
for these observations.

Let $G=(V,E)$ be an undirected graph and $G^*=(V^*,E^*)$ its
{\em suspension graph} with an additional completely connected vertex $0$.
The graph $G^*$ has vertex set $V^*=V\cup\{0\}$ and edge set
$E^*=E\,\cup\,\{(0,v)\mid v\in V\}$.
The {\em $m$-wheel} $W_m$ is the
suspension graph of the $m$-cycle $C_m$; in symbols, $W_m = (C_m)^*$.
We shall compare the Gaussian graphical models
for the graph $G$ and its suspension graph $G^*$.

\begin{thm}
The ML degree of a Gaussian graphical model with underlying graph
$G$ equals the ML degree of a Gaussian graphical model
whose underlying graph is the suspension graph $G^*$.
\end{thm}

\begin{proof}
Let $V = [m]$ and let $S^* \in \s^{m+1}_{\succ 0}$
be a sample covariance matrix on $G^*$,
where the first row and column correspond to the additional vertex $0$.
We denote by $S'$ the lower right $m\times m$ submatrix of $S^*$
corresponding to the vertex set $V$ and by $S$ the Schur complement
of $S^*$ at $S^*_{00}$:
\begin{equation}
\label{SandSprime}
 S \,\,\, := \,\,\,
S'-\frac{1}{S^*_{00}}(S^*_{01},\dots ,S^*_{0m})^T(S^*_{01},\dots ,S^*_{0m}).
\end{equation}
Then $S \in \s^m_{\succ 0}$ is a sample covariance matrix on $G$.
Let $\hat{\Sigma}$ be the MLE for $S$ on the graph $G$.
We claim that the MLE $\hat{\Sigma}^*$ for $S^*$ on
the suspension graph $G^*$ is given by
\[
\hat{\Sigma}^* \,\,= \,\,\left[ \begin{tabular}{c|c c c}
$S^*_{00}$ & $S^*_{01}$ & $\cdots$ & $S^*_{0m}$ \\ \hline $S^*_{01}$ & & & \\ \vdots & \multicolumn{3}{c}{$\hat{\Sigma}+S'-S$} \\ $S^*_{0m}$ & & &
\end{tabular}\right].
\]
Clearly, $\hat{\Sigma}^*$ is positive definite and satisfies $\hat{\Sigma}^*_{ij}=S^*_{ij}$ for all $(i,j)\in E^*$. The inverse of the covariance
 matrix $\,\hat{\Sigma}^*\,$ can be computed by using
the inversion formula based on Schur complements:
\[
(\hat{\Sigma}^*)^{-1} \,\,= \,\,\left[ \begin{tabular}{c|c}
$\frac{1}{S^*_{00}}+ (S^*_{01},\dots ,S^*_{0m})  (\hat{\Sigma})^{-1}(S^*_{01},
\dots ,S^*_{0m})^T$ & $\frac{1}{S^*_{00}}(\hat{\Sigma})^{-1} (S^*_{01},\dots ,
S^*_{0m})^T$  \\ \hline
$\frac{1}{S^*_{00}}(S^*_{01},\dots ,S^*_{0m})
(\hat{\Sigma})^{-1}$ & $ (\hat{\Sigma})^{-1} $
\end{tabular}\right].
\]
Since the lower right block equals $(\hat{\Sigma})^{-1}$,
its entries are indeed zero in all positions $(i,j)\notin E^*$.

We have shown that the MLE $\hat{\Sigma}^*$ is a
rational function of the MLE $\hat{\Sigma}$. This shows
$$\textrm{ML-degree}(G^*)\leq \textrm{ ML-degree}(G).$$
The reverse inequality is also true since we can
compute the MLE on $G$
for any  $S \in S^m_{\succ 0}$
by computing the MLE on $G^*$ for its extension $S^* \in \s^{m+1}_{\succ 0}$
with $S^*_{00} = 1$ and $S^*_{0j} = 0$ for $j \in [m]$.
\end{proof}

We next address the question of how the boundary of
the cone $\mathcal{C}_{G^*}$ can be expressed
 in terms of the boundary of $\mathcal{C}_G$.
We use coordinates $t_{ij}$ for both $\s^m$ and its subspace
$\R^{E}$, and we use the coordinates $u_{ij}$ for
both $\s^{m+1}$ and its subspace $\R^{E^*}$. The Schur complement
(\ref{SandSprime})  defines a rational map
from  $\s^{m+1}$ to $\s^{m}$
which restricts to a rational map from
$\R^{E^*}$ to $\R^E$. The formula is
\begin{equation}
\label{fromTtoU}
t_{ij} \,\, = \,\, u_{ij} - \frac{u_{0i} u_{0j}}{u_{00}}
\qquad \hbox{for} \,\,\,1 \leq i \leq j \leq m.
\end{equation}
A partial matrix $(u_{ij})$ on $G^*$ can be completed
to a positive definite $(m+1) {\times} (m+1)$-matrix
if and only if the partial matrix $(t_{ij})$
on $G$ given by this formula can be completed
to a positive definite $m {\times} m$-matrix.
The rational map (\ref{fromTtoU}) takes the  boundary of the cone
$\mathcal{C}_{G^*}$ onto the boundary of the cone $\mathcal{C}_{G}$.
For our algebraic question, we can derive the following conclusion:

\begin{prop}
The polynomial $H_{G^*}(u_{ij})$ equals the
numerator of the Laurent polynomial obtained
from $H_G(t_{ij})$ by the substitution (\ref{fromTtoU}),
and the same holds for each irreducible factor.
\end{prop}

\begin{example}
The polynomial
$\,H_{W_4}(u_{00},u_{01},u_{02},u_{03},u_{04},
u_{11},u_{22},u_{33},u_{44},u_{12},u_{23},u_{34},u_{14})\,$
for the $4$-wheel $W_4$ has as its main factor an irreducible
polynomial of degree $12$ which is the sum of $813$ terms.
It is obtained from the homogeneous cycle polynomial $\Gamma_4$
by the  substitution (\ref{fromTtoU}).
Recall that $\Gamma_4(t_{11},t_{22},t_{33},t_{44},t_{12},t_{23},t_{34},t_{14})\,$ has only degree $8$ and
is the sum of $19$ terms. \qed
\end{example}

We briefly discuss an issue raised by Question 1, namely,
how to construct the prime ideal $P_{G^*}$ from the
prime ideal $P_G$. Again, we can use the
transformation (\ref{fromTtoU}) to turn every generator
of $P_G$ into a Laurent polynomial whose numerator
lies in  $P_{G^*}$. However, the resulting polynomials
will usually not suffice to generate $P_{G^*}$. This happens
already for the $5$-cycle $G = C_5$ and the
$5$-wheel $G^* = W_5$. The ideal $P_{C_5}$ is generated by
$15$ linearly independent cubics arising as
$3 {\times} 3$-minors of the matrices
$\Sigma_{132,1345},\,\Sigma_{243,2451},\,\Sigma_{354,3512},\,
\Sigma_{415,4123}$ and $\Sigma_{521,5234}$, while
$P_{W_5}$ is generated by
$20$ linearly independent quartics arising as
$4 {\times} 4$-minors of
$\Sigma_{0132,01345},\,\Sigma_{0243,02451},\,\Sigma_{0354,03512},\,
\Sigma_{0415,04123}$ and $\Sigma_{0521,05234}$.
Here is a table that summarizes what we know about the
Gaussian wheels $W_m$:
$$
\begin{array}{l | c c c c }
m &\quad3\quad & \quad4\quad & \quad5\quad & \quad6\quad \\ \hline
\textrm{dimension } d & 10 & 13 & 16 & 19 \\
\textrm{degree} & 1 & 16 & 198 & 2264 \\
\textrm{ML-degree} & 1 & 5 & 17 & 49 \\
\textrm{minimal generators (degree:number)} & 0 & 4\!:\!2 & 4\!:\!20
& 4\!:\!108 \end{array}
$$

\section{Colored Gaussian graphical models} \label{colored_graphs}
We now add a graph coloring to the setup and study colored Gaussian
graphical models. These were introduced by
\cite{HL} who called them {\em RCON-models}.
In the underlying graph $G$, the vertices are colored with $p$
different colors and the edges are colored with $q$ different colors:
$$
\begin{array}{lllll}
V &\,\,=\,\,& V_1\sqcup V_2\sqcup\dots\sqcup V_p, &\quad &p\leq |V| \\
E &\,\,=\,\,& E_1\sqcup E_2\sqcup \dots\sqcup E_q, &\quad &q\leq |E|.
\end{array}
$$
We denote the uncolored graph by $G$ and the colored graph by $\mathcal{G}$.
In addition to the restrictions given by the missing edges in the graph,
the entries of the concentration matrix $K$ are now also restricted
by equating entries in $K$ according to the edge and vertex colorings.
To be precise, the linear space $\mathcal{L}$ of $\s^m$ associated with a
colored graph $\mathcal{G}$
on $m = |V|$ nodes is defined by the following linear equations:
\begin{itemize}
\item For any pair of nodes $\alpha, \beta$ that do not
form an edge we set $\,k_{\alpha \beta} = 0\,$  as before.
\item
For any pair of nodes $\alpha, \beta$ in a common color class $V_i$
we set $\,k_{\alpha\alpha}=k_{\beta \beta}$.
\item
For any pair of edges $(\alpha,\beta),\,(\gamma,\delta)$
in a common color class $E_j$ we set
$\,k_{\alpha \beta}=k_{\gamma \delta}$.
\end{itemize}
The dimension of the model $\mathcal{G}$ is $\,d = p+q$. We
note that, for any sample covariance matrix $S$,
$$
\pi_G(S) \, \in \,  \mathcal{C}_G \qquad \textrm{implies} \qquad \pi_{\mathcal{G}}(S) \, \in \,  \mathcal{C}_{\mathcal{G}}.
$$
Thus, introducing a graph coloring
on $G$ relaxes the question of existence of the MLE.

In this section we shall examine Questions 1-3
for various colorings $ \mathcal{G}$ of the $4$-cycle $G = C_4$.
We begin with an illustration of how colored Gaussian graphical models can be used
in statistical applications.

\begin{figure*}[!h]
\centering
 \includegraphics[width=3cm]{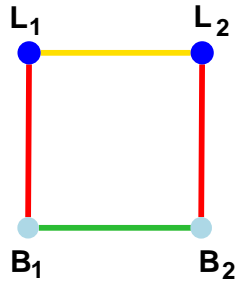}
\caption{Colored Gaussian graphical model for  Frets' heads: $L_i$, $B_i$ denote the length and breadth of the head of son $i$.}
\label{Frets_heads}
\end{figure*}

\begin{example}[Frets' heads]
\label{ex_Fret} We revisite the
 heredity study of head dimensions reported in \cite{MKB}
and known to statisticians as \emph{Frets' heads}.
The data reported in this study consists of the length and breadth
of the heads of $25$ pairs of first and second sons. Because of the
symmetry between the two sons, it makes sense to try to fit
the colored Gaussian graphical model given in Fig.~\ref{Frets_heads}.

This model has $d=5$ degrees of freedom and it consists of all
concentration matrices of the form
$$ K\,\,=\,\,\begin{pmatrix} \lambda_1 & \lambda_3 & 0 & \lambda_4 \\ \lambda_3 & \lambda_1 & \lambda_4 & 0\\ 0 & \lambda_4 & \lambda_2 & \lambda_5 \\ \lambda_4 & 0 & \lambda_5 & \lambda_2 \end{pmatrix}.$$
In Fig.~\ref{Frets_heads},
the first random variable is denoted $L_1$, the second
$L_2$, the third $B_2$, and the fourth $B_1$.
Given a sample covariance matrix $S=(s_{ij})$, the
five sufficient statistics for this model are
\begin{equation}
\label{5suffstat}
t_1=s_{11}+s_{22}, \quad t_2=s_{33}+s_{44},
 \quad t_3=2s_{12}, \quad t_4 = 2(s_{23}+s_{14}), \quad t_5=2s_{34}.
\end{equation}
The ideal of polynomials vanishing on $K^{-1} $ is
generated by four linear forms and one cubic in the $s_{ij}$:
\begin{eqnarray*}
P_{\mathcal{G}} &\,\,=\,\, &\langle \,
s_{11}-s_{22}\,,\,\,\,s_{33}-s_{44}\,,
\,\,\,s_{23}-s_{14}\,,\,\,\,s_{13}-s_{24}\,,
\\ && \,\,\,s_{23}^2 s_{24} - s_{24}^3-s_{22} s_{23} s_{34} + s_{12} s_{24} s_{34} - s_{12} s_{23} s_{44} + s_{22} s_{24} s_{44}  \,\rangle.
\end{eqnarray*}
Note that the four linear constraints on the sample covariance matrix
seen in $P_{\mathcal{G}}$ are also valid
constraints on the concentration matrix. Models with this
property were studied in general by \cite{J} and appear
under the name {\em RCOP-models} in \cite{HL}.

The data reported in the Frets' heads study
results in the following sufficient statistics:
$$t_1=188.256, \quad t_2=95.408,\quad t_3=133.750, \quad t_4=210.062, \quad t_5=67.302.$$
Substituting these values into (\ref{5suffstat})
and solving the equations on $V(P_\mathcal{G})$, we
find the MLE for this data:
$$\hat{\Sigma}=\begin{pmatrix}
\,94.1280\, &  \,66.8750\, &  \,44.3082\, &  \,52.5155\\
   66.8750  & 94.1280 &  52.5155 &  44.3082\\
   44.3082  & 52.5155 &  47.7040  & 33.6510\\
   52.5155 &  44.3082 &  33.6510  & 47.7040  \end{pmatrix}.
$$

Both the degree and
the ML-degree of this colored Gaussian graphical model is $3$,
which answers Questions 1 and 2. It remains to
describe the boundary of the cone $\mathcal{C}_{\mathcal{G}}$
and to determine its defining polynomial $H_{\mathcal{G}}$. The variety
of rank one concentration matrices has four irreducible components:
$$\langle k_2, k_4, k_5, k_1+k_3\rangle, \,\, \langle k_2, k_4, k_5, k_1-k_3\rangle, \,\, \langle k_1, k_3, k_4, k_2+k_5\rangle, \,\, \langle k_1, k_3, k_4, k_2-k_5\rangle.$$
These are points in $\PP^4$ and the ideals of their dual hyperplanes are
$\,\langle t_1-t_3\rangle, \, \langle t_1+t_3\rangle,
\, \langle t_2-t_5\rangle, \, \langle t_2+t_5\rangle$.
The variety of rank two concentration matrices is irreducible.
Its prime ideal and the dual thereof are
$$ \begin{matrix}
& \langle k_2 k_3+k_1 k_5, \,\,\,\, k_1 k_2-k_4^2+k_3 k_5, \,\,\,\, k_3 k_4^2+k_1^2 k_5-k_3^2 k_5\rangle & \\
& \langle 4t_2^2 t_3^2-4 t_1 t_2 t_4^2+t_4^4+8t_1t_2t_3t_5-4t_3t_4^2t_5+4t_1^2t_5^2\rangle. & \end{matrix} $$
This suggests that the hypersurface
$\,\partial \mathcal{C}_{\mathcal{G}}\,$ is given by
the polynomial
\begin{equation}
\label{eq_ex}
H_{\mathcal{G}}\,=\,
(t_1-t_3)(t_1+t_3)(t_2-t_5)(t_2+t_5)(4t_2^2 t_3^2-4 t_1 t_2 t_4^2+t_4^4+8t_1t_2t_3t_5-4t_3t_4^2t_5+4t_1^2t_5^2).
\end{equation}
Using {\tt CVX} as in Example \ref{ex_5graph}
we checked that all five factors meet
$\,\partial \mathcal{C}_{\mathcal{G}}$, so
(\ref{eq_ex}) is indeed correct. \qed
\end{example}

We performed a similar analysis for all colored Gaussian graphical models
on the 4-cycle $C_4$, which have the property that edges in the same color
class connect the same vertex color classes. The results are presented in
Table \ref{RCOR1}, \ref{RCOR2} and \ref{RCOP}. These models are of special interest
because they are invariant under rescaling of variables in the same vertex
color class. Such models were introduced and studied by \cite{HL}.
For models with an additional permutation property
(these are the {\em RCOP-models}), we explicitly list the polynomial
$\,H_{\mathcal{G}}$. A census of
 these models appears in Table \ref{RCOP}.

\begin{table}
\caption{Results on Questions 1, 2, and 3 for all colored
Gaussian graphical models with some symmetry restrictions
(namely, edges in the same color
class connect the same vertex color classes) on the 4-cycle.}
\label{RCOR1}
\medskip
\centering
\begin{tabular}{l | c | c | c | c | c | c }
Graph & $K$ &dim $d$ & degree & mingens $P_\mathcal{G}$
& ML-degree & deg $H_\mathcal{L}$ \\ \hline
\parbox{1cm}{\rule[-0.5mm]{0pt}{11mm}\includegraphics[width=1cm]{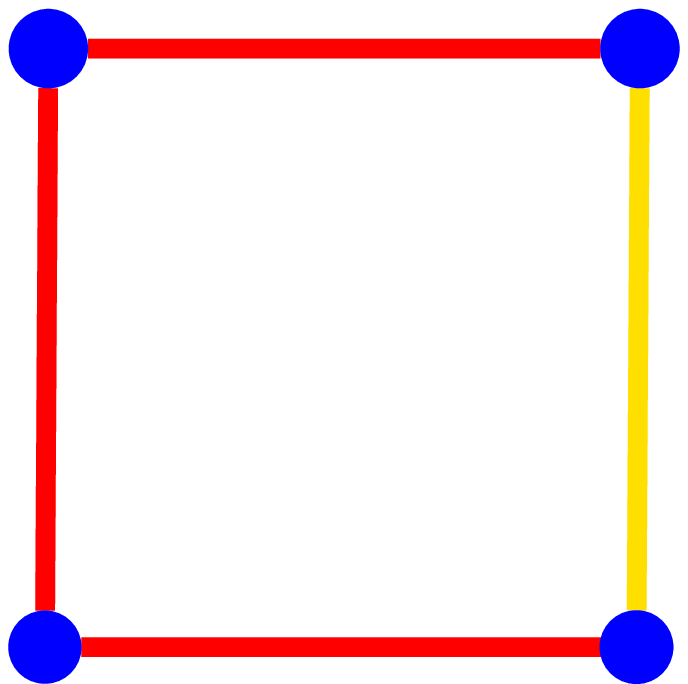}} & \rule[-7mm]{0pt}{15mm}$\begin{pmatrix} \lambda_1 & \lambda_2 & 0 & \lambda_2\\
\lambda_2 & \lambda_1 & \lambda_3 & 0\\
0 & \lambda_3 & \lambda_1 & \lambda_2\\
\lambda_2 & 0 & \lambda_2 & \lambda_1 \end{pmatrix}$ & 3 & 5 & 1:4, 2:5 & 5 & $2_2+2_3+2_3$\\
\parbox{1cm}{\rule[-0.5mm]{0pt}{10mm}\includegraphics[width=1cm]{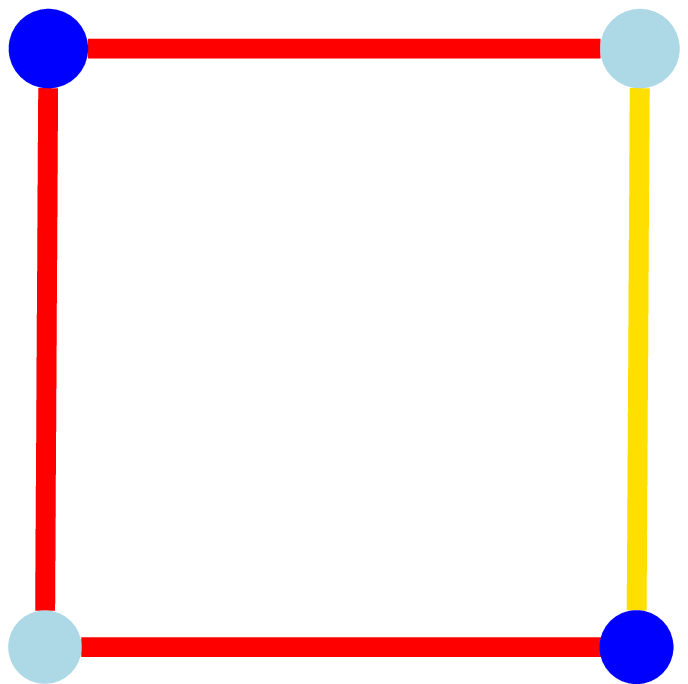}} & \rule[-7mm]{0pt}{15mm}$\begin{pmatrix} \lambda_1 & \lambda_3 & 0 & \lambda_3\\
\lambda_3 & \lambda_2 & \lambda_4 & 0\\
0 & \lambda_4 & \lambda_1 & \lambda_3\\
\lambda_3 & 0 & \lambda_3 & \lambda_2 \end{pmatrix}$ & 4 & 11 & 1:1, 2:10 & 5 & $2_2+4_3$\\
\parbox{1cm}{\rule[-0.5mm]{0pt}{10mm}\includegraphics[width=1cm]{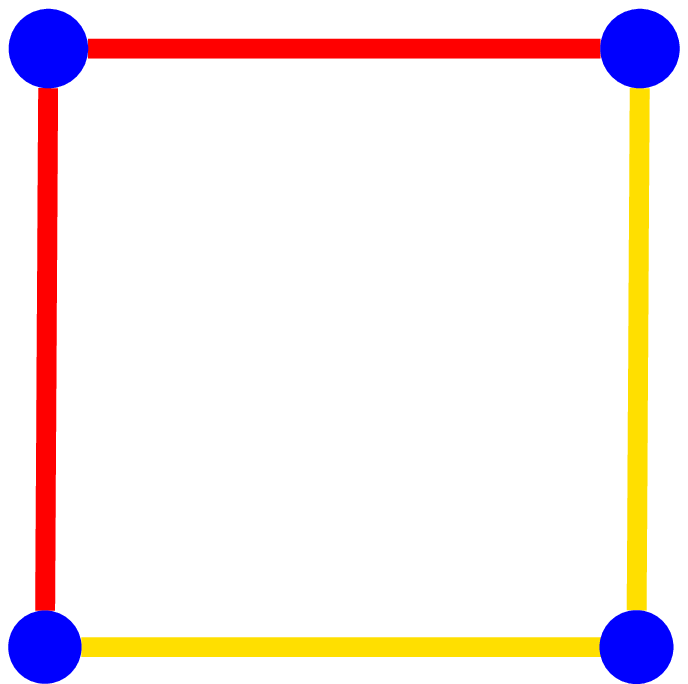}} & \rule[-7mm]{0pt}{15mm}$\begin{pmatrix} \lambda_1 & \lambda_2 & 0 & \lambda_2\\
\lambda_2 & \lambda_1 & \lambda_3 & 0\\
0 & \lambda_3 & \lambda_1 & \lambda_3\\
\lambda_2 & 0 & \lambda_3 & \lambda_1 \end{pmatrix}$ & 3 & 4 & 1:4, 2:6 & 2 & $2_3$\\
\parbox{1cm}{\rule[-0.5mm]{0pt}{10mm}\includegraphics[width=1cm]{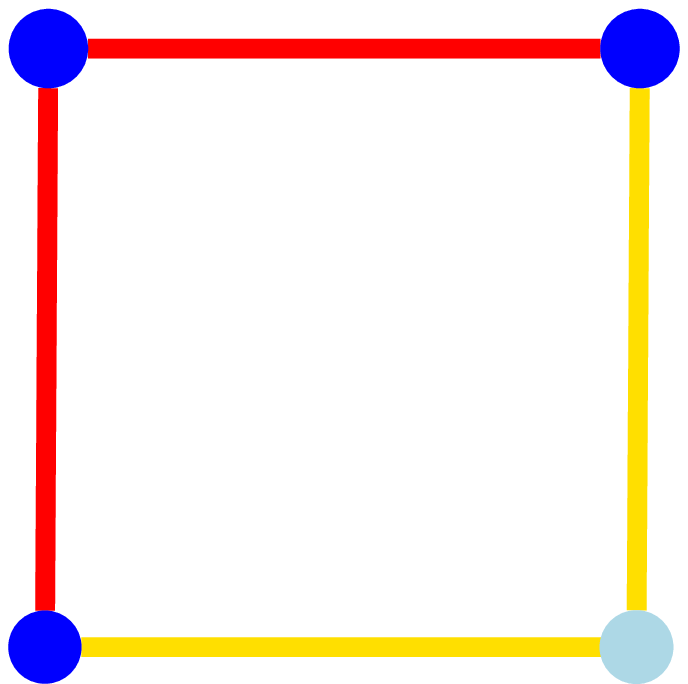}} & \rule[-7mm]{0pt}{15mm}$\begin{pmatrix} \lambda_1 & \lambda_3 & 0 & \lambda_3\\
\lambda_3 & \lambda_1 & \lambda_4 & 0\\
0 & \lambda_4 & \lambda_2 & \lambda_4\\
\lambda_3 & 0 & \lambda_4 & \lambda_1 \end{pmatrix}$ & 4 & 6 & 1:3, 2:4 & 3 & $2_2+4_3$\\
\parbox{1cm}{\includegraphics[width=1cm]{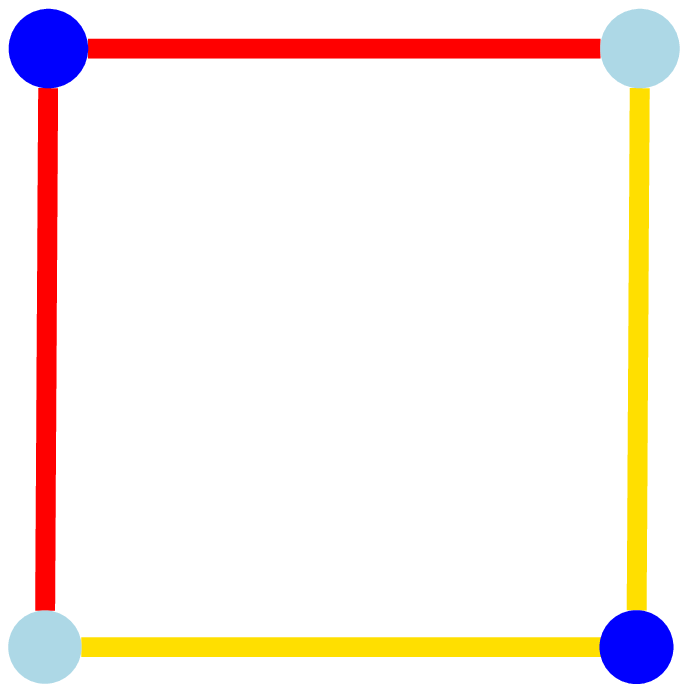}} & \rule[-7mm]{0pt}{15mm}$\begin{pmatrix} \lambda_1 & \lambda_3 & 0 & \lambda_3\\
\lambda_3 & \lambda_2 & \lambda_4 & 0\\
0 & \lambda_4 & \lambda_1 & \lambda_4\\
\lambda_3 & 0 & \lambda_4 & \lambda_2 \end{pmatrix}$ & 4 & 8 & 1:3, 2:2, 3:4 & 2 & $2_3$\\
\parbox{1cm}{\includegraphics[width=1cm]{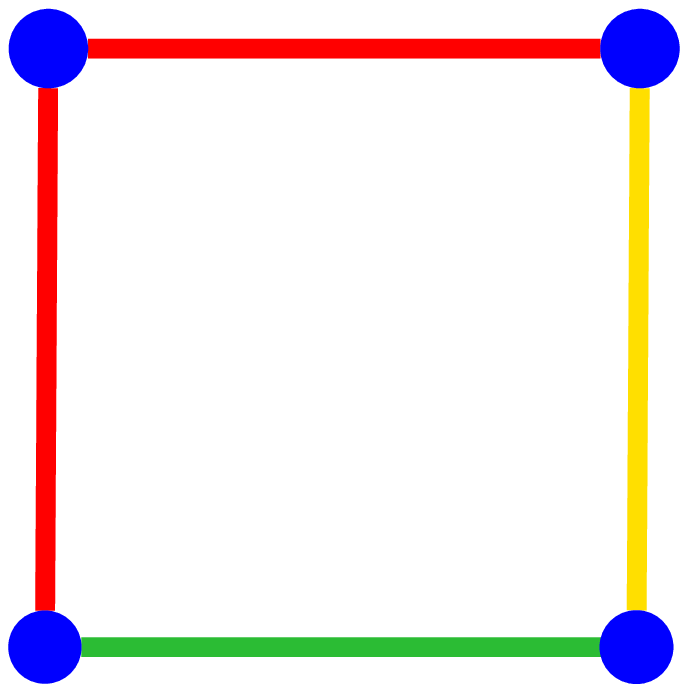}} & \rule[-7mm]{0pt}{15mm}$\begin{pmatrix} \lambda_1 & \lambda_2 & 0 & \lambda_2\\
\lambda_2 & \lambda_1 & \lambda_3 & 0\\
0 & \lambda_3 & \lambda_1 & \lambda_4\\
\lambda_2 & 0 & \lambda_4 & \lambda_1 \end{pmatrix}$ & 4 & 11 & 1:1, 2:10, 3:1 & 6 & $2_2+2_2+4_3$\\
\parbox{1cm}{\includegraphics[width=1cm]{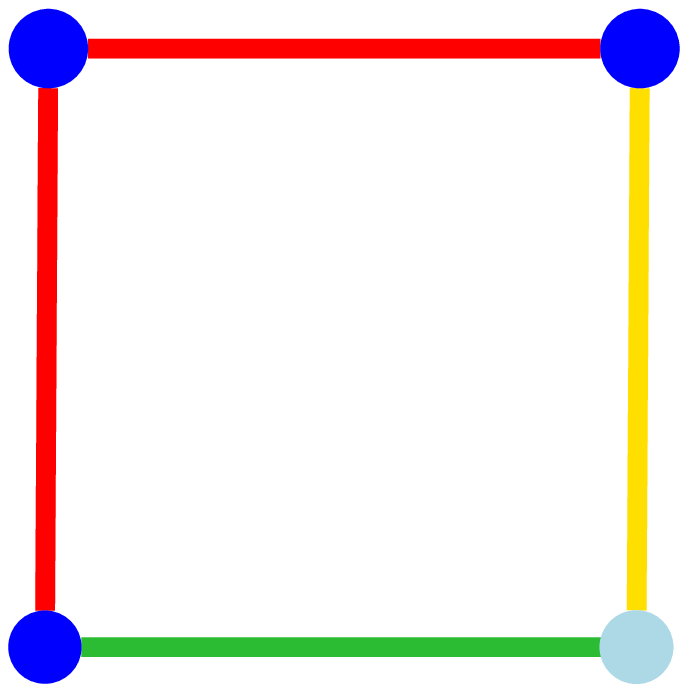}} & \rule[-7mm]{0pt}{15mm}$\begin{pmatrix} \lambda_1 & \lambda_3 & 0 & \lambda_3\\
\lambda_3 & \lambda_1 & \lambda_4 & 0\\
0 & \lambda_4 & \lambda_2 & \lambda_5\\
\lambda_3 & 0 & \lambda_5 & \lambda_1 \end{pmatrix}$ & 5 & 13 & 2:8, 3:3 & 3 & $4_2+4_3$\\
\parbox{1cm}{\includegraphics[width=1cm]{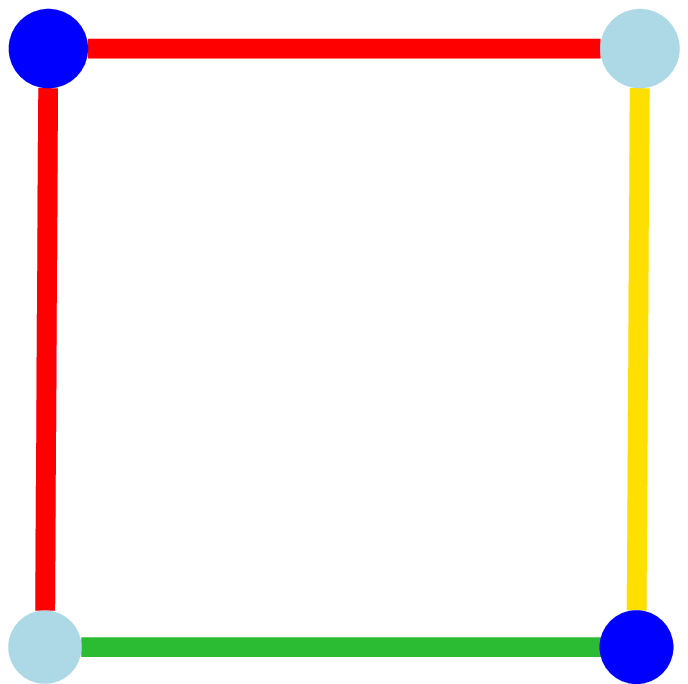}} & \rule[-7mm]{0pt}{15mm}$\begin{pmatrix} \lambda_1 & \lambda_3 & 0 & \lambda_3\\
\lambda_3 & \lambda_2 & \lambda_4 & 0\\
0 & \lambda_4 & \lambda_1 & \lambda_5\\
\lambda_3 & 0 & \lambda_5 & \lambda_2 \end{pmatrix}$ & 5 & 21 & 2:5, 3:10 & 6 & $2_2+2_2+4_3$\\
\parbox{1cm}{\includegraphics[width=1cm]{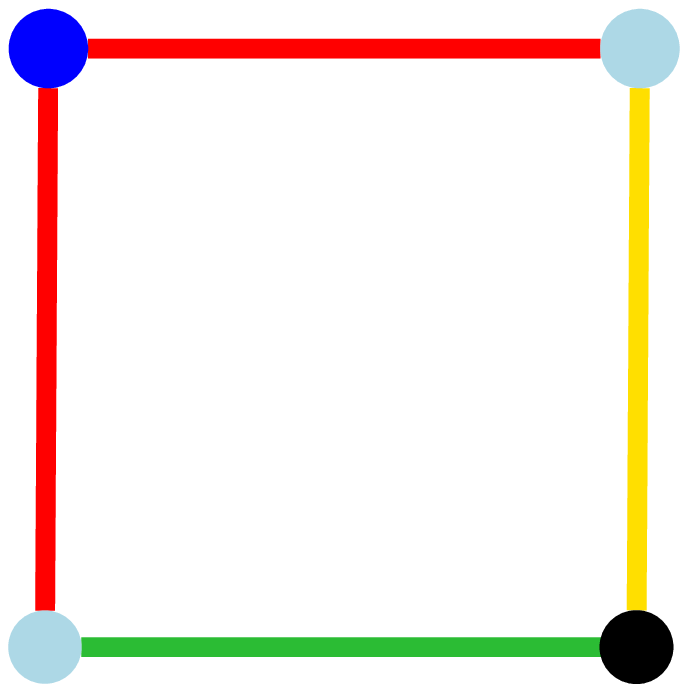}} & \rule[-7mm]{0pt}{15mm}$\begin{pmatrix} \lambda_1 & \lambda_4 & 0 & \lambda_4\\
\lambda_4 & \lambda_2 & \lambda_5 & 0\\
0 & \lambda_5 & \lambda_3 & \lambda_6\\
\lambda_4 & 0 & \lambda_6 & \lambda_2 \end{pmatrix}$ & 6 & 15 & 2:5, 3:1 & 3 & $2_2+3_2$\\
\parbox{1cm}{\includegraphics[width=1cm]{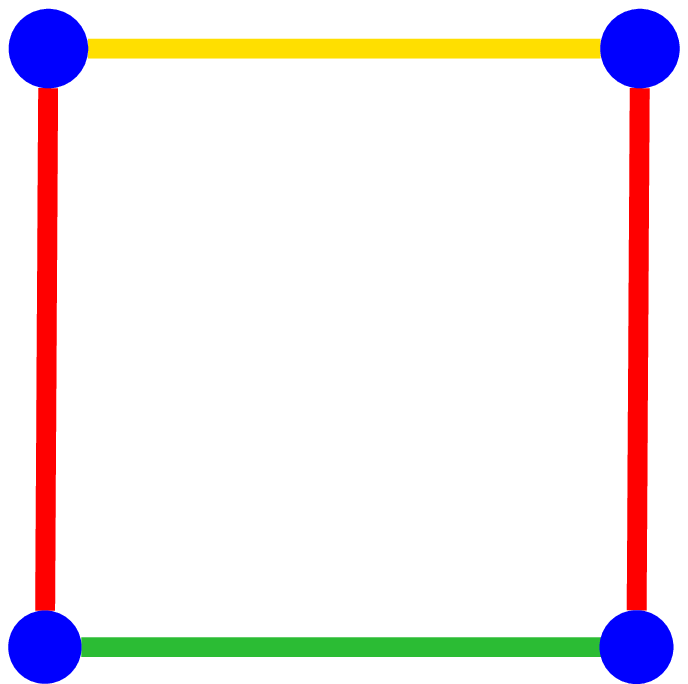}} & \rule[-7mm]{0pt}{15mm}$\begin{pmatrix} \lambda_1 & \lambda_2 & 0 & \lambda_3\\
\lambda_2 & \lambda_1 & \lambda_3 & 0\\
0 & \lambda_3 & \lambda_1 & \lambda_4\\
\lambda_3 & 0 & \lambda_4 & \lambda_1 \end{pmatrix}$ & 4 & 5 & 1:4, 2:1, 3:2 & 3 & $1_2+1_2+2_2$\\
\parbox{1cm}{\includegraphics[width=1cm]{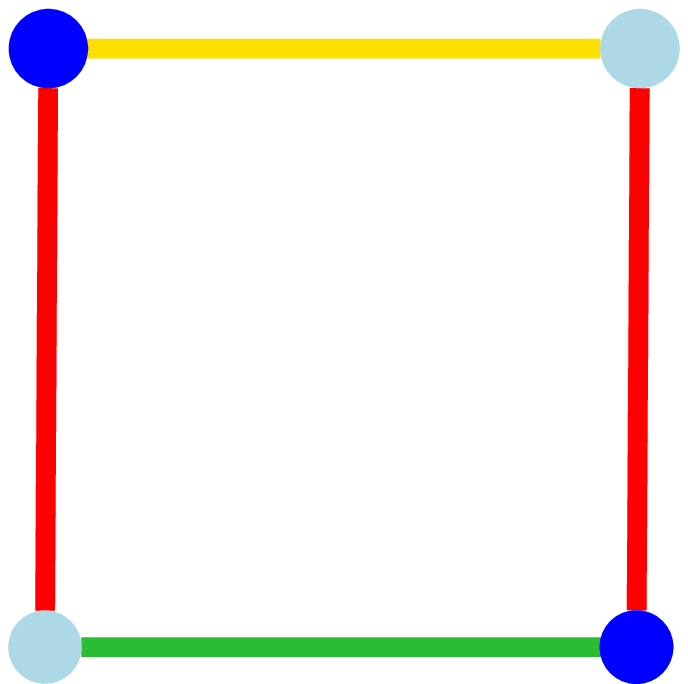}} & \rule[-7mm]{0pt}{15mm}$\begin{pmatrix} \lambda_1 & \lambda_3 & 0 & \lambda_4\\
\lambda_3 & \lambda_2 & \lambda_4 & 0\\
0 & \lambda_4 & \lambda_1 & \lambda_5\\
\lambda_4 & 0 & \lambda_5 & \lambda_2 \end{pmatrix}$ & 5 & 11 & 1:1, 2:5, 3:4 & 3 & $2_2+2_2$\\
\parbox{1cm}{\includegraphics[width=1cm]{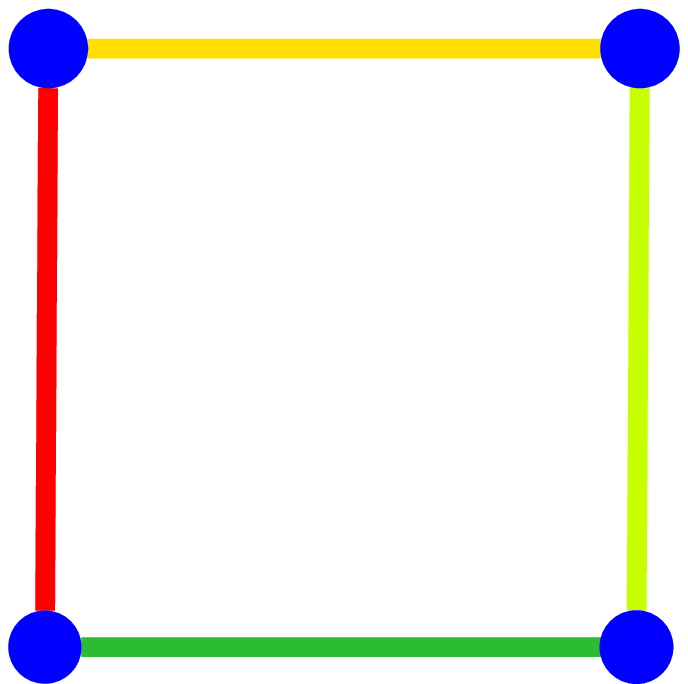}} & \rule[-7mm]{0pt}{15mm}$\begin{pmatrix} \lambda_1 & \lambda_2 & 0 & \lambda_5\\
\lambda_2 & \lambda_1 & \lambda_3 & 0\\
0 & \lambda_3 & \lambda_1 & \lambda_4\\
\lambda_5 & 0 & \lambda_4 & \lambda_1 \end{pmatrix}$ & 5 & 11 & 1:1, 2:5, 3:4 & 3 & $2_2+2_2$\\
\parbox{1cm}{\includegraphics[width=1cm]{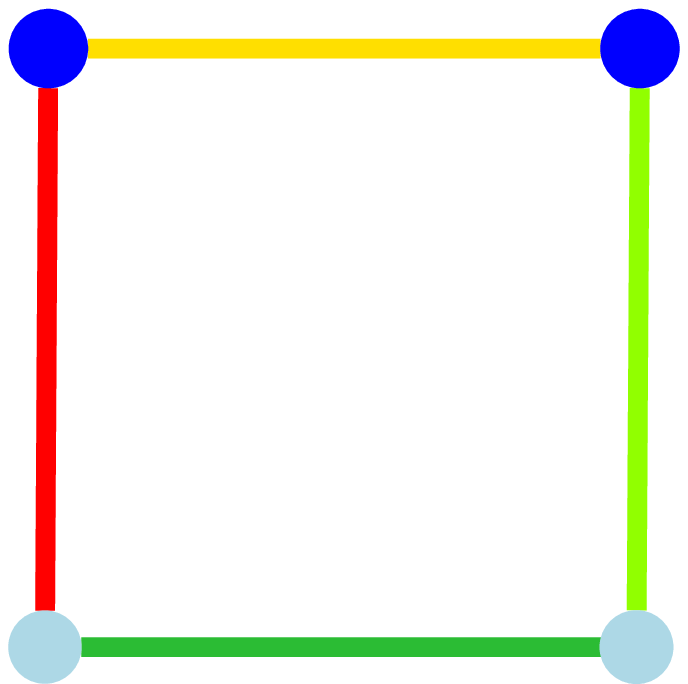}} & \rule[-7mm]{0pt}{15mm}$\begin{pmatrix} \lambda_1& \lambda_3& 0& \lambda_6\\
\lambda_3& \lambda_1& \lambda_4& 0\\
0& \lambda_4& \lambda_2& \lambda_5\\
\lambda_6& 0& \lambda_5& \lambda_2 \end{pmatrix}$ & 6 & 17 & 2:3, 3:4 & 5 & $1_1+1_1+1_1+1_1+4_2+4_2$\\
\end{tabular}
\end{table}
\begin{table}
\caption{Continuation of Table \ref{RCOR1}.}
\label{RCOR2}
\centering
\begin{tabular}{l | c | c | c | c | c | c }
\parbox{1cm}{\includegraphics[width=1cm]{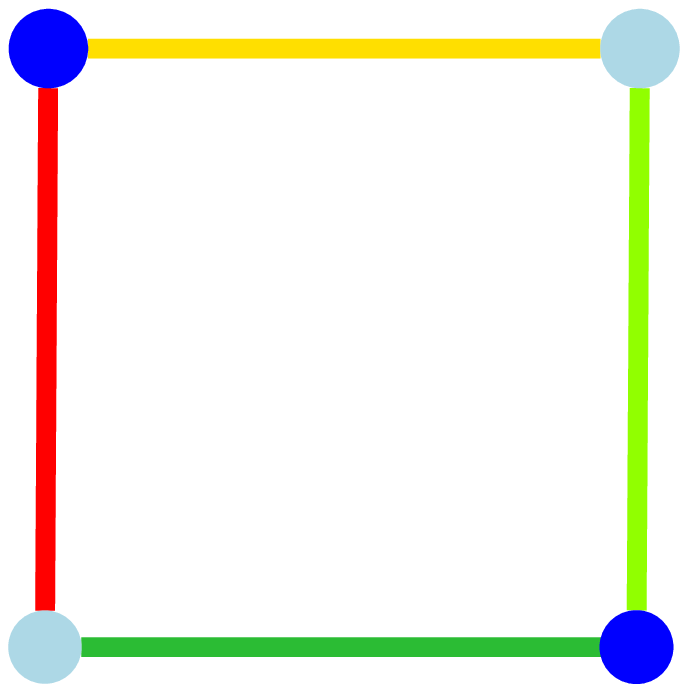}} & \rule[-7mm]{0pt}{15mm}$\begin{pmatrix} \lambda_1 & \lambda_3 & 0 & \lambda_6\\
\lambda_3 & \lambda_2 & \lambda_4 & 0\\
0 & \lambda_4 & \lambda_1 & \lambda_5\\
\lambda_6 & 0 & \lambda_5 & \lambda_2 \end{pmatrix}$ & 6 & 21 & 3:10, 4:12 & 3 & $2_2+2_2$\\
\parbox{1cm}{\includegraphics[width=1cm]{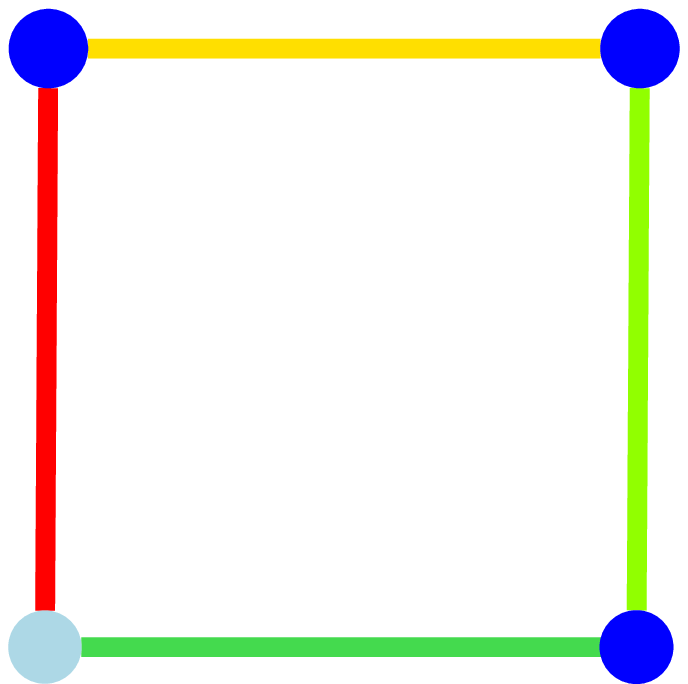}} & \rule[-7mm]{0pt}{15mm}$\begin{pmatrix} \lambda_1 & \lambda_3 & 0 & \lambda_6\\
\lambda_3 &\lambda_1 &\lambda_4 &0\\
0 &\lambda_4 &\lambda_1 &\lambda_5\\
\lambda_6 &0 &\lambda_5 &\lambda_2 \end{pmatrix}$ & \quad 6 \qquad & \quad 17 \qquad & 2:2, 3:8, 4:1 & \quad 4 \qquad & \quad $10_2$ \qquad\\
\parbox{1cm}{\includegraphics[width=1cm]{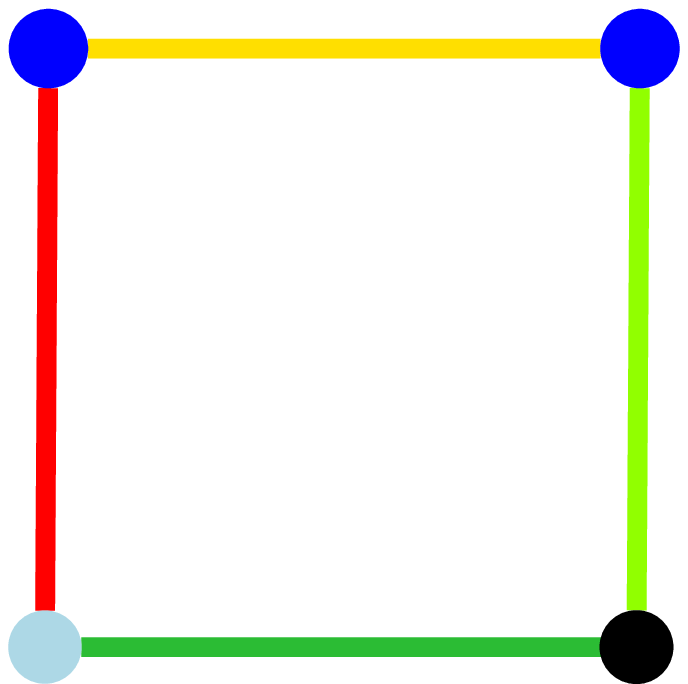}} & \rule[-7mm]{0pt}{15mm}$\begin{pmatrix} \lambda_1 &\lambda_4& 0& \lambda_7\\
\lambda_4 &\lambda_1 &\lambda_5 &0\\
0 &\lambda_5 &\lambda_2 &\lambda_6\\
\lambda_7 &0& \lambda_6 &\lambda_3 \end{pmatrix}$ & 7 & 13 & 2:1, 3:3 & 5 & $1_1+1_1+2_1+12_2$\\
\parbox{1cm}{\includegraphics[width=1cm]{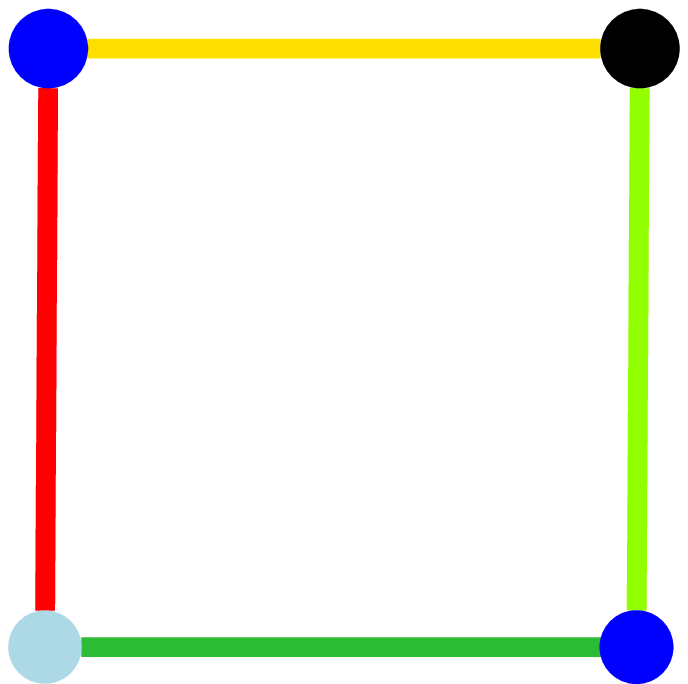}} & \rule[-7mm]{0pt}{15mm}$\begin{pmatrix} \lambda_1&\lambda_4& 0& \lambda_7\\
\lambda_4 &\lambda_2 &\lambda_5 &0\\
0 &\lambda_5 &\lambda_1 &\lambda_6\\
\lambda_7 &0& \lambda_6 &\lambda_3 \end{pmatrix}$ & 7 & 17 & 3:3, 4:6 & 3 & $4_2$\\
\parbox{1cm}{\includegraphics[width=1cm]{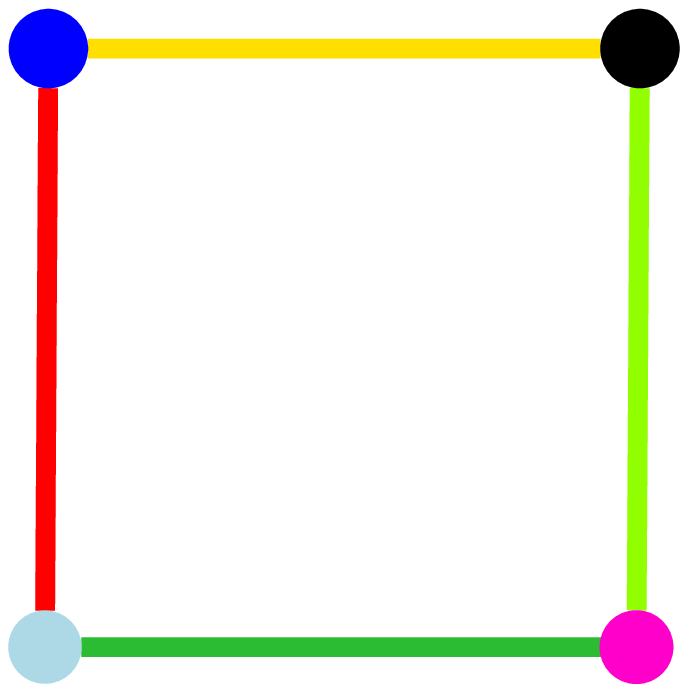}} & \rule[-7mm]{0pt}{15mm}$\begin{pmatrix} \lambda_1& \lambda_5& 0& \lambda_8\\
\lambda_5 &\lambda_2 &\lambda_6 &0\\
0 &\lambda_6 &\lambda_3 &\lambda_7\\
\lambda_8 &0 &\lambda_7 &\lambda_4 \end{pmatrix}$ & 8 & 9 & 3:2 & 5 & $2_1+2_1+2_1+2_1+8_2$
\end{tabular}
\end{table}

\vskip -0.5cm

\begin{table}
\caption{All RCOP-models \citep{HL}
when the underlying graph is the $4$-cycle.}
\label{RCOP}
\centering
\begin{tabular}{l | c | c | c | c | c| c}
Graph & $K$ &dim $d$ & degree & mingens $P_\mathcal{G}$
 & ML-degree & $H_\mathcal{L}$  \\ \hline
\parbox{1cm}{\includegraphics[width=1cm]{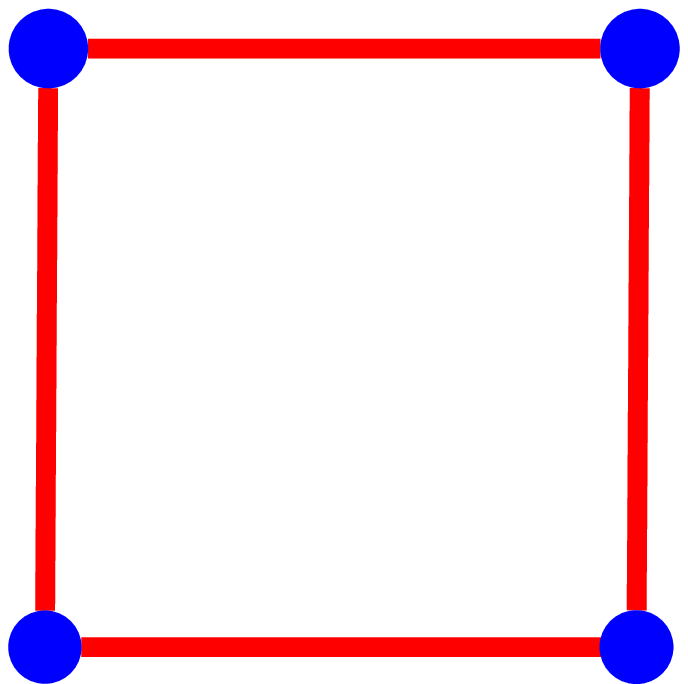}} & \rule[-7mm]{0pt}{15mm}$\begin{pmatrix} \lambda_1 & \lambda_2 & 0 & \lambda_2 \\
\lambda_2 & \lambda_1 & \lambda_2 & 0 \\
0 & \lambda_2 & \lambda_1 & \lambda_2 \\
\lambda_2 & 0 & \lambda_2 & \lambda_1 \end{pmatrix}$ & 2 & 2 & 1:7, 2:1 &2& $(2t_1-t_2)(2t_1+t_2)$  \\
\parbox{1cm}{\includegraphics[width=1cm]{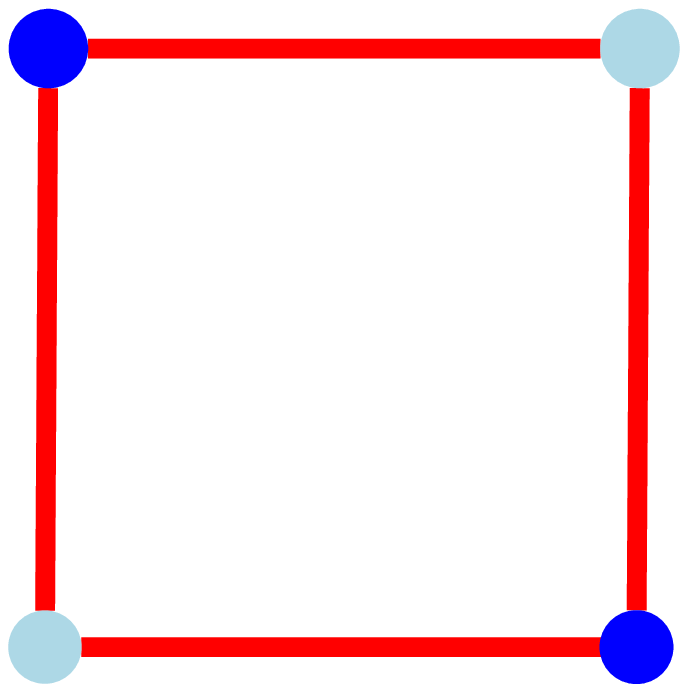}} & \rule[-7mm]{0pt}{15mm}$\begin{pmatrix} \lambda_1 & \lambda_3 & 0 & \lambda_3\\
\lambda_3 & \lambda_2 & \lambda_3 & 0\\
0 & \lambda_3 & \lambda_1 & \lambda_3\\
\lambda_3 & 0 & \lambda_3 & \lambda_2 \end{pmatrix}$ & 3 & 4 & 1:5, 2:2 & 2 & $16t_1t_2-t_3^2$  \\
\parbox{1cm}{\includegraphics[width=1cm]{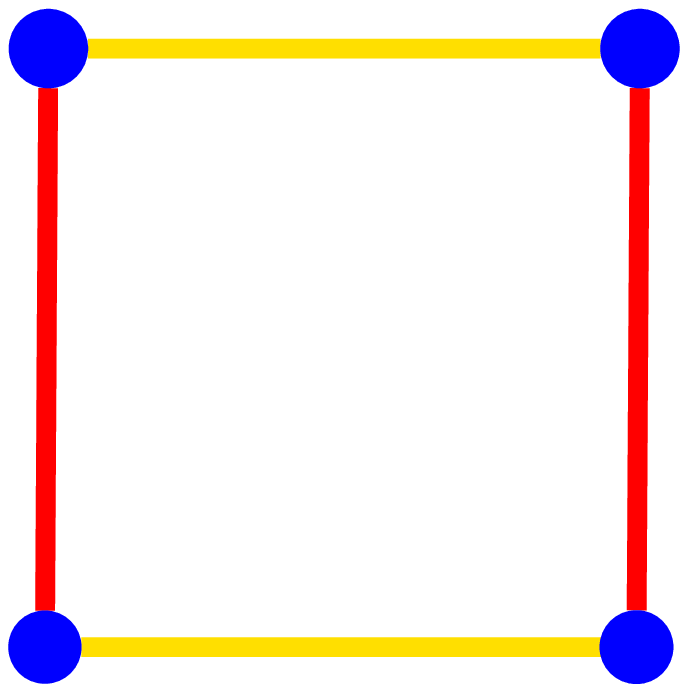}} & \rule[-7mm]{0pt}{15mm}$\begin{pmatrix} \lambda_1 & \lambda_2 & 0 & \lambda_3\\
\lambda_2 & \lambda_1 & \lambda_3 & 0\\
0 & \lambda_3 & \lambda_1 & \lambda_2\\
\lambda_3 & 0 & \lambda_2 & \lambda_1 \end{pmatrix}$ & 3 & 3 & 1:6, 3:1 & 3 & $(t_1-t_2)(t_1+t_2)(t_1-t_3)(t_1+t_3)$  \\
\parbox{1cm}{\includegraphics[width=1cm]{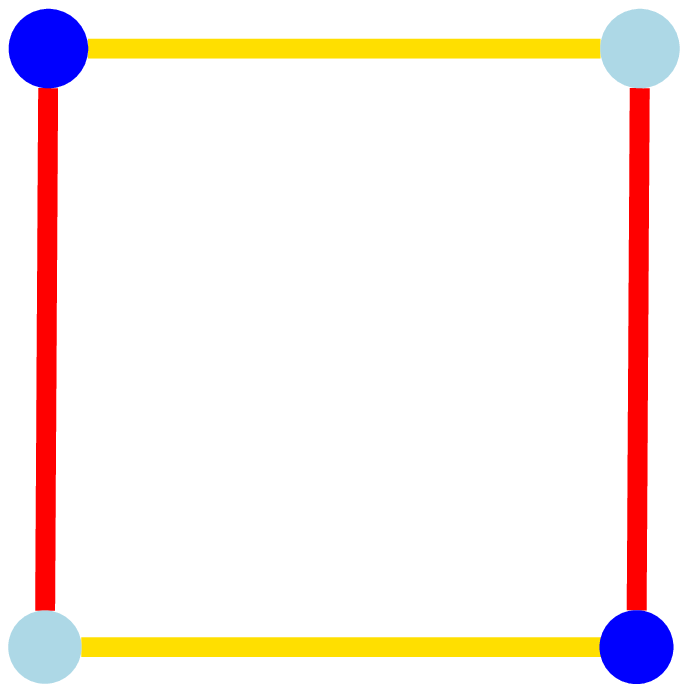}} & \rule[-7mm]{0pt}{15mm}$\begin{pmatrix} \lambda_1 & \lambda_3 & 0 & \lambda_4\\
\lambda_3 & \lambda_2 & \lambda_4 & 0\\
0 & \lambda_4 & \lambda_1 & \lambda_3\\
\lambda_4 & 0 & \lambda_3 & \lambda_2 \end{pmatrix}$ & 4 & 5 & 1:4, 2:1, 3:2 & 3 & $(4t_1t_2-t_3^2)(4t_1t_2-t_4^2)$  \\
\parbox{1cm}{\includegraphics[width=1cm]{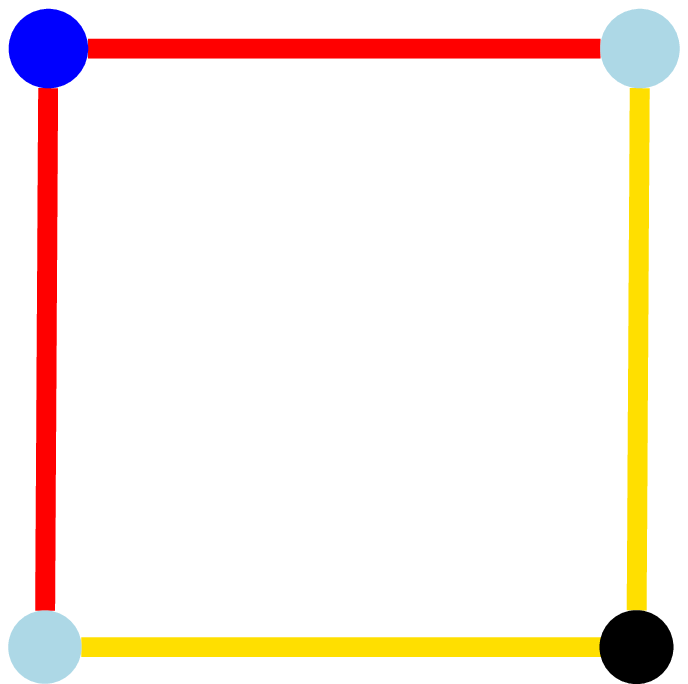}} & \rule[-7mm]{0pt}{15mm}$\begin{pmatrix} \lambda_1 & \lambda_4 & 0 & \lambda_4\\
\lambda_4 & \lambda_2 & \lambda_5 & 0\\
0 & \lambda_5 & \lambda_3 & \lambda_5\\
\lambda_4 & 0 & \lambda_5 & \lambda_2 \end{pmatrix}$ & 5 & 6 & 1:3, 2:1, 3:1 & 3 & $(8t_1t_2-t_4^2)(8t_2t_3-t_5^2)$ \\
\parbox{1cm}{\includegraphics[width=1cm]{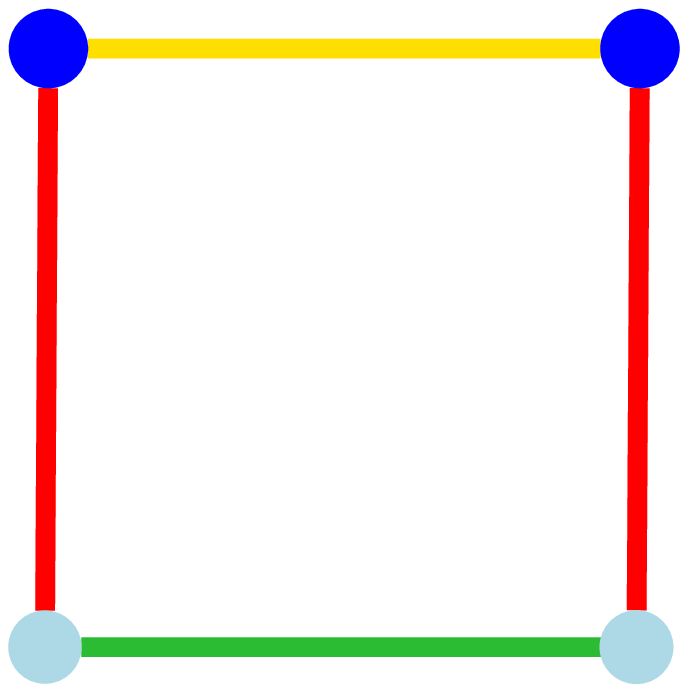}} & \rule[-7mm]{0pt}{15mm}$\begin{pmatrix} \lambda_1 & \lambda_3 & 0 & \lambda_4\\
\lambda_3 & \lambda_1 & \lambda_4 & 0\\
0 & \lambda_4 & \lambda_2 & \lambda_5\\
\lambda_4 & 0 & \lambda_5 & \lambda_2 \end{pmatrix}$ & 5 & 3 & 1:4, 3:1 & 3 & (\ref{eq_ex}) in Example \ref{ex_Fret}
\end{tabular}
\end{table}

\begin{example}
We can gain a different perspective on the proof of Lemma \ref{MLispositive}
by considering colored Gaussian graphical models.
Under the assumption (\ref{allequal}) that all parameters in the partial
matrix (\ref{partialmatrix}) are equal to some fixed value $x$,
the MLE $\hat K$ for the concentration matrix has the same structure.
Namely, all diagonal entries of $\hat K$ are equal,
and all non-zero off-diagonal entries of $\hat K$ are equal.
This means that we can perform our MLE computation for the colored Gaussian graphical model with the chordless $m$-cycle as
underlying graph, where all vertices and all edges have the same color:
\begin{equation}
\label{Kl1l2}
K \,\, =\,\,
\begin{pmatrix}
    \lambda_1 & \lambda_2 &  0   & 0     &  \cdots & \lambda_2 \\
\lambda_2 & \lambda_1    & \lambda_2 & 0     &  \cdots &   0 \\
 0   & \lambda_2 & \lambda_1    & \lambda_2 &   \cdots  & 0 \\
 \vdots & \vdots & \ddots & \ddots & \ddots &  0  \\
   0        & 0 &     0   & \lambda_2 & \lambda_1 & \lambda_2 \\
   \lambda_2  & 0  &    0  &   0    & \lambda_2 & \lambda_1
  \end{pmatrix}.
 \end{equation}
In contrast to the approach in the proof of Lemma \ref{MLispositive},
in this representation we only need to solve a system of
two polynomial equations in two unknowns,
regardless of the cycle size $m$.  The equations are
$$ (K^{-1})_{11} = 1 \quad \hbox{and} \quad (K^{-1})_{12} = x . $$
By clearing denominators we obtain two polynomial equations
in the unknowns $\lambda_1$ and $\lambda_2$.
We need to express these in terms of the parameter $x$,
but there are many extraneous solutions.
The ML degree is algebraic degree of the special
solution $\, (\hat \lambda_1 (x),\hat \lambda_2 (x))\,$
 which makes (\ref{Kl1l2}) positive definite. \qed
\end{example}

\begin{figure*}[!h]
\centering
\includegraphics[width=5.8cm]{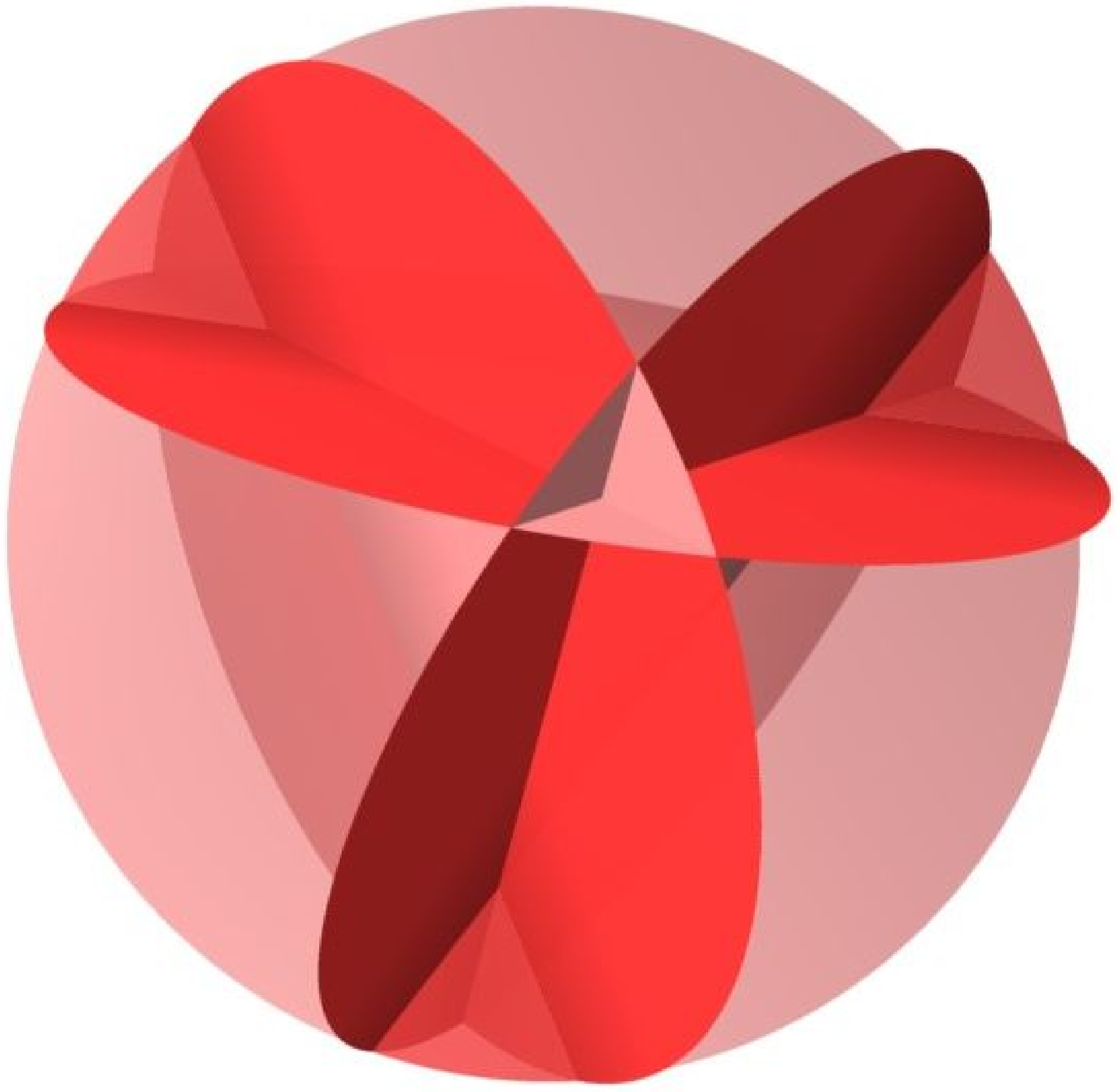}
\quad\qquad\qquad \includegraphics[width=5.9cm]{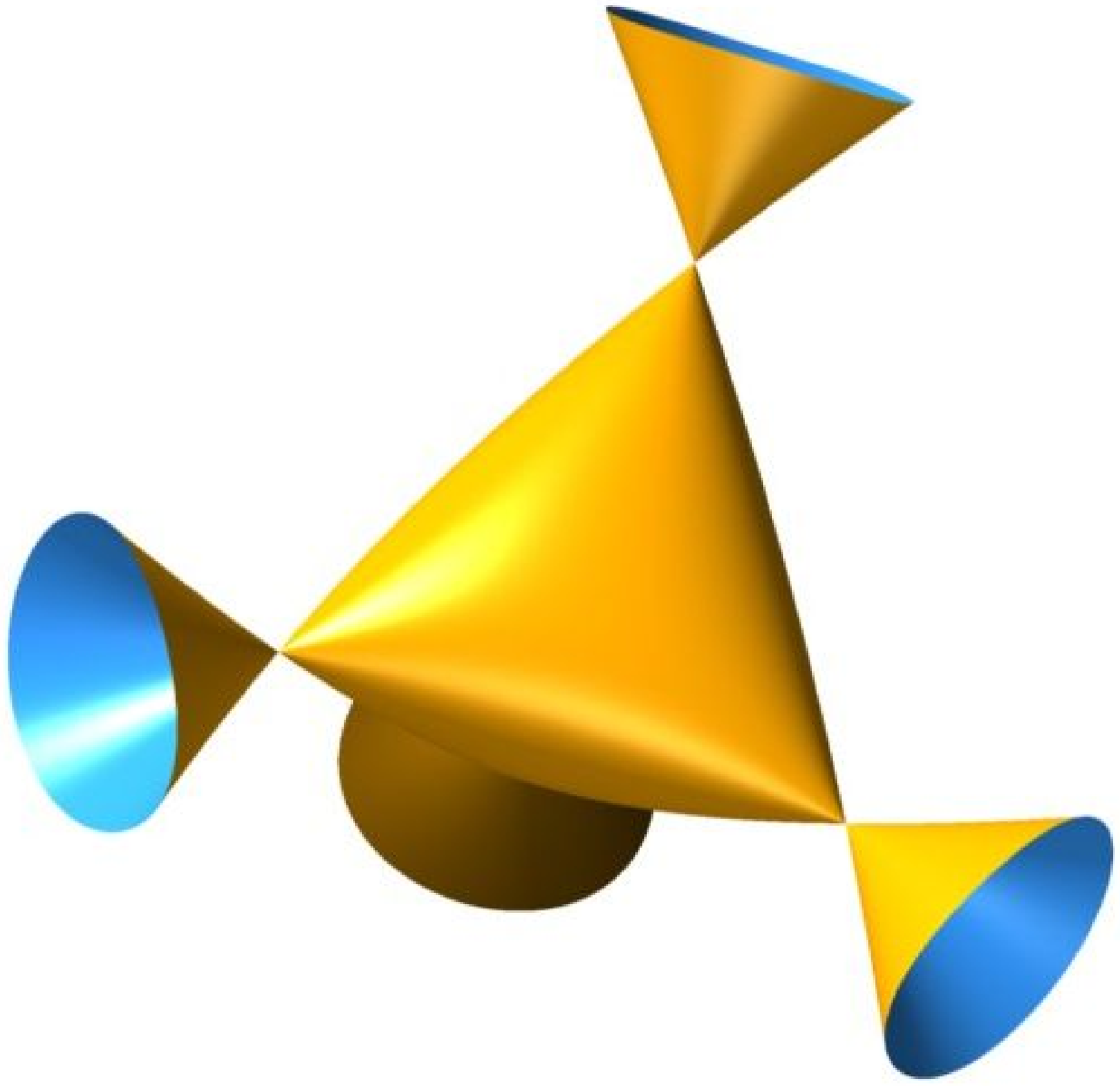} \\
\includegraphics[width=5.7cm]{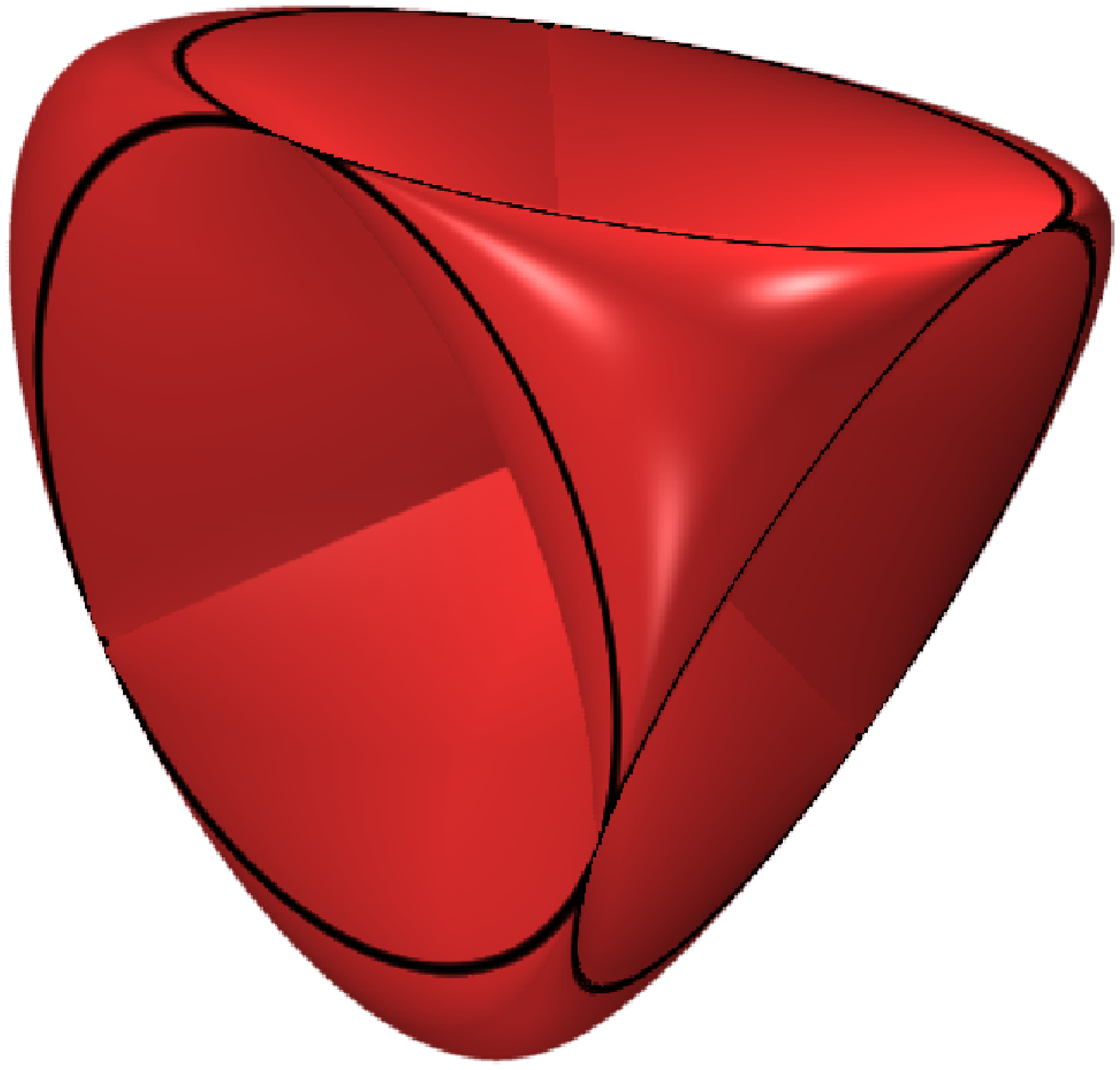}
\quad\qquad\qquad \includegraphics[width=5.5cm]{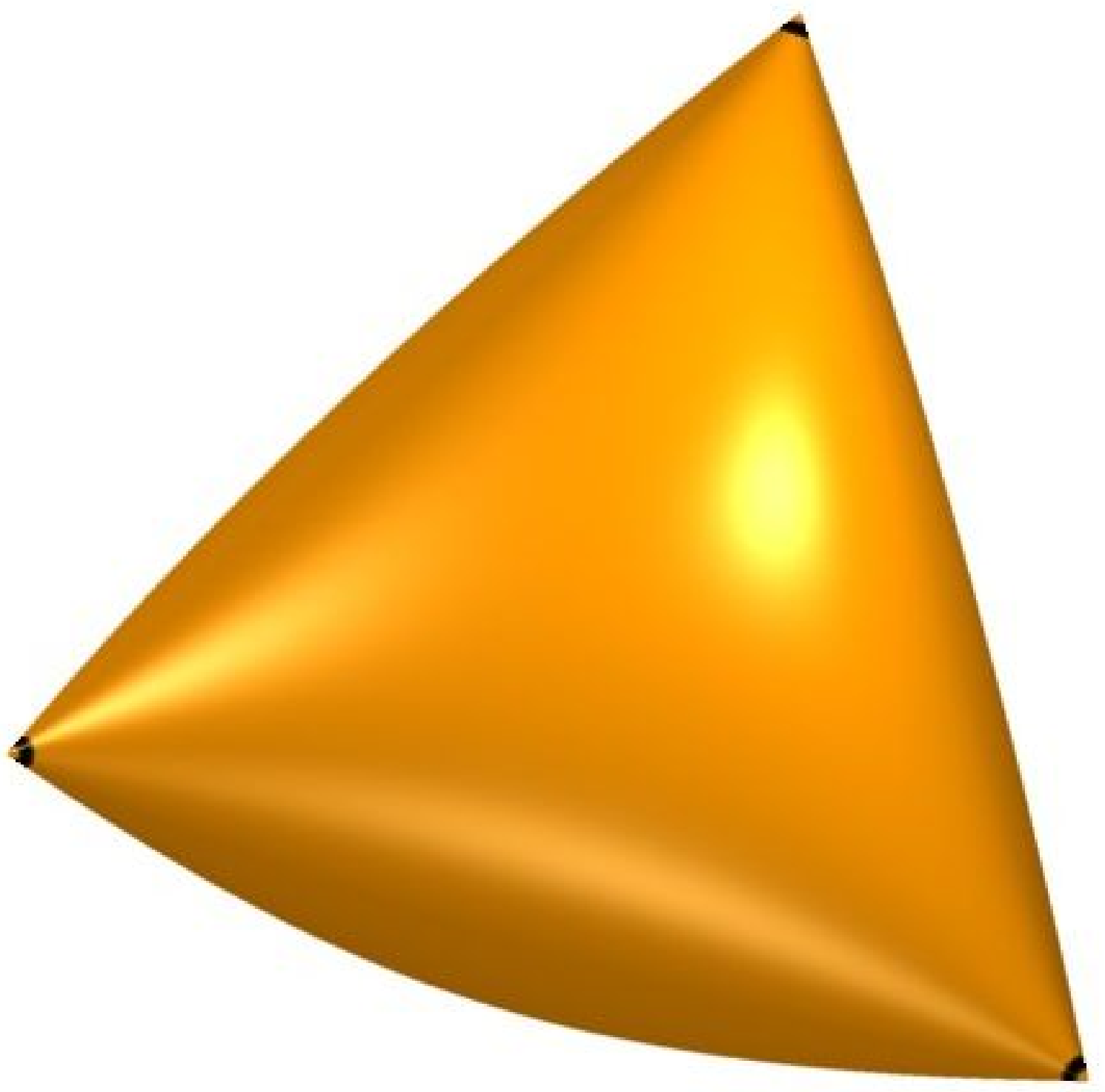} \\
\caption{The cross section of the cone of sufficient statistics
in Example \ref{ex:dualpillow} is the red convex body shown in the left figure.
It is dual to Cayley's cubic surface, which is shown in
yellow in the right figure
and also in Fig.~\ref{fig_labs} on the left.}
\label{dualpillow}
\end{figure*}

\begin{example}
\label{ex:dualpillow}
Let $\mathcal{G}$ be the colored triangle with the
same color for all three vertices and three distinct colors
for the edges. This is an RCOP model with $m=3$ and $d=4$.
The corresponding subspace $\mathcal{L}$ of $\s^3$ consists of
all concentration matrices
$$ K \quad = \quad
\begin{pmatrix}
\lambda_4 & \lambda_1 &  \lambda_2 \\
\lambda_1 & \lambda_4 &  \lambda_3 \\
\lambda_2 & \lambda_3 &  \lambda_4
\end{pmatrix}.
$$
This linear space $\mathcal{L}$ is
generic enough so as to exhibit the
geometric behavior described in Subsection 2.2.
The four-dimensional cone  $\mathcal{K}_{\mathcal{L}}$
is the cone over the $3$-dimensional spectrahedron  bounded by
{\em Cayley's cubic surface} as shown on
the right in Fig.~\ref{dualpillow}. Its dual
$\mathcal{C}_{\mathcal{L}}$ is the cone over the $3$-dimensional
convex body shown on the left in Fig.~\ref{dualpillow}.
The boundary of this convex body consists of four flat $2$-dimensional
circular faces (shown in black)
and four curved surfaces whose common Zariski closure is a
{\em quartic Steiner surface}. Fig.~\ref{dualpillow}
was made with {\tt surfex}\footnote{www.surfex.algebraicsurface.net/}, a
software package for visualizing algebraic surfaces.

Here, the
inequalities (\ref{pataki}) state $2 \leq p \leq 3$, and the
algebraic degree of SDP is
$\,\delta(3,3,2) = \delta(3,3,1) = 4$.
We find that $H_\mathcal{L}$
is a polynomial of degree $8$ which factors  into
four linear forms and one quartic:
$$ H_\mathcal{L} =
(t_1-t_2+t_3-t_4)(t_1+t_2-t_3-t_4) (t_1-t_2-t_3+t_4) (t_1+t_2+t_3+t_4)
(t_1^2 t_2^2+t_1^2 t_3^2+t_2^2 t_3^2-2 t_1 t_2 t_3 t_4) $$
By Theorem \ref{deg=MLdeg},
 both the degree and the ML degree of this
 model are also equal to  $\,\phi(3,4) = 4$. \qed
\end{example}

\begin{acknowledgements}
We wish to thank Seth Sullivant,  Bernd Ulrich
and Ruriko Yoshida for helpful comments.
\end{acknowledgements}



\end{document}